\documentclass{amsart}

\usepackage{amsmath,amssymb}
\usepackage{tikz}
\usepackage[unicode=true]{hyperref}

\usetikzlibrary{arrows,intersections}

\pgfarrowsdeclare{acts}{acts}{
\pgfarrowsleftextend{0pt}
\pgfarrowsrightextend{0pt}
}{
\pgfsetdash{}{0pt} 
\pgfsetroundjoin 
\pgfsetroundcap 
\pgfpathmoveto{\pgfpoint{0pt}{0pt}}
\pgfpatharc{90}{-180}{1.5pt}
\pgfpathlineto{\pgfpoint{-1.5pt}{3pt}}
\pgfusepathqstroke
}

\newtheorem{thm}{Theorem}[section]
\newtheorem{cor}[thm]{Corollary}
\newtheorem{prop}[thm]{Proposition}
\newtheorem{lem}[thm]{Lemma}

\theoremstyle{definition}
\newtheorem{defn}[thm]{Definition}
\newtheorem{notn}[thm]{Notation}
\newtheorem{conv}[thm]{Convention}
\newtheorem{hyp}[thm]{Hypothesis}

\newtheorem{rem}[thm]{Remark}
\newtheorem{ex}[thm]{Example}

\newcommand{\bbA}{\mathbb{A}}
\newcommand{\bbC}{\mathbb{C}}
\newcommand{\bbF}{\mathbb{F}}
\newcommand{\grass}{\mathbb{G}}
\newcommand{\bbN}{\mathbb{N}}
\newcommand{\bbP}{\mathbb{P}}
\newcommand{\bbQ}{\mathbb{Q}}
\newcommand{\bbR}{\mathbb{R}}

\newcommand{\bbZ}{\mathbb{Z}}

\newcommand{\calA}{\mathcal{A}}

\newcommand{\pibar}{\overline{\pi}}
\newcommand{\picomp}{\pi^{\vee}}
\newcommand{\rhobar}{\overline{\rho}}

\newcommand{\del}{\partial}
\newcommand{\surject}{\twoheadrightarrow}
\newcommand{\inject}{\hookrightarrow}
\newcommand{\iso}{\mathrel{\overset{\sim}{\smash{\longrightarrow}\vrule height.3ex
width0ex\relax}}}
\newcommand{\dcup}{\mathrel{\dot{\cup}}}
\newcommand{\bigdcup}{\mathop{\dot{\bigcup}}}

\newcommand{\mult}{^{\times}}

\newcommand{\eq}{^{\mathrm{eq}}}
\newcommand{\auf}[1]{\mathord{|_{#1}}}

\newcommand{\topcl}[1]{\overline{#1}^{\mathrm{top}}}

\newcommand{\1}{^{-1}}

\newcommand{\valring}{\mathcal{O}_{\!K}}
\newcommand{\maxid}{\mathcal{M}_K}

\newcommand{\ball}[1]{B(#1)}

\newcommand{\CN}[1][N]{\mathcal{C}_{#1}} 

\newcommand{\Lx}{\mathcal{L}}
\newcommand{\LHen}{\Lx_{\mathrm{Hen}}}
\newcommand{\LHenA}{\Lx_{\mathrm{Hen}, \calA}}
\newcommand{\Lring}{\Lx_{\mathrm{ring}}}
\newcommand{\Labsring}{\Lx_{\mathrm{absring}}}
\newcommand{\Loag}{\mathcal{L}_{\mathrm{oag}}}
\newcommand{\Loagpar}{\Loag(\mathrm{par})}

\newcommand{\Tx}{\mathcal{T}}
\newcommand{\THen}{\Tx_{\mathrm{Hen}}}
\newcommand{\THenA}{\Tx_{\mathrm{Hen}, \calA}}

\newcommand{\Zp}{\bbZ_p}

\newcommand{\Kacl}{\tilde{K}}

\newcommand{\RVacl}{\tilde{\RV}}

\newcommand{\anfct}{\mathbf{A}}

\newcommand{\aci}{\operatorname{\overline{ac}}}
\newcommand{\acl}{\operatorname{acl}}
\newcommand{\Aut}{\operatorname{Aut}}
\newcommand{\dcl}{\operatorname{dcl}}

\newcommand{\dir}{\operatorname{dir}}
\newcommand{\dirRV}{\dir_{\RV}}
\newcommand{\affdir}{\operatorname{affdir}}

\newcommand{\GL}{\operatorname{GL}}

\newcommand{\id}{\operatorname{id}}
\newcommand{\im}{\operatorname{im}}
\newcommand{\Jac}{\operatorname{Jac}}

\newcommand{\rado}{\operatorname{rad}_{\mathrm{o}}}
\newcommand{\radc}{\operatorname{rad}_{\mathrm{c}}}
\newcommand{\res}{\operatorname{res}}
\newcommand{\RV}{\operatorname{RV}}
\newcommand{\rv}{\operatorname{\hat{rv}}}
\newcommand{\rvi}{\operatorname{rv}}

\newcommand{\Tr}{\operatorname{T}}

\newcommand{\tsp}{\operatorname{tsp}}  
\newcommand{\val}{\operatorname{\hat{v}}}
\newcommand{\vali}{\operatorname{v}}
\newcommand{\vRV}{\val_{\RV}}
\newcommand{\viRV}{\vali_{\RV}}

\newbox\nsbox
\newcommand{\ns}[1]{%
  \setbox\nsbox=\hbox{\ensuremath{#1}}%
  {\vrule width0mm height\ht\nsbox}^*\!%
  \box\nsbox%
}
\newcommand{\nss}[1]{%
  \setbox\nsbox=\hbox{\ensuremath{\scriptstyle{#1}}}%
  {\vrule width0mm height\ht\nsbox}^*\!%
  \box\nsbox%
}
\newcommand{\nsN}{\ns\bbN}
\newcommand{\nsQ}{\ns\bbQ}
\newcommand{\nsR}{\ns\bbR}
\newcommand{\nsu}{\ns u}

\newcommand{\nsS}{\ns {\!S}}
\newcommand{\nsX}{\ns {\!X}}
\newcommand{\nsY}{\ns {\!Y}}

\newbox\gnBoxA
\newdimen\gnCornerHgt
\setbox\gnBoxA=\hbox{$\ulcorner$}
\global\gnCornerHgt=\ht\gnBoxA
\newdimen\gnArgHgt

\def\code #1{%
	\setbox\gnBoxA=\hbox{$#1$}%
	\gnArgHgt=\ht\gnBoxA%
	\ifnum \gnArgHgt<\gnCornerHgt
		\gnArgHgt=0pt%
	\else
		\advance \gnArgHgt by -\gnCornerHgt%
	\fi
	\raise\gnArgHgt\hbox{$\ulcorner$} \box\gnBoxA %
		\raise\gnArgHgt\hbox{$\urcorner$}}

\newcommand{\labelledclaim}[2]{%
\begin{itemize}
 \item[(#1)] #2
\end{itemize}
}

\begin{document}

\title{Non-Archimedean Whitney stratifications}
\author{Immanuel Halupczok}
\thanks{The author was partially supported by the Fondation Sciences Math\'ematiques
de Paris and partly by the SFB~478 and the SFB~878 of the Deutsche Forschungsgemeinschaft}
\address{Institut f\"ur Mathematische Logik und Grundlagenforschung
Universit\"at M\"unster\\
Einsteinstra\ss{}e 62\\
48149 M\"unster\\
Germany}
\email{math@karimmi.de}
\urladdr{http://www.immi.karimmi.de/en.math.html}

\begin{abstract}
We define ``t-stratifications'', a strong notion of stratifications for Henselian valued fields $K$ of equi-characteristic $0$, and prove that they exist. In contrast to classical stratifications in Archimedean fields,
t-stratifications also contain non-local information about the stratified sets. For example,
they do not only see the singularities in the valued field, but also those in the residue field.

Like Whitney stratifications, t-stratifications exist for different classes of subsets of $K^n$,
e.g.\ algebraic subvarieties or certain classes of analytic subsets. The general framework
are definable sets (in the sense of model theory) in a language that satisfies certain hypotheses.

We give two applications. First, we show that t-stratifications in suitable valued fields $K$ induce classical Whitney stratifications in $\bbR$ or $\bbC$; in particular, the existence of t-stratifications implies the existence of Whitney stratifications. This uses methods of non-standard analysis.

Second, we show how, using t-stratifications, one can determine the ultra-metric isometry type of definable subsets of $\bbZ_p^n$ for $p$ sufficiently big. For those $p$, this proves a conjecture stated in a previous article. In particular, this yields a new, geometric proof of the rationality of Poincar\'e series.
\end{abstract}


\maketitle

\section{Introduction}

Over the fields $\bbR$ and $\bbC$, a very useful tool to describe singularities of algebraic or analytic sets are Whitney stratifications; see e.g.\ \cite{Whi.strat}, \cite{BCR.realGeom}.
The definition of a Whitney stratification can be translated in a straightforward way to
non-Archimedean local fields and in \cite{CCL.cones}, it has been proven that in this sense, Whitney stratifications also exist in the $p$-adics $\bbQ_p$. In the present article, we introduce t-stratifications:
another kind of stratifications which exist in Henselian valued fields and whose regularity condition is,
in a certain sense, much stronger and in particular strictly non-local.
To get a first impression, suppose that we have a set $X \subseteq K^n$ for some valued field $K$.
Any ball $B \subseteq K$ around $0$ is an additive subgroup of $K$ and we can consider the image $X_B$ of
$X$ in the quotient $(K/B)^n$. A t-stratification for $X$ simultaneously describes the singularities of 
all those images $X_B$ (in some sense which we will make precise).

If we take $B$ to be the maximal ideal of the valuation ring, then the residue field $k$ of $K$ is a subset of
$K/B$, so in particular, if $k \subseteq K$ and $X = V(K) \subseteq K^n$ is a variety defined over $k$,
then a t-stratification of $X$ induces a
stratification of $V(k) = X_B \cap k^n$.
If $k$ is $\bbR$ or $\bbC$, and under some additional assumptions about $K$,
one can show that this stratification of $V(k)$ is a Whitney stratification in the classical sense.
In this way, although t-stratifications and classical Whitney stratifications live in different worlds,
we will prove that the existence of t-stratifications genuinely implies the existence of 
Whitney stratifications.

A different point of view is that the existence of a t-stratification of a set $X \subseteq K^n$ is a very strong statement about the isometry type of $X$ (with respect to the ultrametric induced by the valuation on $K$).
This leads to the original motivation for this article. The main conjecture of \cite{i.QpB} is a
complete classification of algebraic (or, more precisely, definable) sets $X \subseteq \bbZ_p^n$
up to isometry. The present results yield a proof of this conjecture for big $p$ (depending on $X$).
To such a set $X\subseteq \bbZ_p^n$, one can associate a
\emph{Poincar\'e series} -- a formal powers series, which has been proven to be a rational function \cite{Den.rat}.
This Poincar\'e series only depends on the isometry type of $X$ and the classification
of isometry types yields a new, geometric proof of its rationality for big $p$.

\medskip

The main result of this article is the existence of t-stratifications in various contexts. To formulate this precisely, we start by fixing some notation.
Let $K$ be a Henselian valued field of equi-characteristic $0$ (i.e., both, $K$ and its residue field have characteristic $0$). We will later fix a suitable class $\mathcal{C}$ of subsets of $K^n$.
One example of a suitable class is the class of all
subvarieties of $K^n$ (not necessarily closed or irreducible; so ``subvariety'' means: locally closed in the Zariski topology); other possible classes $\mathcal C$ are the class of definable subsets in a suitable first order language. In particular, there are classes $\mathcal C$ including analytic subsets of $K^n$.

The goal is to understand the singular locus of sets $X \in \mathcal{C}$.
Roughly, the main theorem states that given such a set $X \subseteq K^n$, $K^n$ can be partitioned into subsets $S_0, \dots, S_n \in \mathcal{C}$ with $\dim S_d = d$ such that
near any point $x \in S_d$, $X$ is ``non-singular in $d$ directions''. Classically, this statement is formalized using
local trivializations. Our result is of a similar nature; we obtain that
on a suitable ball $B$ around $x \in S_d$, the family of sets $(S_d, S_{d+1}, \dots, S_n, X)$ is ``$d$-translatable'', which is a pretty strong notion of trivialization. However, in contrast to local trivializations in the classical sense,
we will specify the size of $B$, which makes our result strictly non-local.

To define $d$-translatability we first have to introduce 
``risometries'', which play a central role in the whole article.
We first fix some more notation:
$\valring$ is the valuation ring of $K$, $k$ is the residue field,
$\Gamma$ is the value group, $\vali\colon K \to \Gamma \cup \{\infty\}$ is the valuation map,
and for $x = (x_1, \dots, x_n) \in K^n$, we define a ``valuation'' $\val(x) := \inf_{i} \vali(x_i)$;
this corresponds to working with the maximum norm on $K^n$. When we speak of a ball in $K^n$, we
mean ball with respect to this maximum norm.
A risometry is something intermediate between an isometry (which preserves the valuation of differences) and a translation (which completely preserves differences); for a ball $B \subseteq K^n$,
a bijection $\phi\colon B \to B$ is a risometry iff it preserves differences up to the leading term, i.e., if
\[\val\big((\phi(y) - \phi(y')) - (y - y')\big) > \val(y - y')\] for all $y, y' \in B$, $y \ne y'$.

We will only be interested in maps $\phi$ that are ``morphisms in the right category''.
If $\mathcal C$ is a class of definable sets, then this simply means that $\phi$, too, should be definable.
However, if $\mathcal C$ is the class of varieties, it doesn't make sense to require $\phi$ to be algebraic,
since balls are not varieties. In that case, we have to work in the category of sets definable in the valued field language (where balls indeed are objects). Expressing this entirely in an algebraic language would be
rather cumbersome but let me at least say that definable means something like ``piecewise algebraic'';
in particular, if $\phi$ is a definable map with $n$-dimensional domain, then the Zariski closure of its graph is also $n$-dimensional.

We call the family of sets $(S_d, S_{d+1}, \dots, S_n, X)$ ``$d$-translatable'' on a ball $B \subseteq K^n$ if
there exists a definable risometry $\phi\colon B \to B$ and a 
$d$-dimensional vector space $V \subseteq K^n$ such that each
set $\phi(S_d \cap B), \dots, \phi(S_n \cap B), \phi(X \cap B)$ is, as a subset of $B$, translation invariant in direction $V$, i.e., it is the intersection of $B$ with a union of cosets of $V$.
In less formal terms, $d$-translatability means that there exists a local trivialization via a risometry.

Now we can formulate a first precise version of the main theorem.

\begin{thm}\label{thm:weakA}
For every set $X \subseteq K^n$ in the class $\mathcal{C}$,
there exists a ``t-stratification of $K^n$ reflecting $X$'',
i.e., a partition $(S_i)_{0 \le i \le n}$ of $K^n$ with
$S_i \in \mathcal{C}$ such that for each $d \le n$, we have the following:
\begin{itemize}
\item
$\dim S_d = d$ or $S_d = \emptyset$
\item
For any ball $B \subseteq S_d \cup \dots \cup S_n$,
the family $(S_d, \dots, S_n, X)$ is $d$-translatable on $B$.
\end{itemize}
\end{thm}

The ``full version'' of this theorem (formulated in the language of model theory) is
Theorem~\ref{thm:main}.
Corollary~\ref{cor:pos-char} is a reformulation which is uniform in the field $K$
and which also works in sufficiently large positive characteristic.
For readers not familiar with model theory,
Theorem~\ref{thm:alg} is a reformulation of Corollary~\ref{cor:pos-char} in a purely algebraic context.
Note that these full versions (\ref{thm:main}, \ref{cor:pos-char}, \ref{thm:alg})
also yield the existence of t-stratifications uniformly in families.

One of the big strengths of Theorem~\ref{thm:weakA} is that we obtain $d$-translatability
not just on a neighborhood of each point of $S_d$ but on any
ball disjoint from $S_0 \cup \dots \cup S_{d-1}$. However, a drawback of this is that
the choice of the $S_d$ can be pretty uncanonical. Here are some examples.

\begin{ex}\label{ex:alm-sing}
For $X = \{(x,y) \in K^2 \mid xy = 0\}$, the obvious stratification works:
$S_0 = \{(0,0)\}$, $S_1 = X \setminus S_0$, and $S_2 = K^2 \setminus X$. 
Now consider the curve $X' = \{(x,y) \in K^2\mid xy = a\}$ for some $a \in K\setminus \{0\}$. Then $X'$ is smooth
(in some naive sense), and this fits to the fact that
each $x \in X'$ has a 1-translatable neighborhood. However, one can check that for any ball
$B$ around $(0,0)$ with $B \cap X' \ne \emptyset$, $X'$ is not 1-translatable
on $B$. Thus to obtain a t-stratification, we are forced to choose a point $s_0 \in B$ and to set $S_0 := \{s_0\}$. This choice
is very uncanonical: if, say, $\vali(a) = 0$, then any $s_0 \in \valring^2$ works.
(After that, $S_1$ and $S_2$ can be defined as expected: $S_1 := X'$, $S_2 :=$ the remainder.)
The intuition behind this phenomenon is that ``from far away, $X'$ looks very similar to $X$'', which has a singularity at $(0,0)$.
Thus $X'$ has an ``almost singularity'' near $(0,0)$, and such almost singularities
are detected by t-stratifications.
A more precise formulation of the similarity between $X$ and $X'$ is that for $B$ as above, 
the images $X_B, X'_B \subseteq (K/B)^2$ from the beginning of the introduction are equal.
\end{ex}

In a model theoretic setting, where $\mathcal{C}$ the class of $\emptyset$-definable sets in a suitable
language, even for $X \subseteq K$, the existence of t-stratifications is not entirely trivial,
as the next example shows.

\begin{ex}\label{ex:ball}
Suppose that 
$X$ is a $\emptyset$-definable ball in $K$. Then $X$ is not translatable on any ball $B$ which strictly contains $X$, hence for $(S_0, S_1)$ to be a t-stratification reflecting $X$,
any such $B$ must contain an element of $S_0$.
This could be achieved by putting some element of $X$ into $S_0$, but it is well possible
that $X \cap \acl(\emptyset) = \emptyset$, i.e., $X$ is disjoint from every finite $\emptyset$-definable set.
However, one can show (e.g.\ using some version of cell decomposition) that this can only happen if
$X$ is an open ball, say, $X = \{x \in K \mid \vali(x - a) > \lambda\}$, and moreover, one can find
a finite $\emptyset$-definable set $S_0$ containing an element $x$ with $\vali(x - a) \ge \lambda$.
In particular, $S_0 \cap B \ne \emptyset$ for every ball $B\supsetneq X$, so this does the job.
\end{ex}

Now let us compare Theorem~\ref{thm:weakA} to classical Whitney stratifications (their definition is recalled
in Subsection~\ref{subsect:defn-whit}).
Fix a point $x$ in some stratum $S_d$, and suppose that $B$ is a ball containing $x$
on which we have $d$-translatability, i.e., we have a $d$-dimensional vector space $V \subseteq K^n$ as explained before Theorem~\ref{thm:weakA}.
In Whitney stratifications, each stratum is smooth, so tangent spaces exist.
With t-stratifications, we do not know whether actual tangent spaces exist, but $V$ can be seen as an ``approximate tangent space'' of $S_d$ at $x$ (approximate tangent spaces are made more precise in Subsection~\ref{subsect:trans}).
Moreover, for any $j \ge d$ and any $y \in B \cap S_j$,
$d$-translatability on $B$ also implies that an approximate tangent space of $S_j$ at $y$ approximately contains $V$.
Thus, formulated sloppily, we have: for any $x \in S_d$, any $j \ge d$, and any $y \in S_j$ close enough to $x$, $T_yS_j$ approximately contains $T_xS_d$. If one replaces ``approximately contains'' by ``contains in the limit for $y \to x$'',
then this statement is essentially the classical Condition~(a) of Whitney.

``Containing
approximately'' sounds like a weaker condition and indeed, t-stratifications do not necessarily
satisfy the straightforward translation of the Whitney conditions to non-Archimedean fields
(in contrast to the stratifications from \cite{CCL.cones}).
However, to be able to apply methods from non-standard analysis, this is exactly the right statement.
As a consequence, if we let
$K$ be a non-standard model of $\bbR$ or $\bbC$ (i.e., a particular valued field whose residue field $k$ is equal to $\bbR$ or $\bbC$, respectively), then any t-stratification of $K^n$ induces a
stratification of $k^n$ that satisfies Condition~(a).

Using this method, we will prove (Theorem~\ref{thm:t->whit}) that t-stratifications induce classical
Whitney stratifications. For this, we also need our t-stratifications to satisfy a
non-standard version of Whitney's Condition~(b). A priori, this is not true: using a kind of non-Archimedean logarithmic spiral, one can construct a t-stratification violating Condition~(b).
However, we are only considering t-stratifications consisting of sets in the class $\mathcal{C}$, and inside this class,
such counter-examples are excluded by the following result.

\begin{thm}\label{thm:weakB}
For every set $X \subseteq K^n$ in our class $\mathcal{C}$ and every $x \in K^n$,
there exists a finite subset $M_x \subseteq \Gamma$ of the
value group such that for any $y \in K^n$ with $\val(y - x) \notin M_x$ and any
ball $B$ containing $y$ but not $x$, $X$ is ``translatable on $B$ in direction
$K\cdot(y-x)$'', i.e., there exists a definable risometry $\phi\colon B \to B$ such that
$\phi(X \cap B)$ is translation invariant in direction $K\cdot(y-x)$.
\end{thm}

The ``full version'' of this theorem is Theorem~\ref{thm:kegel} and Corollary~\ref{cor:whit-b} is 
a non-standard version of Whitney's Condition~(b).
Whereas Theorem~\ref{thm:weakA} only yields the existence of translatability, Theorem~\ref{thm:weakB}
is a strong result about translatability in specific directions. Indeed, formulated sloppily, it says the following.
Fix any $x \in K^n$. Then for almost any $y \in X$ (more precisely: for $y \in X$ at almost any distance from $x$), the approximate tangent space $T_yX$ approximately contains the line $K\cdot(y-x)$.
The order of the quantifiers is important here, i.e., the set of permitted distances depends on $x$;
otherwise we would obtain that almost every tangent space approximately contains almost every line,
which, of course, is absurd.
%
%
%
%
%
%

\medskip

In the Archimedean setting, given a finite family of subsets of $\bbR^n$ or $\bbC^n$,
one can find a single Whitney stratification that simultaneously fits to all those sets.
In valued fields, we can even treat ``small'' infinite families with a single t-stratification.
Here, ``small'' is not meant in the sense of cardinality; instead, a family of sets is small roughly
if it is parametrized by a product of subsets of the residue field $k$
and the value group $\Gamma$. (In contrast, a family parametrized by the valued
field $K$ would be large.) To make sense of this, one needs the language of model theory, i.e., the class
$\mathcal{C}$ should be a suitable class of definable sets. Being able to treat
such infinite families will be crucial in the proof of Theorem~\ref{thm:weakA} to make
the induction work.

\medskip

Now let us come back to the original motivation for this article, namely understanding
sets up to isometry.
The existence of a t-stratification $(S_i)_i$ reflecting a set $X$ implies that
the isometry type of $X$ is rather simple (and even its risometry type), since
all risometries appearing in Theorem~\ref{thm:weakA} can be pieced together to
one single risometry $\phi\colon K^n \to K^n$ such that $\phi(X)$ translation invariant in $d$ directions on
any ball $B \subseteq S_d \cup \dots \cup S_n$. (Note however that this $\phi(X)$ will almost never
lie in $\mathcal C$.) In particular, only few different risometry types occur at all
in $\mathcal{C}$, and t-stratifications help understanding them, in the following precise sense.

Suppose we have a uniform family of sets $X_q \subseteq K^n$ in $\mathcal{C}$,
parametrized by $q \in Q \in \mathcal C$ and suppose we want to decide for which $q,q' \in Q$
there exists a risometry $K^n \to K^n$ sending $X_q$ to $X_{q'}$.
A priori, this is a difficult task;
in model theoretic terms, the induced equivalence relation on $Q$ is not
definable in general. However, if we assume that each $X_q$ is equipped with
a t-stratification $(S_{i,q})_i$ and we ask that these t-stratifications are also
respected by the risometries, then the corresponding refined equivalence relation on $Q$ is definable
(Proposition~\ref{prop:sak}), and the risometry type of $(X_q, (S_{i,q})_i)$ can be described by
a ``finite amount of data living only in $k$ and $\Gamma$ (and not in $K$)''.
A slightly weaker but purely algebraic version of this statement is given
in Corollary~\ref{cor:alg-isotyp}; roughly, there exist finitely many regular functions on $Q$
such that the risometry type only depends on the valuation and the leading term of these functions.

The exact data needed to describe the risometry type of $(X_q, (S_{i,q})_i)$
can be extracted from the proof of Proposition~\ref{prop:sak}. In this way,
one could obtain a complete (but long and technical) classification of all possible risometry types.
In \cite{i.QpB}, such a classification statement has been formulated 
for isometries instead of risometries in the case $K = \bbQ_p$.
As already mentioned, we will prove that conjecture for $p$ sufficiently big
(Theorem~\ref{thm:QpB}). In fact, our proof works for more general $K$;
Definition~\ref{defn:lev} is the classification statement in that setting.
To make ``for $p$ sufficiently big'' precise, we assume
that $X(K) \subseteq K^n$ is given uniformly for all $K$; then we obtain the result
for all $K$ with big enough residue characteristic, where the bound may depend on $X$.

\medskip

Now let me describe for which classes $\mathcal{C}$ Theorems~\ref{thm:weakA} and \ref{thm:weakB} hold. In most of the article
(in particular for the main Theorems~\ref{thm:main} and \ref{thm:kegel} and their proofs),
$\mathcal{C}$ will be the class of definable sets in a suitable language $\Lx$ expanding the language $\LHen$ of valued fields (see Definition~\ref{defn:LHen}).
Hypotheses~\ref{hyp:general} and \ref{hyp:orth} list the precise conditions on $\Lx$
we will need (Hypothesis~\ref{hyp:orth} is only needed for Theorem~\ref{thm:weakB}).
We will prove that these Hypotheses hold in a quite general setting, namely in any
expansion $\Lx$ of $\LHen$ by an analytic structure in the sense of \cite{CL.analyt} (see Propositions~\ref{prop:an->hyp0}, \ref{prop:an->jac}, and \ref{prop:an->orth}).
A concrete example of such an analytic structure on a complete valued field of rank one is given in Example~\ref{ex:an}; it includes all functions given by restricted power series with integer coefficients.
Many more examples of analytic structures can be found in \cite{CL.analyt}.

Most of the assumptions in Hypotheses~\ref{hyp:general} and \ref{hyp:orth} are immediate consequences of
\cite{CL.analyt}; in the case $\LHen$, they also follow easily from classical results like quantifier elimination and cell decomposition. However, the ``higher-dimensional Jacobian property'' is more subtle
and to prove it, we will inductively use the existence of t-stratifications in lower dimensions. To make this possible, the formulation of Theorem~\ref{thm:main} states precisely in which dimension the Jacobian property is needed as a prerequisite. Even in the base case $\LHen$, I do not know a better proof than the above one,
where we work with the ``trivial'' analytic structure (whose existence, by the way, is far from trivial; see
\cite[Section~4.6]{CL.analyt}).

I promised that for $\mathcal{C}$, we can also take the class of varieties. Varieties
are definable in any of the above languages, but the question is whether we can obtain a
t-stratification consisting of varieties and not just of definable sets.
We will show that for $\Lx = \LHen$, this is indeed true---even if the set $X$ we started with is only definable; Corollary~\ref{cor:alg} is a model-theoretic formulation of that result, and
Theorem~\ref{thm:alg} is essentially a precise algebraic formulation of Theorem~\ref{thm:weakA} in the case of varieties.

\medskip

Here is an overview over the article.

In Section~\ref{sect:setting}, we introduce the basic notations and tools.
The first three subsections are independent of model theory. In particular, in Subsection~\ref{subsect:RVlin},
we define and describe the ``higher dimensional leading term structures'' $\RV^{(n)}$, which are ubiquitous in this article.
In Subsections~\ref{subsect:modth} and \ref{subsect:hyp}, we
fix the model theoretic setup and assumptions.

The main purpose of Section~\ref{sect:t-strat} is to define translatability and t-stratifications and to prove some first properties. The most important ones are that the restriction of a t-stratifications to a suitable
affine subspace of a ball is again a t-stratification (Lemma~\ref{lem:t-strat-ball}), that
being a t-stratification is a first-order property (Proposition~\ref{prop:sak}~(1)) and that t-stratifications
can be used to understand risometry types (Proposition~\ref{prop:sak}~(2)).

The next section contains the main part of the proof of Theorem~\ref{thm:main}, namely that
under Hypothesis~\ref{hyp:general}, t-stratifications do exist. There is a sketch of the proof at the beginning of the section. Subsection~\ref{subsect:cor} contains some direct corollaries and the last subsection
gives other characterizations of what it means for a t-stratification to reflect a set; these will be useful for applications and for the inductive proof of the higher-dimensional Jacobian property.

Up to there, we do not yet know whether Hypothesis~\ref{hyp:general} can be satisfied at all;
in Section~\ref{sect:an->hyp}, we will show that it holds in any field with analytic structure in the sense of \cite{CL.analyt}. The easy part is Proposition~\ref{prop:an->hyp0}, which treats everything except the Jacobian property; the latter is proven in Proposition~\ref{prop:an->jac}.

The remaining sections give some variants and applications of the main result, mostly under some additional assumptions.
In Section~\ref{sect:alg}, we show how to obtain t-stratifications
such that for each $d$, $S_0 \cup \dots \cup S_d$ is closed in a suitable topology.
In the pure valued field language $\LHen$, this can be applied to the Zariski
topology, which yields the algebraic version of the main result.

In the next section, we show how our result implies the existence of classical Whitney
stratifications (Theorem~\ref{thm:t->whit}). To this end, we first prove the valued field version of Whitney's Condition~(b)
(Theorem~\ref{thm:kegel}, Corollary~\ref{cor:whit-b}); this needs an additional
hypothesis on the language $\Lx$ (Hypothesis~\ref{hyp:orth}), which also holds in any field with analytic structure (Proposition~\ref{prop:an->orth}).

Finally, we show how the existence of t-stratifications implies the main conjecture of
\cite{i.QpB} about sets up to isometry in $\bbQ_p$ for $p \gg 1$ (Section~\ref{sect:QpB})
and we list some open
questions concerning enhancements of the main result (Section~\ref{sect:open}).

\subsection{Acknowledgment}

I am grateful to numerous people for many fruitful conversations
and also for concrete suggestions concerning this article. A person I want to thank
particularly is Raf Cluckers. At some point, we started to work together on
other questions which I thought would be useful for the present article; it
turned out that after all, these other results are not needed here (only Lemma~\ref{lem:banach}
has been stolen from our common project), but we continued to collaborate,
which always stayed to be a great pleasure.
I also want to thank an anonymous referee, whose comments helped to
enhance some arguments and improve the presentation.

From a financial point of view, I want to thank the
Fondation Sciences Math\'ematiques de Paris which supported me during the first
year I worked on this project and the Deutsche Forschungsgemeinschaft which supported
me during the remainder of the time, via the projects 
SFB~478 ``Geometrische Strukturen in der Mathematik''
and SFB~878 ``Groups, Geometry \& Actions''.

\section{The setting}
\label{sect:setting}


\subsection{Basic notation}
\label{subsect:notn}

In most of the article, we will work in a fixed valued field. We use the following notation.

\begin{notn}\label{notn:VF}
$K$ is a valued field, $\valring$ is the valuation ring, $\maxid \subseteq \valring$ is
its maximal ideal, $k$ is the residue field, $\Gamma$ is the value group,
$\vali\colon K \to \Gamma \cup \{\infty\}$ is the valuation, and $\res\colon \valring \to k$ is the residue map.

Moreover, we define $\val\colon K^n \to \Gamma \cup \{\infty\}, \val(x_1, \dots, x_n) = \min_i \vali(x_i)$
(as in the introduction). We also write $\res$ for the canonical map $\valring^n \to k^n$,
and, more generally, for $X(\valring) \to X(k)$ where $X$ is any variety defined over $\valring$.
\end{notn}

We apply the map $\res$ to sub-vector spaces of $K^n$ as follows.

\begin{defn}\label{defn:resV}
A vector space $V \subseteq K^n$ can be considered as an element of the Grassmanian
$\grass_{n,d}(K) = \grass_{n,d}(\valring)$ (where $d = \dim V$). We write
$\res(V)$ for its image in $\grass_{n,d}(k)$, considered as a sub-vector space of $k^n$.
Equivalently, we have
$\res(V) = \{\res(x) \mid x \in V \cap \valring^n\}$.

Vice versa, if
$V \subseteq k^n$ is a vector space, then any vector space
$\tilde{V} \subseteq K^n$ with $\res(\tilde{V}) = V$ will be called a \emph{lift}
of $V$.
\end{defn}

The map $\val$ on $K^n$ satisfies the ultrametric triangle inequality.
Balls in $K^n$ are defined using this ``metric'', i.e., a ball is the same
as a ``cube'': a product of $n$ balls in $K$ of the same radius. Here is the precise
notion of balls we will use.

\begin{defn}\label{defn:ball}
\begin{enumerate}
\item
An \emph{open ball} in $K^n$ is a set of the form $\ball{a, >\delta} := \{x \in K^n \mid \val(x - a) > \delta\}$ for $a \in K^n$ and $\delta \in \Gamma \cup \{-\infty\}$.
\item
A \emph{closed ball} is a set of the form
$\ball{a, \ge \delta} := \{x \in K^n \mid \val(x - a) \ge \delta\}$
for $a \in K^n$ and $\delta \in \Gamma$.
\item
A \emph{ball} is either an open or a closed ball.
\item
The \emph{radius} of a ball $B$ is the above $\delta$; we denote it by
$\rado(B)$ if $B$ is an open ball and by $\radc(B)$ if $B$ is a closed ball.
\end{enumerate}
\end{defn}

Thus: we do consider $K^n$ as a ball (an open one), but we do not consider points as balls,
and neither do we allow arbitrary cuts in $\Gamma$ as radii of balls.
The reason to have two different notations $\rado$ and $\radc$ is
that if $\Gamma$ is discrete, then any ball $B \ne K^n$ can be considered both as
an open or as a closed ball and $\rado(B) < \radc(B)$.

From time to time, given a ball $B$ we will need to consider the ball $B'$ of the same radius containing the origin.
We do not introduce a special notation for this;
instead, note that $B' = B - B = \{b - b' \mid b, b' \in B\}$.

We will work a lot with projections
$\pi\colon K^n \surject K^d$ to some subset of the coordinates. The corresponding
projection $k^n \surject k^d$ at the level of the residue field will
be denoted by $\pibar$. By $\picomp\colon K^n \surject K^{n-d}$, we will denote
the projection to the complementary set of coordinates.
Often, we will consider restricted coordinate projections $\pi \colon B \to K^d$ for some
subset $B \subseteq K^n$ (most of the time, a ball); in that case,
$\pibar$ still denotes the entire projection $k^n \surject k^d$.

Given a coordinate projection $\pi \colon B \to K^n$,
any fiber $\pi\1(x)$ (for $x \in \pi(B)$) can be identified with a subset of $K^{n-d}$ via $\picomp$. Using this, any definition made for $K^{n-d}$ can also be applied to fibers of coordinate projections.
As an example, this yields a notion of a ball inside a fiber $\pi\1(x)$.

\subsection{Higher dimensional leading term structures and their linear algebra}
\label{subsect:RVlin}

The ``leading term structure'' is usually defined as $\RV := \{0\} \cup K\mult / (1 + \maxid)$.
We will need the following higher dimensional version of it.

\begin{defn}\label{defn:RVn}
Define $\RV^{(n)} := K^n/\mathord{\sim}$, where $x \sim y \iff (\val(x - y) > \val(x) \,\vee\, x = y = 0)$; we
write $\rv\colon K^n \surject \RV^{(n)}$ for the canonical map
and $\vRV$ for the map from $\RV^{(n)}$ to $\Gamma\cup \{\infty\}$ satisfying $\vRV \circ \rv = \val$.
Instead of $\RV^{(1)}$, we also write $\RV$.
\end{defn}

The following is a (more general) coordinate free version of this definition.

\begin{defn}\label{defn:RVL}
Let $L$ be a free $\valring$-module and set $V_L := K \otimes_{\valring} L$.
First, define the valuation $\val\colon V_L \to \Gamma \cup \{\infty\}$
by setting $\val(rx) := \val(r)$ for any $r \in K, x \in L \setminus \maxid L$. (This is well-defined and satisfies the ultrametric triangle inequality.)
Then set $\RV_L:= V_L/\mathord{\sim}$ where $x \sim y \iff (\val(x - y) > \val(x) \,\vee\, x = y = 0)$,
write $\rv\colon V_L \surject \RV_L$ for the canonical projection, and 
write $\vRV$ for the map from $\RV_L$ to $\Gamma\cup \{\infty\}$ satisfying $\vRV \circ \rv = \val$.
\end{defn}

It is easy to check that we have $\RV^{(n)} = \RV_{\valring^n}$ and
that the two definitions of $\val$ and $\rv$ on $K^n$ coincide.
In most of the article,
we will work with coordinates anyway, so we will not bother giving coordinate free definitions
of everything. Note however that using Definition~\ref{defn:RVL}, one obtains the following for free.

\begin{lem}\label{lem:RVfunct}
If $L$ is a free $\valring$-module, then 
any $\phi \in \Aut(L)$ induces a map $\phi\colon \RV_L \to \RV_L$ satisfying
$\phi(\rv(x)) = \rv(\phi(x))$ for $x \in V_L$ (where we also write $\phi$ for the induced map $V_L \to V_L$).
In particular, any $M \in \GL_n(\valring)$ induces a map $M\colon \RV^{(n)} \to \RV^{(n)}$.
\end{lem}

We will need one more notation.

\begin{defn}\label{defn:dir}
For $x \in K^n \setminus \{0\}$, let the \emph{direction} of $x$ be the one-dimensional
subspace $\dir(x) := \res(K \cdot x)$ of $k^n$, considered
as an element of the projective space $\bbP^{n-1} k$.
Notationally, we will almost always treat $\dir(x)$ as a representative $y \in \res(K \cdot x)$ of the actual direction.
Whenever we will use this notation, we will make sure that the particular choice of $y$ doesn't
matter.

One easily verifies that the direction map factors over $\RV^{(n)}$; we write $\dirRV$ for the corresponding map
$\RV^{(n)} \to \bbP^{n-1} k$ (i.e., $\dirRV \circ \rv = \dir$).
\end{defn}

A lot of commutative diagrams can be drawn, showing how all these maps fit together. The
following two lemmas shows only some of them. Lemma~\ref{lem:Un} in particular
shows that $\RV^{(n)}$ can also be defined as the quotient of $K^n$ by a suitable group action, generalizing
the one-dimensional case $\RV = K/(1+\maxid)$.

\begin{lem}\label{lem:Un}
Let $U_n$ be the kernel of the map $\res\colon \GL_n(\valring) \surject \GL_n(k)$.
Then we have the following commutative diagrams, where $G \looparrowright X$ means that
$G$ acts on $X$, and each straight line $G \looparrowright X \surject Y$ is exact in the sense
that $Y$ is the quotient of $X$ by the action of $G$.
\begin{center}
\begin{tikzpicture}[x=1cm,y=1.5cm,baseline=1.5cm]
\node (U) at (0,2) {$U_n$};
\node (GLO) at (2,2) {$\GL_n(\valring)$};
\node (GLk) at (5,2) {$\GL_n(k)$};
\node (RV) at (3.8,.8) {$\RV^{(n)}$};
\node (Ga) at (3,0) {$\Gamma \cup \{\infty\}$};
\path [name path=p1] (U) -- (RV);
\path [name path=p2] (GLO) -- (Ga);
\node [name intersections={of=p1 and p2}] (K) at (intersection-1) {$K^n$};

\draw[right hook->] (U) -- (GLO); \draw[->>] (GLO) -- node[above]{$\res$} (GLk);
\draw[acts->] (GLO) -- (K); \draw[->>] (K) -- node[left]{$\val$} (Ga);
\draw[acts->] (U) -- (K); \draw[->>] (K) -- node[above]{$\rv$} (RV);
\draw[acts->] (GLk) -- (RV); \draw[->>] (RV) -- node[right]{$\vRV$} (Ga);
\end{tikzpicture}
\begin{tikzpicture}[x=1cm,y=1.5cm]
\node (U) at (0,1) {$U_n$};
\node (PO) at (1.8,1) {$\grass_{n,d}(\valring)$};
\node (Pk) at (4,1) {$\grass_{n,d}(k)$};
\draw[acts->] (U) -- (PO); \draw[->>] (PO) -- node[above]{$\res$} (Pk);
\end{tikzpicture}
\end{center}
\end{lem}

\begin{proof}
None of this is difficult to show. As an example, let us verify that if
$V, V' \subseteq K^n$ are $d$-dimensional vector spaces with $\res(V) = \res(V')$,
then there exists $M \in U_n$ with $MV = V'$.

Choose any basis $(b_i)_{i \le d}$ of $\res(V)$ and extend it to a basis $(b_i)_{i \le n}$ of $k^n$.
Choose preimages $(v_i)_{i \le n}$ and $(v'_i)_{i \le n}$ of $(b_i)_{i \le n}$ in $\valring^n$
such that for $i \le d$, we have $v_i \in V$ and $v'_i \in V'$. Then the linear map
sending $(v_i)_{i \le n}$ to $(v'_i)_{i \le n}$ sends $V$ to $V'$, and it lies in $U_n$
since induces the identity on $k^n$.
\end{proof}

\begin{lem}
We also have the following commutative diagram, where each straight line is exact. Here,
$\RV\setminus \{0\} \looparrowright \RV^{(n)}\setminus \{0\}$
is an action induced by the scalar multiplication $K\mult \looparrowright K^n\setminus \{0\}$ and
$\Gamma \looparrowright \Gamma$ is the action by translation.
The middle horizontal line is exact in the sense that $k^n \setminus \{0\} = \vRV\1(0)$.

\begin{center}
\begin{tikzpicture}[x=1cm,y=1.2cm]
\node (k) at (0.5,2) {$k\mult$};
\node (RV) at (3,2) {$\RV\setminus \{0\}$};
\node (Ga) at (5,2) {$\Gamma$};
\node (kn) at (0.5,1) {$k^n \setminus \{0\}$};
\node (RVn) at (3,1) {$\RV^{(n)} \setminus \{0\}$};
\node (Ga2) at (5,1) {$\Gamma$};
\node (Pk) at (0.5,0) {$\bbP^{n-1}k$};
\node (Pk2) at (3,0) {$\bbP^{n-1}k$};

\draw[right hook->] (k) -- (RV); \draw[->>] (RV) -- node[above]{$\viRV$} (Ga);
\draw[right hook->] (kn) -- (RVn); \draw[->>] (RVn) -- node[above]{$\vRV$} (Ga2);
\draw[double distance=1.5pt] (Pk) -- (Pk2);
\draw[acts->] (k) -- (kn); \draw[->>] (kn) -- (Pk);
\draw[acts->] (RV) -- (RVn); \draw[->>] (RVn) -- node[right]{$\dirRV$} (Pk2);
\draw[acts->] (Ga) -- (Ga2);
\end{tikzpicture}
\end{center}
\end{lem}

\begin{proof}
Easy.
\end{proof}
Note that the top right square of the diagram implies that $\vRV\colon \RV^{(n)} \setminus \{0\}
\surject \Gamma$ is a fibration with fibers ``isomorphic'' to $k^n \setminus \{0\}$.

Here are some more basic properties of $\RV^{(n)}$ and the maps defined above.

\begin{lem}\label{lem:vf}
\begin{enumerate}
\item\label{it:rv-sum}
If $a_1, a_2 \in K^n$ satisfy $\val(a_1 + a_2) = \min\{\val(a_1), \val(a_2)\}$, then
$\rv(a_1)$ and $\rv(a_2)$ together determine $\rv(a_1 + a_2)$,
i.e., for any other $a'_1, a'_2 \in K^n$ with $\rv(a'_i) = \rv(a_i)$,
we have $\rv(a'_1 + a'_2) = \rv(a_1 + a_2)$.
\item\label{it:dir-sum}
If $a_1, a_2 \in K^n$ satisfy $\dir(a_1) \ne \dir(a_2)$, then
$\dir(a_1 + a_2)$ lies in the $k$-vector space spanned by $\dir(a_1)$ and $\dir(a_2)$.
\item\label{it:dir-pi}
Suppose that $\pi \colon K^n \surject K^d$ is a coordinate projection,
$\pibar \colon k^n \surject k^d$ is the corresponding projection at the level
of the residue field, and $a \in K^n \setminus \{0\}$. Then we have $\val(\pi(a)) = \val(a)$ iff $\pibar(\dir(a)) \ne 0$.
Moreover, in that case $\pibar(\dir(a)) = \dir(\pi(a))$,
and if $a' \in K^n$ is another element with $\pi(a') = \pi(a)$ and $\dir(a') = \dir(a)$,
then we have $\rv(a') = \rv(a)$.
\item\label{it:dir-res}
For any sub-vector space $V \subseteq K^n$, we have $\dir(V) = \res(V)$.
\item\label{it:dir-scal}
Let $\langle \cdot, \cdot \rangle$ denote the standard scalar product, both
on $K^n$ and on $k^n$. Then for any $a, b \in K^n$ we have
$\val(\langle a, b \rangle) \ge \val(a) + \val(b)$ and we have the equivalence
$\val(\langle a, b \rangle) > \val(a) + \val(b) \iff
\langle \dir(a), \dir(b) \rangle = 0$.
\end{enumerate}
\end{lem}
\begin{proof}
Easy.
\end{proof}

\subsection{Risometries}
\label{subsect:riso}

Let us now have a look at the notion of risometry, which already appeared
in the introduction. Its definition can be written down nicely using the multidimensional $\rv$-map introduced in Definition~\ref{defn:RVn}.

\begin{defn}\label{defn:riso}
For $X_1, X_2 \subseteq K^n$, a \emph{risometry} from $X_1$ to $X_2$ is a
bijection $\phi\colon X_1 \to X_2$ satisfying $\rv(\phi(x) - \phi(x')) = \rv(x - x')$
for any $x, x' \in X_1$. If $\phi$ is such a risometry, we use the following terminology.
\begin{itemize}
\item
For maps $\chi_i$ with domain $X_i$ (for $i=1,2$),
we say that $\phi$ is a risometry from $\chi_1$ to $\chi_2$ (or: $\phi$ sends $\chi_1$ to $\chi_2$)
if $\chi_1 = \phi \circ \chi_2$.
More generally, if $(\theta_{1,\nu})_\nu$ and $(\theta_{2,\nu})_\nu$ are tuples where
for each $\nu$, $\theta_{1,\nu}$ and $\theta_{2,\nu}$ are either maps with domain $X_1$ and $X_2$ or subsets of $X_1$ and $X_2$, then
we say that $\phi$ sends $(\theta_{1,\nu})_\nu$
to $(\theta_{2,\nu})_\nu$ (and we sometimes write $\phi\colon (\theta_{1,\nu})_\nu \to (\theta_{2,\nu})_\nu$) if $\phi$ sends
$\theta_{1,\nu}$ to $\theta_{2,\nu}$ for each $\nu$.
\item
If $\chi$ is a map whose domain contains $X_1 \cup X_2$, we also say that $\phi$ respects $\chi$
if it sends $\chi\auf{X_1}$ to $\chi\auf{X_2}$, and $\chi$ respects a set $Y \subseteq K^n$
if it sends $Y \cap X_1$ to $Y \cap X_2$.
\end{itemize}
\end{defn}

As in ``$\rvi$'', the ``r'' in ``risometry'' stands for ``residue field''.
The condition about $\rv$ in Definition~\ref{defn:riso} already implies injectivity, so any map satisfying
that condition is a risometry from its domain to its image.
Note also that the composition of risometries is again a risometry;
in particular, the risometries from a set to itself form a group.

\begin{rem}\label{rem:risoGLO}
Lemma~\ref{lem:RVfunct} implies that 
if $\phi\colon X \to Y$ is a risometry and $M \in \GL_n(\valring)$,
then we also have a risometry $M \circ \phi \circ M\1 \colon M(X) \to M(Y)$.
This will be used from time to time to ``without loss change coordinates''.
In particular, any matrix $\bar{M} \in \GL_n(k)$ can be lifted to a matrix
$M \in \GL_n(\valring)$, so we can apply any coordinate transformation at the level
of the residue field.
\end{rem}

\begin{rem}\label{rem:risoU}
The group $U_n \subseteq \GL_n(\valring)$ introduced in Lemma~\ref{lem:Un}
consists exactly of those linear maps $K^n \to K^n$ which are risometries.
In particular, if $V_1, V_2 \subseteq K^n$ are vector spaces
with $\res(V_1) = \res(V_2)$, then
there exists a risometry $K^n\to K^n$ sending $V_1$ to $V_2$.
\end{rem}

Cartesian products of risometries are again risometries;
the following lemma strengthens this a bit.
Each of its statements is almost trivial (so we omit the proof), but together,
they will be useful to construct risometries.

\begin{lem}\label{lem:risoKomb}
Let $V_1, V_2 \subseteq K^n$ be sub-vector spaces such that
we have a direct sum decomposition $\res(V_1) \oplus \res(V_2) = k^n$; write $\pi_i\colon K^n \to K^n$
for the projection with image $V_i$ and kernel $V_{3-i}$ (for $i = 1,2$).
\begin{enumerate}
\item 
Suppose that $X \subseteq K^n$ and that $\phi_1,\phi_2\colon X \to K^n$
are maps satisfying
\[
\tag{$*_i$}
\val(\phi_i(x) - \phi_i(x') - \pi_i(x - x')) > \val(x - x') \qquad \text{for every } x,x' \in X, \, x \ne x'
.
\]
Then the map $x \mapsto \phi_1(x) + \phi_2(x)$ is a risometry from $X$ to its image.
\item
If $\phi$ is a map satisfying $(*_i)$ and $\psi$ is a risometry (with suitable domain and image), then
$\phi \circ \psi$ and $\psi \circ \phi$ also satisfy $(*_i)$;
in particular, $\pi_i \circ \psi$ and $\psi \circ \pi_i$ satisfy $(*_i)$.
\end{enumerate}
\end{lem}

Next, we describe risometries between finite sets and how
such risometries can be extended to larger sets.
In the following, for $x \in K^n$ and $T \subseteq K^n$, the notation
$\rv(x - T)$ means $\{\rv(x - t) \mid t \in T\}$.

\begin{lem}\label{lem:fin&iso}
Let $T \subseteq K^n$ be a finite set.
\begin{enumerate}
\item
The only risometry $T \to T$ is the identity. (In particular, between
two different finite sets, there is at most one risometry.)
\item
For $x_1, x_2 \in K^n$ $x_1 \ne x_2$, the following are equivalent:
\begin{enumerate}
\item
There exists a risometry $\phi\colon K^n \to K^n$ with $\phi(T) = T$ and $\phi(x_1) = x_2$.
\item
$\ball{x_1, \ge \val(x_1 - x_2)} \cap T = \emptyset$.
\item
$\rv(x_1 - T) = \rv(x_2 - T)$.
\end{enumerate}
In particular, the fibers of the map sending $x \in K$ to the set $\rv(x - T)$ are exactly
the singletons $\{t\}$ for $t \in T$ and
the maximal balls $B \subseteq K$ that are disjoint from $T$.
\item
A map $\phi\colon K^n \to K^n$ which is the identity on $T$ is a risometry
if and only if for each maximal ball $B \subseteq K^n \setminus T$,
the restriction $\phi\auf{B}$ is a risometry from $B$ to itself.
\end{enumerate}
\end{lem}
\begin{proof}
(2)
``(a) $\Rightarrow$ (c)'' and ``(b) $\Rightarrow$ (a)'' are trivial.
(For the latter, define $\phi$ to be the translation by $x_2 - x_1$ on $\ball{x_1, \ge \val(x_1 - x_2)}$
and the identity everywhere else.)

``(c) $\Rightarrow$ (b)'': Without loss, $\val(x_1 - x_2) = 0$ and $x_1, x_2 \in \valring^n$. Suppose for contradiction that $T_0 := T \cap \valring^n$ is non-empty. The assumption (c) implies $\rv(x_1 - T_0) = \rv(x_2 - T_0)$ and hence $\res(x_1 - T_0) = \res(x_2 - T_0)$.
This implies 
\[
\sum_{\bar{t} \in \res(T_0)}(\res(x_1) - \bar{t})
=
\sum_{\bar{t} \in \res(T_0)}(\res(x_2) - \bar{t})
.
\]
Adding $\sum_{\bar{t} \in \res(T_0)} \bar{t}$ and then dividing by $|\res(T_0)|$ on both sides yields $\res (x_1) = \res (x_2)$, which contradicts $\val(x_1 - x_2) = 0$.

The ``in particular'' part of (2) follows from (b) $\iff$ (c).

(1) If $\phi\colon T \to T$ is a risometry, then for any $t \in T$ we have
$\rv(t - T) = \rv(\phi(t) - \phi(T)) = \rv(\phi(t) - T)$.
Suppose that $\phi(t) \ne t$. Then (2) ``(c) $\Rightarrow$ (b)'' yields
$\ball{t, \ge \val(t - \phi(t))} \cap T = \emptyset$, which contradicts $t \in T$.

(3)
``$\Rightarrow$'' follows from (2), (a) $\Rightarrow$ (b). For ``$\Leftarrow$'', suppose
that $\phi\auf{B}$ is a risometry $B \to B$ for each maximal ball $B \subseteq K^n \setminus T$;
we have to verify that $\rv(x - x') = \rv(\phi(x) - \phi(x'))$ for every
$x, x' \in K^n$. If $x$ and $x'$ lie in the same maximal ball $B$,
then there is nothing to show. Otherwise, we have $\val(x - x') < \val(x - \phi(x))$ and
$\val(\phi(x) - x') < \val(x'-\phi(x'))$, which implies $\rv(x - x') = \rv(\phi(x) - x')
= \rv(\phi(x) - \phi(x'))$.
\end{proof}

\subsection{Model theoretic conventions and setting}
\label{subsect:modth}

Unless specified otherwise, ``definable'' will always mean definable with parameters.
There will be some results concerning $\emptyset$-definable sets (of the form:
for some $\emptyset$-definable $X$, there exists a $\emptyset$-definable $Y$\dots).
Our general assumptions about the language and the theory will always allow to add parameters to the language
(see Remark~\ref{rem:param}), so the reason to write ``$\emptyset$-definable''
is only to emphasize that $Y$ is definable over the same parameters as $X$.

We start by fixing a basic language $\LHen$
for valued fields and a corresponding theory $\THen$.
In almost all of the article, we only care about the language up to interdefinability;
however, we will have to be precise about the sorts of the language.
We use the notation introduced in Subsections~\ref{subsect:notn} and \ref{subsect:RVlin}. 

\begin{defn}\label{defn:LHen}
\begin{enumerate}
\item
Let $\LHen$ be the language consisting of one sort $K$ for the valued field with the ring language,
all sorts $\RV\eq$ with the corresponding canonical maps between them, and the map $\rvi\colon K \to \RV$.
More precisely, the sorts of $\RV\eq$ are all sets of the form $X/\mathord{\sim}$, where
$X \subseteq \RV^k$ is $\emptyset$-definable and $\sim$ is a $\emptyset$-definable equivalence relation on $X$,
and the corresponding canonical map is $X \surject X/\mathord{\sim}$.
\item
We call $\RV\eq$ the \emph{auxiliary sorts}.
By an \emph{auxiliary set} resp.\ \emph{element}, we mean a subset resp.\ element of an auxiliary sort.
\item
Let $\THen$ be the theory of Henselian valued fields of equi-characteristic $0$ in the language $\LHen$.
\end{enumerate}
\end{defn}


Notationally, we will often treat $\RV\eq$ as the union of all auxiliary sorts.
In particular, by a ``definable map $\chi \colon K^n \to \RV\eq$'', we mean a
definable map whose target is an arbitrary auxiliary sort (and similarly for
definable sets $Q \subseteq \RV\eq$).

\begin{rem}
One easily checks that $k$, $\Gamma$, and $\RV^{(n)}$ are auxiliary sorts. (For the latter, note that
the map $K^n \surject \RV^{(n)}$ factors over $\RV^n$.)
\end{rem}

\begin{notn}\label{notn:code}
If $(X_q)_{q \in Q}$ is a definable family of sets (or maps), then $\code{X_q}$ denotes a ``code'' for $X_q$.
More precisely, if $X_q$ is defined by a formula $\phi(x, q)$,
then there exists a definable map $f\colon Q \to Q'$ for some definable set $Q'$ (possibly imaginary) and a formula $\psi(x,y)$ such that $\psi(x, f(q))$
also defines $X_q$ and such that $f(q)$ is a canonical parameter for $X_q$.
We set $\code{X_q} := f(q)$. (Of course, this involves some choices.)
\end{notn}

Most of the time when we will use Notation~\ref{notn:code}, we will make sure that $Q'$ can be chosen
in a non-imaginary sort of $\LHen$.


\subsection{Requirements on the theory}
\label{subsect:hyp}

In most of the article, we will not work with $\THen$ and $\LHen$ themselves, but with
an expansion $\Tx$ of $\THen$ in a language $\Lx \supseteq \LHen$ (which has the same sorts as $\LHen$).
In particular,
the main theorem will be proven in any expansion $\Tx$ of $\THen$ having certain
properties, which will be listed in Hypothesis~\ref{hyp:general}.
Variants of these properties have already been introduced and described in \cite{CL.bmin} and \cite{CL.analyt}:
``b-minimality'' is a list of axioms designed to yield cell decomposition
and a notion of dimension, and the ``Jacobian property'' imposes
conditions on definable functions in one variable. Our version of the Jacobian property
includes definable functions in several variables. In addition to these two properties, we will require
that zero-dimensional sets are finite and that $\RV$ is stably embedded.

We start by defining our version of the Jacobian property. Note that even in the one-variable
case, it does not entirely agree with \cite[Definitions 6.3.5, 6.3.6]{CL.analyt}.

\begin{defn}\label{defn:jac}
\begin{enumerate}
\item
Let $X \subseteq K^n$ be a set. We say that a map $f \colon X \to K$ has the \emph{Jacobian property} (on $X$), if either it is constant, or there exists a $z \in K^n \setminus \{0\}$ such that for every $x, x' \in X$ with $x \ne x'$, we have
\[
\vali(f(x) - f(x') - \langle z, x - x'\rangle) > \val(z) + \val(x - x')
.
\]
\item
We say that an expansion $\Tx$ of $\THen$ has the \emph{Jacobian property} if for every model $K \models \Tx$, for every set $A \subseteq K \cup \RV\eq$, for every $n \in \bbN$, and
for every $A$-definable map $f\colon K^n \to K$, there exists
an $A$-definable map $\chi\colon K^n \surject Q \subseteq \RV\eq$ such that
for each $q \in Q$, if $\chi\1(q)$ contains a ball, then $f\auf{\chi\1(q)}$ has the Jacobian property.
\item
If (2) only holds for $n \le n_0$
(where $n_0 \in \bbN$), then we say that $\THen$ has the \emph{Jacobian property up to dimension $n_0$}.
\end{enumerate}
\end{defn}

We will associate dimensions to definable sets in Definition~\ref{defn:dim};
the condition in (2) that $\chi\1(q)$ contains a ball will be equivalent to $\dim(\chi\1(q))=n$.

\begin{rem}\label{rem:jac}
In Definition~\ref{defn:jac}~(1), replacing $z$ by any other $z'\in K^n$ satisfying $\rv(z) = \rv(z')$ does not change the validity of the inequation, since $\vali(\langle z, x - x'\rangle - \langle z', x - x'\rangle) > \val(z) + \val(x - x')$ by Lemma~\ref{lem:vf} (\ref{it:dir-scal}).
\end{rem}

Now we can summarize the prerequisites needed for Theorem~\ref{thm:main}.
We also introduce a notation for a weakening of the hypothesis, which we will need in an inductive argument.

\begin{hyp}\label{hyp:general}
We assume that $\Tx$ is an expansion of $\THen$ in a language
$\Lx \supseteq \LHen$ that has the same sorts as $\LHen$, with the following properties.
\begin{enumerate}
\item\label{it:st-emb}
$\RV$ is stably embedded, i.e., in every model of $\Tx$, every definable subset of $\RV^n$ is definable using
only parameters from $\RV$.
\item\label{it:fin}
In every model $K \models \Tx$, every definable map from $\RV$ to $K$ has finite image.
\item\label{it:strat1}
For every model $K \models \Tx$, for every set $A \subseteq K \cup \RV\eq$ of parameters,
and for every $A$-definable set $X \subseteq K$, there exists a finite, $A$-definable
set $S_0 \subseteq K$ such that every ball $B \subseteq K \setminus S_0$ is either contained
in $X$ or disjoint from $X$.
\item\label{it:jac}
$\Tx$ has the Jacobian property (Definition~\ref{defn:jac}).
\end{enumerate}
For $n \in \bbN$, we write ``Hypothesis~\ref{hyp:general}$_n$'' for the following weakening of this hypothesis.
\begin{enumerate}
\item[(\ref{it:st-emb})] -- (\ref{it:strat1}) as above.
\item[($4'$)] $\Tx$ has the Jacobian property up to dimension $n$.
\item[($4''$)]
For every model $K \models \Tx$, for every set $A \subseteq K \cup \RV\eq$, and
for every $A$-definable map $f\colon K \to K$, there exists
an $A$-definable map $\chi\colon K^n \surject Q \subseteq \RV\eq$ such that for each $q \in Q$,
$f\auf{\chi\1(q)}$ is either injective or constant.
\end{enumerate}
\end{hyp}

Note that Condition~($4''$) is
relevant only in the case $n = 0$, since it follows from the Jacobian property in dimension $1$.


\begin{rem}\label{rem:param}
All conditions in Hypothesis~\ref{hyp:general} remain true if we add constant symbols to the language.
In particular,
any result proven for $\emptyset$-definable sets automatically also holds
over any parameter set $A$. This will be used throughout the paper without further
mentioning.
\end{rem}

\begin{rem}\label{rem:code}
By Hypothesis~\ref{hyp:general}~(\ref{it:st-emb}), for any definable set $Q \subseteq \RV\eq$ we may
assume $\code{Q} \in \RV\eq$, and similarly $\code{f} \in \RV\eq$ for a definable map $f \colon Q \to \RV\eq$.
\end{rem}

\begin{rem}\label{rem:1-dim}
Hypothesis~\ref{hyp:general}~(\ref{it:strat1}) exactly says that t-stratifications
exists for subsets of $K$. Note also that by Lemma~\ref{lem:fin&iso}, the condition relating $S_0$ and $X$
is equivalent to: $X$ is a union of fibers of the map $K \to \RV\eq, x \mapsto \code{\rvi(x - S_0)}$.
\end{rem}

Hypothesis~\ref{hyp:general} does not mention the notion of b-minimality from \cite{CL.bmin} explicitly.
Since we will use results from \cite{CL.bmin} about the existence of a good notion of dimension,
we conclude this subsection by proving that Hypothesis~\ref{hyp:general}$_0$ implies
b-minimality. The following definition is \cite[Definition~2.2.1]{CL.bmin}, applied
to the context of valued fields with auxiliary sorts $\RV\eq$.
Note that it is not exactly the same as \cite[Definitions~6.3.1]{CL.analyt},
since there, only $\RV$ is used as an auxiliary sort. We will come back to this difference when it becomes an issue,
namely in Section~\ref{sect:an->hyp}, when we prove that analytic structures in the sense
of \cite{CL.analyt} satisfy Hypothesis~\ref{hyp:general}.


\begin{defn}\label{defn:b-min-eq}
An expansion $\Tx$ of $\THen$ is \emph{b-minimal over $\RV\eq$} if for every model $K \models \Tx$
and every set $A \subseteq K \cup \RV\eq$ of parameters, the following holds.
\begin{enumerate}
\item
For every $A$-definable set $X \subseteq K$, there exists an $A$-definable
map $\chi\colon X \surject Q \subseteq \RV\eq$
such that every fiber $\chi\1(q)$ (for $q \in Q$) is either a point or an open ball.
\item
There exists no surjective definable map from an auxiliary set
to an open ball $B \subseteq K$.
\item
For every $A$-definable $X\subseteq K$ and $f\colon X \to K$, there exists
an $A$-definable map $\chi\colon X \surject Q \subseteq \RV\eq$ such that for each $q \in Q$,
$f\auf{\chi\1(q)}$ is either injective or constant.
\end{enumerate}
\end{defn}

\begin{lem}\label{lem:b-min-eq}
Hypothesis~\ref{hyp:general}$_0$ implies b-minimality over $\RV\eq$.
\end{lem}

\begin{proof}
(2) follows from Hypothesis~\ref{hyp:general}~(\ref{it:fin}), (3) is exactly (4'') from Hypothesis~\ref{hyp:general}$_0$,
and (1) can be deduced as follows.
Let $X \subseteq K$ be given and let $S_0 \subseteq K$ be a finite set as in the Hypothesis~\ref{hyp:general}~(\ref{it:strat1}).
Then by Lemma~\ref{lem:fin&iso},
$X$ is a union of fibers of the map $\chi\colon K \to \RV\eq, x \mapsto \code{\rv(x - S_0)}$ and
each such fiber is either a point or a ball; thus $\chi\auf{X}$ does the job.
\end{proof}

As mentioned in the introduction, we will prove (Proposition~\ref{prop:an->jac}) that 
Henselian valued fields with analytic structure in the sense of \cite{CL.analyt} satisfy Hypothesis~\ref{hyp:general}.
The proof of the Jacobian property will inductively use the existence of t-stratifications in lower dimensions.
To make this precise, everywhere in
the proof of the existence of t-stratifications,
we will specify for which $n$ Hypothesis~\ref{hyp:general}$_n$ is needed. (Most of the time, Hypothesis~\ref{hyp:general}$_0$
will already be enough.)

\subsection{First consequences: dimension and spherically completeness}
\label{subsect:dim}

We assume Hypothesis~\ref{hyp:general}$_0$. By Lemma~\ref{lem:b-min-eq}, this implies b-minimality (over $\RV\eq$), and 
by \cite{CL.bmin}, this implies the existence of a good notion of dimension of
definable sets, which in particular satisfies the axioms given in \cite{Dri.dimDef}.

\begin{defn}\label{defn:dim}
Let $X \subseteq K^n$ be a definable set.
The \emph{dimension} $\dim X$ is the maximal $d$ such
that there exists a coordinate projection $\pi\colon K^n \surject K^d$
such that $\pi(X)$ contains a ball. We set $\dim \emptyset := -\infty$.
For $x \in K^n$, the \emph{local dimension} of $X$ at $x$ is
$\dim_x X := \min\{\dim (X \cap \ball{x, >\gamma}) \mid \gamma \in \Gamma\}$.
\end{defn}

\begin{rem}\label{rem:0fin}
By Hypothesis~\ref{hyp:general}~(\ref{it:strat1}), any subset of $K$ is either finite or it contains a ball,
hence $0$-dimensional subsets of $K$ are finite. This is also true for subsets of $K^n$, as
one sees by applying coordinate projections $K^n \surject K$.
\end{rem}

It is clear that dimension is definable, i.e., if $X_q \subseteq K^n$ is a $\emptyset$-definable family
of sets (for $q \in Q$), then $\{q \in Q \mid \dim X_q = d\}$ is $\emptyset$-definable
for every $d$. Moreover, we have the following.

\begin{lem}[\cite{CL.bmin}, \cite{Dri.dimDef}]\label{lem:dim}
Dimension has the following properties:
\begin{enumerate}
\item\label{it:dim-u} If $(X_q)_{q \in Q}$ is a definable family of subsets of $K^n$ and
$Q \subseteq \RV\eq$ is auxiliary, then $\dim \bigcup_{q \in Q}X_q =
\max_{q}\dim X_q$.
\item\label{it:dim-f} If $X \subseteq K^m$, $Y \subseteq K^n$ and $f\colon X \to Y$ are definable and
each fiber $f\1(y)$ has dimension $d$, then $\dim X = \dim Y + d$.
\end{enumerate}
\end{lem}

The following property of local dimension is \cite[Theorem~3.1]{iF.dim}. (That Theorem only requires dimension to satisfy some very general axioms which follow directly from our definition and Lemma~\ref{lem:dim}~(\ref{it:dim-f}).)
\begin{lem}\label{lem:loc-dim}
Let $X \subseteq K^n$ be a definable set and set $Y := \{x \in X \mid \dim_x X < \dim X\}$. Then $\dim Y < \dim X$.
\end{lem}

Hypothesis~\ref{hyp:general}$_0$ also implies that $K$ is
``definable spherically complete'' in the following sense.

\begin{lem}\label{lem:dsc}
For every definable family $(B_q)_{q \in Q}$ of balls $B_q \subseteq K$ which
form a chain with respect to inclusion, the intersection $\bigcap_{q \in Q} B_q$
is non-empty.
\end{lem}

\begin{proof}
Let such a family $(B_q)_{q \in Q}$ be given. We can assume that there is no smallest ball.

Let $S_q$ be the finite set obtained by applying Hypothesis~\ref{hyp:general} (\ref{it:strat1})
to $B_q$ and set $S'_q := \{x\in S_q \mid \exists \xi \in \RV\colon x + \rvi\1(\xi) \subseteq B_q\}$.
Using compactness, we assume that $S_q$ (and hence $S'_q$) is definable uniformly in $q$;
thus the union $S' := \bigcup_{q \in Q} S'_q$ is $0$-dimensional by Lemma~\ref{lem:dim} (\ref{it:dim-u})
and hence finite by Remark~\ref{rem:0fin}.

Using that $B_q$ is a union of fibers of the map $x \mapsto \rvi(x - S_q)$ (by Lemma~\ref{lem:fin&iso}),
one obtains $S'_q \ne \emptyset$ for every $q \in Q$.
Choose an element $x_0 \in S'$ such that
$Q' := \{q \in Q \mid x_0 \in S'_q \}$ is co-final in $Q$ (with respect to inclusion of the corresponding balls).
Now consider any $q' \in Q'$ and choose $\xi \in \RV$ with
$C := x_0  + \rvi\1(\xi) \subseteq B_{q'}$. Then every ball strictly containing $C$ contains $x_0$,
hence in particular $x_0 \in B_{q}$ for every $q \in Q$ with $B_{q} \supsetneq B_{q'}$.
Since $Q'$ is co-final in $Q$, we obtain $x_0 \in \bigcap_{q \in Q} B_q$.
\end{proof}

It will be important to us is that there is no risometry from a ball $B$
to a proper subset of $B$. In spherically complete valued fields, this is true in general.
We will require it only for definable risometries; to obtain that, definable spherically completeness
is be enough. The proof goes via the following ``definable Banach fixed point theorem''.

\begin{lem}\label{lem:banach}
Let $B \subseteq K^n$ be a ball and
suppose that $f\colon B \to B$ is definable and
contracting in the sense that
for any $x_1, x_2 \in B$ with $x_1\not=x_2$, $\val(f(x_1) - f(x_2)) > \val(x_1 - x_2)  $.
Then $f$ has (exactly) one fixed point.
\end{lem}

\begin{proof}
Suppose that $f(x) \ne x$ for all $x \in B$. For $x \in B$, set
\[
B_x := \ball{x,  \ge \val(x - f(x))}
.
\]
For two different points $x,x' \in B$, the assumption $\val(f(x) - f(x')) > \val(x - x')$
implies $\val(x - x') \ge \min\{\val(x - f(x)), \val(x' - f(x'))\}$ and hence either
$B_x$ contains $x'$ or vice versa. In particular, $B_x \cap B_{x'} \ne \emptyset$,
so all balls $B_x$ form a chain under inclusion. By Lemma~\ref{lem:dsc} (applied to each coordinate
projection), their intersection $\bigcap_{x \in B} B_x$ contains an element $x_0$.
Then $f(x_0) = x_0$, since otherwise, the assumption
$\val(f(x_0) - f(f(x_0))) > \val(x_0 - f(x_0))$ implies $x_0 \notin B_{f(x_0)}$.
\end{proof}

\begin{lem}\label{lem:isu}
Let $B \subseteq K^n$ be a ball and let $f\colon B \to X$ be a definable
risometry with $X \subseteq B$. Then $X = B$.
\end{lem}
\begin{proof}
Let $x_0 \in B$ be given; the idea is to find a preimage of $x_0$ by Newton-approximation
(although $f$ might not be differentiable, it behaves as if the derivative would be approximately 1).
For $x \in B$, set $g(x) := x + x_0 - f(x)$. Obviously, a fixed point of $g$ is a preimage of
$x_0$, so we just need to verify that $g$ is contracting. Indeed,
\[
g(x_1) - g(x_2) = (x_1 - x_2) - (f(x_1) - f(x_2)),
\]
and since $f$ is a risometry, we have $\rv(x_1 - x_2) = \rv(f(x_1) - f(x_2))$ and thus $\val(g(x_1) - g(x_2)) > \val(x_1 - x_2)$.
\end{proof}

\section{t-stratifications}
\label{sect:t-strat}

In this section, we will make the definition of t-stratification precise and we will prove
a bunch of basic properties. Throughout the section, we assume that $\Tx$
is a theory satisfying Hypothesis~\ref{hyp:general}$_0$ and that $K$ is a model of $\Tx$.
We start by looking more closely at the notion of translatability.

\subsection{Translatability}
\label{subsect:trans}

Recall that a lift of a sub-space $V \subseteq k^n$ is any sub-space $\tilde{V} \subseteq K^n$ with $\res(\tilde{V}) = V$.

\begin{defn}\label{defn:trans}
Suppose that $B_0$ is any definable set, $\chi\colon B_0 \to \RV\eq$ is a definable map, and
$B \subseteq B_0$ is a ball (open or closed).
\begin{enumerate}
\item\label{it:tr-inv}
For a sub-space $\tilde{V} \subseteq K^n$, we say that $\chi$ is
\emph{$\tilde{V}$-translation invariant} on $B$ if
for any $x, x' \in B$ with $x - x' \in \tilde{V}$, we have
$\chi(x) = \chi(x')$.
\item\label{it:V-tr}
For a sub-space $V \subseteq k^n$,
we say that $\chi$ is \emph{$V$-translatable} on $B$ if
there exists a lift $\tilde{V} \subseteq K^n$ of $V$ and a
definable risometry $\phi\colon B \to B$ such that
$\chi \circ \phi$ is $\tilde{V}$-translation invariant on $B$;
$\phi$ will be called a \emph{straightener} (of $\chi$ on $B$).
\item\label{it:d-tr}
For an integer $d \in \{0, \dots, n\}$, we say that $\chi$ is
\emph{$d$-translatable} on $B$ if there exists a $d$-dimensional
$V \subseteq k^n$ such that $\chi$ is $V$-translatable on $B$.
\end{enumerate}
\end{defn}

By Remark~\ref{rem:risoU}, in (\ref{it:V-tr}) the choice of $\tilde{V}$ doesn't matter, i.e.,
if $\chi$ is $V$-translatable, then for any lift $\tilde{V}$
of $V$, we can find a straightener $\phi$ such that $\chi \circ \phi$ is $\tilde{V}$-translation invariant.

One can easily modify Definition~\ref{defn:trans} to obtain a notion of translatability (on a ball $B$)
for definable sets $X \subseteq K^n$. More generally, we will use the following convention.

\begin{conv}\label{conv:trans}
Several Definitions of properties $P$ of definable maps
$\chi\colon B_0 \to \RV\eq$ for $B_0 \subseteq K^n$ (like Definition~\ref{defn:trans}) will also be applied
to subsets of $B_0$ and to tuples of maps and sets, by first turning such an object into a map, as follows.
\begin{enumerate}
\item A definable set $X \subseteq B_0$ has property $P$ iff the map $\chi \colon B_0 \to \RV\eq$
sending $X$ to $0 \in k$ and $B_0 \setminus X$ to $1 \in k$ has property $P$.
\item If $\theta = (\theta_1, \dots, \theta_\ell)$ is a tuple of maps and sets,
then we first replace each set $\theta_i$ by the corresponding map as in (1) and then
consider the map $\chi\colon x \mapsto (\theta_1(x), \dots, \theta_\ell(x))$;
$\theta$ has property $P$ iff $\chi$ has it.
\end{enumerate}
\end{conv}
Note that with this convention, the notion of translatability of a tuple of sets
agrees with the one introduced before Theorem~\ref{thm:weakA}.

It is clear that if a map $\chi$ is $V$-translatable on a ball $B$ for some $V \subseteq K^n$,
then it is also $V'$-translatable on $B$ for any sub-space $V' \subseteq V$.
Also, since risometries preserve balls, $V$-translatability on $B$
implies $V$-translatability on $B'$ for any sub-ball $B' \subseteq B$.
The following fact is less obvious than it looks. (Its proof needs definability of
the involved maps.)

\begin{lem}
Suppose that a definable map $\chi\colon B \to \RV\eq$ is both, $V_1$ and $V_2$-translatable
for some $V_1, V_2 \subseteq k^n$. Then $\chi$ is $(V_1 + V_2)$-translatable.
\end{lem}

\begin{proof}
Without loss, $V_1 \cap V_2 = 0$. Choose a complement $V_3$ of $V_1$ in $k^n$ with $V_2 \subseteq V_3$,
and let $\tilde V_i$ be lifts of $V_i$ (for $i=1,2,3$) with $\tilde V_2 \subseteq \tilde V_3$.
For $i = 1, 3$, let $\pi_i\colon K^n \to K^n$ be the projection to $\tilde V_i$ with kernel $\tilde V_{4-i}$.

Without loss, $\chi$ is $\tilde{V}_1$-translation invariant and $0 \in B$.
Let $\phi\colon B \to B$ be a definable risometry such that $\chi \circ \phi$ is $\tilde{V}_2$-translation invariant
and define $\psi\colon B \to B, x \mapsto \phi(\pi_3(x)) + \pi_1(x)$.
By Lemma~\ref{lem:risoKomb}, $\psi$ is a risometry. Its image is contained in $B$,
so by Lemma~\ref{lem:isu}, it is equal to $B$.
We claim
that $\chi \circ \psi$ is $(\tilde V_1+\tilde V_2)$-translation invariant.
Indeed, suppose that $x, x' \in B$ satisfy $x - x' \in \tilde V_1+\tilde V_2$
(or, equivalently, $\pi_3(x - x') \in \tilde V_2$).
Then $\chi(\psi(x)) = \chi(\phi(\pi_3(x)) + \pi_1(x)) = \chi(\phi(\pi_3(x)))
= \chi(\phi(\pi_3(x'))) = \chi(\phi(\pi_3(x')) + \pi_1(x')) = \chi(\psi(x'))$.
\end{proof}

By this lemma, for every definable map $\chi\colon  B_0 \to \RV\eq$ and every ball $B \subseteq B_0$, there exists a (unique) maximal space in which $\chi$ is translatable on $B_0$; we fix a notation for it.

\begin{defn}
Let $\chi\colon B_0 \to \RV\eq$ be a definable map for some $B_0 \subseteq K^n$ and let $B \subseteq B_0$ be a ball.
We write $\tsp_B(\chi)$ for the maximal sub-vector space of $k^n$
such that $\chi$ is $\tsp_B(\chi)$-translatable on $B$ (``$\tsp$'' stands for ``translatability space'').
Using Convention~\ref{conv:trans}, we also define $\tsp_B(\theta)$ when $\theta$ is a set or
a tuple of sets and maps.
\end{defn}

Using this definition, we have: $\chi$ is $V$-translatable on $B$ iff $V \subseteq \tsp_B(\chi)$,
and $\chi$ is $d$-translatable iff $\dim \tsp_B(\chi) \ge d$.

\medskip

To understand a $V$-translatable map $\chi$, we will often choose a projection $\pi\colon K^n \surject K^d$
and work fiberwise. This only works well if the fibers of $\pi$ are ``sufficiently transversal'' to lifts of $V$.
It will be handy to fix, once and for all, a finite set of $\emptyset$-definable projections
which work for all $V$. This is the purpose of the following definition.

\begin{defn}
Let $V \subseteq k^n$ be a sub-vector space.
An \emph{exhibition} of $V$ is a coordinate projection
$\pi\colon K^n \surject K^d$ inducing an isomorphism $\pibar \colon V \iso k^d$ (in particular, $d = \dim V$). We also say that $\pi$ \emph{exhibits} $V$.
If $B \subseteq K^n$ is a subset (usually a ball),
then the restriction $\pi\auf{B}$ will also be
called an exhibition of $V$.
\end{defn}

The following lemma summarizes basic facts needed to work fiberwise.

\begin{lem}\label{lem:tr-faser}
Suppose that $B \subseteq K^n$ is a ball, $V\subseteq k^n$ is a sub-vector space of dimension $d$, and
$\chi \colon B \to \RV\eq$ is a $V$-translatable definable map. Fix an exhibition $\pi\colon B \to K^d$ of $V$
and a lift $\tilde{V} \subseteq K^n$. Then we have the following.
\begin{enumerate}
\item
There exists a definable risometry $\phi\colon B \to B$
satisfying $\pi \circ \phi = \pi$
such that $\chi \circ \phi$ is $\tilde{V}$-translation invariant.
(In other words, $\phi$ is a straightener respecting the fibers of $\pi$.)
\item
For any definable risometry $\psi\colon B \to B$ and any $\pi$-fiber $F = \pi\1(x) \subseteq B$ (where $x \in \pi(B)$), there exists
a definable risometry $\psi'\colon F \to F$ such that
$(\chi \circ \psi)\auf{F} = (\chi\auf{F}) \circ \psi'$.
\end{enumerate}
\end{lem}

In (2), one can think of $\chi$ and $\chi\circ \psi$ as
two different but risometric maps; from that point of view, the statement is that
the restrictions of risometric maps to a $\pi$-fiber are also risometric.

\begin{proof}[Proof of Lemma~\ref{lem:tr-faser}]
Let $\pi_1, \pi_2\colon K^n \to K^n$ be the projections with $\im \pi_1 = \ker \pi_2 = \tilde V$ and
$\ker \pi_1 = \im \pi_2 = \ker \pi$.

(1)
Let $\psi$ be a straightener of $\chi$ and consider
the map $x \mapsto \pi_1(x) + \pi_2(\psi\1(x))$. By Lemma~\ref{lem:risoKomb},
it is a risometry and by Lemma~\ref{lem:isu}, it goes from $B$ onto $B$.
Its inverse is the desired straightener $\phi$.

(2)
Let $\phi$ be a straightener of $\chi$ satisfying $\pi \circ \phi = \pi$ and
define $\phi'(x) := \pi_1(x) + \pi_2(\phi\1(\psi(x)))$.
By Lemmas~\ref{lem:risoKomb} and~\ref{lem:isu}, $\phi'$ is a risometry from $B$ to $B$
and by definition, we have $\pi \circ \phi' = \pi$ and (using that $\chi \circ \phi$ factors over $\pi_2$)
$(\chi \circ \phi) \circ \phi\1 \circ \psi = (\chi \circ \phi) \circ \phi'$, so we can define $\psi'$ to be the restriction
of $\phi \circ \phi'$ to $F$.
\end{proof}

Using this, we can give an alternative characterization of translatability.
Recall that for a ball $B$, $B - B$ is the ball of the same radius containing
the origin and that $\dir$ was introduced in Definition~\ref{defn:dir}.

\begin{lem}\label{lem:translater}
Let $\chi \colon B \to \RV\eq$ be a definable map, $V \subseteq k^n$ a sub-space and $\pi\colon K^n \to K^d$
an exhibition of $V$. Then $\chi$ is $V$-translatable if and only if there exists
a definable family of risometries $\alpha_x \colon B \to B$, where
$x$ runs over $\pi(B - B)$, with the following properties (for all $x, x' \in \pi(B - B)$ and all $z \in B$):
\begin{enumerate}
\item\label{it:chi} $\chi \circ \alpha_x = \chi$;
\item\label{it:circ} $\alpha_x \circ \alpha_{x'} = \alpha_{x + x'}$;
\item\label{it:pi} $\pi(\alpha_x(z) - z) = x$;
\item\label{it:dir} $\dir(\alpha_x(z) - z) \in V$ if $x \ne 0$.
\end{enumerate}
\end{lem}

\begin{proof}
``$\Rightarrow$'': Choose a straightener $\phi$ respecting the fibers of $\pi$
(using Lemma~\ref{lem:tr-faser}~(1)) and let
$\tilde{V}$ be the corresponding lift of $V$. For any $x \in \pi(B - B)$, denote by $\alpha'_x\colon B \to B$
the translation by the unique element of $\pi\1(x) \cap \tilde{V}$.
Then $\alpha'_x$ satisfies $\chi \circ \phi \circ \alpha'_x = \chi\circ \phi$ and (2)~--~(4), and from
this, one deduces that $\alpha_x := \phi \circ \alpha'_x \circ \phi\1$ satisfies (1)~--~(4).

``$\Leftarrow$'': Without loss, $0 \in B$, $V = k^d \times \{0\}^{n-d}$, and
$\pi$ is the projection to the first $d$ coordinates.
Write elements of $B$ as $(x, y) \in K^d \times K^{n-d}$.
We claim that $\phi(x,y) := \alpha_{x}(0, y)$ is a straightener.

By (1), $\chi \circ \phi$ is $(K^d \times \{0\}^{n-d})$-translation invariant.
To check that $\phi$ is a risometry, consider $(x_1, y_1), (x_2, y_2) \in B$ and set $x := x_2 - x_1$.
We have $\phi(x_1,y_1) = \alpha_{x_1}(0, y_1)$ and
(2) implies $\phi(x_2,y_2) = \alpha_{x_1}(\alpha_{x}(0, y_2))$, so since
$\alpha_{x_1}$ is a risometry, it suffices to check that
$\rv((x_1, y_1) - (x_2, y_2)) = \rv((0, y_1) - \alpha_{x}(0, y_2))$;
but this follows from $\val(y_2 - \picomp(\alpha_x(0, y_2))) > \val(x)$, which in turn follows
from (3) and (4). (Recall that $\picomp\colon B \to K^{n-d}$ is the ``complement'' of $\pi$.)
\end{proof}

\begin{defn}\label{defn:translater}
Let $B \subseteq K^n$ be a ball, $\pi\colon B \to K^d$ a coordinate projection,
and $\chi\colon B \to \RV\eq$ a definable map.
A definable family of risometries $(\alpha_x)_{x \in \pi(B - B)}$ from $B$ to itself
satisfying (1)~--~(4) of Lemma~\ref{lem:translater}
will be called a \emph{translater} of $\chi$ (on $B$, with respect to $\pi$).
We also apply Convention~\ref{conv:trans}.
\end{defn}

Characterizing translatability via translaters has the disadvantage of being more technical, but one advantage is that it avoids the (uncanonical) lift $\tilde{V}$
appearing in Definition~\ref{defn:trans}.

The following lemma says how (and when) translatability is preserved under the
restriction to an affine subspace.

\begin{lem}\label{lem:trans&fiber}
Suppose that $\chi\colon B \to \RV\eq$ is a $V$-translatable definable map (where $B \subseteq K^n$ is a ball
and $V \subseteq k^n$) and $\rho\colon B \to K^{d}$ is a coordinate projection
with $\rhobar(V) = k^{d}$ (in particular $\dim V \ge d$). Then the restriction of
$\chi$ to any fiber $\rho\1(y)$ (for $y \in \rho(B)$)
is $(V \cap \ker\rhobar)$-translatable.
\end{lem}

\begin{proof}
Choose an exhibition $\pi\colon B \to K^{d'}$ of $V$ satisfying $\ker \pi \subseteq \ker \rho$
and let $\phi \colon B \to B$ be a straightener of $\chi$ respecting the fibers of $\pi$.
Then $\phi$ sends any $\rho$-fiber $\rho\1(y)$ to itself and thus
$\phi\auf{\rho\1(y)}$ is a straightener of $\chi\auf{\rho\1(y)}$ proving 
$(V \cap \ker\rho)$-translatability.
\end{proof}

The next two lemmas state that translatability behaves as one would expect with respect
to dimension and topological closure (using the valued field topology);
we write $\topcl{X}$ for the topological closure of a set $X$.

\begin{lem}\label{lem:I&D}
Suppose that $B \subseteq K^n$ is a ball, that
$X \subseteq B$ is a definable set which is $V$-translatable on $B$
for some $V \subseteq k^n$, and that $\pi\colon B \to K^d$ exhibits $V$.
Then for any $x \in \pi(B)$, we have
\[
\dim X = \dim (X \cap \pi\1(x)) + d
.
\]
In particular, $\dim X \ge \dim V$.
\end{lem}

\begin{proof}
The translaters of Lemma~\ref{lem:translater} can be restricted to definable bijections
between the fibers $X \cap \pi\1(x)$, so all of them have the same dimension.
Now use Lemma~\ref{lem:dim}~(\ref{it:dim-f}).
\end{proof}

\begin{lem}\label{lem:I&T}
If $X \subseteq K^n$ is $V$-translatable on a ball $B\subseteq K^n$, then so is $(X, \topcl{X})$.
\end{lem}
\begin{proof}
Since risometries are homeomorphisms, a straightener for $X$ also straightens
$\topcl{X}$.
\end{proof}

\subsection{Definition of t-stratifications}

We now give the general definition of a t-stratification and prove basic properties.
(The ``t'' in ``t-stratification'' stands for ``translatable''.)
Recall that translatability has been defined precisely in Definition~\ref{defn:trans}
and Convention~\ref{conv:trans}.

\begin{defn}\label{defn:t-strat}
Let $B_0 \subseteq K^n$ be a ball.
A \emph{$t$-stratification} of $B_0$ is a partition of $B_0$
into definable sets $S_0, \dots, S_n$ with the properties listed below.
We write $S_{\le d}$ for $S_0 \cup \dots \cup S_d$ and
$S_{\ge d}$ for $S_d \cup \dots \cup S_n$.
\begin{enumerate}
\item $\dim S_d \le d$
\item For each $d$ and each ball $B \subseteq S_{\ge d}$ (open or closed),
$(S_i)_{i \le n}$ is $d$-translatable on $B$.
\end{enumerate}
We say that a t-stratification $(S_i)_{i \le n}$ \emph{reflects}
a definable map $\chi\colon B_0 \to \RV\eq$ if the following stronger version of (2) holds.
\begin{enumerate}
\item[(2')] For each $d$ and each ball $B \subseteq S_{\ge d}$ (open or closed),
$((S_i)_{i \le n}, \chi)$ is $d$-translatable on $B$.
\end{enumerate}
We define when a t-stratification reflects a set or a tuple of sets and maps
using Convention~\ref{conv:trans}.
\end{defn}

In other words, for any ball $B \subseteq B_0$, we let $d$ be minimal with $B \cap S_d \ne \emptyset$
and require $d$-translatability on $B$. Note that this is as much as one can get:
since $\dim S_d \le d$, Lemma~\ref{lem:I&D} implies that $(S_i)_{i}$ is not $(d+1)$-translatable on $B$.
In particular, we have $\tsp_B(S_d) = \tsp_B((S_i)_{i}) = \tsp_B((S_i)_{i}, \chi)$.

In general, if, for some $V \subseteq k^n$, two maps $\chi$ and $\chi'$ are both $V$-translatable, this does not imply that
$(\chi, \chi')$ is also $V$-translatable. However, we will see in Remark~\ref{rem:reflect} that
if a t-stratification $(S_i)_{i}$ reflects both $\chi$ and $\chi'$, then it also reflects $(\chi, \chi')$.

\begin{rem}\label{rem:ball->all}
If $(S_i)_i$ is a t-stratification of $B_0$ (reflecting $\chi$), then
the restriction to any subball of $B_0$ is also a t-stratification (reflecting
the restriction of $\chi$). In the other direction, a t-stratification
of $B_0 \subseteq K^n$ can be extended to a t-stratification of $K^n$
by replacing $S_n$ with $S_n \cup (K^n \setminus B_0)$, but only under the assumption
that $S_0 \ne \emptyset$. This assumption is needed because in general,
$(S_i)_i$ will not be $1$-translatable on any ball strictly bigger than $B_0$.
\end{rem}

In general, even if $\chi\colon K^n \to \RV\eq$ is definable, the map $B \mapsto \tsp_B(\chi)$
does not need to be definable (see Example~\ref{ex:tsp-ndef} below). However, for t-stratifications,
the corresponding map is definable, as the following lemma states.

\begin{lem}\label{lem:trans-def}
For a fixed t-stratification $(S_i)_i$ of $K^n$, the map
$B \mapsto \tsp_B((S_i)_i)$ is definable, uniformly for all models
$K \models \Tx$ (i.e., if each $S_i$ is given by a formula, then
then there exists a formula defining $B \mapsto \tsp_B((S_i)_i)$
in all models).
\end{lem}

\begin{proof}
The dimension of $\tsp_B((S_i)_i)$ can be defined as the minimal
$d$ with $B \cap S_d \ne \emptyset$. We claim that then, for any
$x \in S_d \cap B$ and any sufficiently small ball $B'$ containing $x$, we have
$\tsp_B((S_i)_i) = \{\dir(x_1 - x_2) \mid x_1, x_2 \in S_d \cap B'\}$.

To prove this claim, fix an exhibition $\pi\colon B \to K^d$ of $W := \tsp_B((S_i)_i)$.
For each $\pi$-fiber $F = \pi\1(y)$ (with $y \in \pi(B)$),
$S_d \cap F$ is zero-dimensional
by Lemma~\ref{lem:I&D} and hence finite by Remark~\ref{rem:0fin}, so for any $x \in S_d \cap F$, we can find
a ball $B' \subseteq B$ such that $B' \cap S_d \cap F = \{x\}$.
By $W$-translatability on $B'$, $B' \cap S_d \cap \pi\1(y')$
is a singleton for any $y' \in \pi(B')$ and we obtain the claim.
\end{proof}

The following is an example of a set $X \subseteq K^3$ such that $\tsp_B(X)$ does not depend definably on $B$. The key ingredient
is that if the residue field $k$ is ``sufficiently evil'', then whether a definable map $k \to k$ can be lifted
to a definable map $\valring \to \valring$ is not a definable condition.

\begin{ex}\label{ex:tsp-ndef}
Recall that in the field $\bbQ$, the subset $\bbN$ and the exponential map $\bbQ \times \bbN \to \bbQ, (x, y) \mapsto x^y$
are definable. Fix an elementary extension $\nsQ \succ \bbQ$ and write $\nsN$ for the interpretation of $\bbN$
in $\nsQ$. Set $K := \nsQ((t))$ and
define 
\[
X := \{(t\1x, y, tz) \in t\1\valring \times \valring \times t\valring \mid \res(x) \in \nsN \wedge \res(z) = \res(y)^{\res(x)}\}.
\]
We claim that for $a \in \nsN$, $X$ is $(\nsQ^2\times\{0\})$-translatable on the ball $B_a := t\1\res\1(a) \times \valring^2$ iff $a \in \bbN$, which of course is not a definable set.

Fix $a \in \nsN$.
Since on $B_a$, $X$ is $(K \times \{0\}^2)$-translation invariant, the question is only whether
the fiber $X_a = \{(y, tz) \in \valring \times t\valring \mid \res(z) = \res(y)^a\}$
is $\nsQ \times \{0\}$-translatable on $\valring^2$.
If $a \in \bbN$, then an easy computation shows that the map $\valring^2 \to \valring^2, (y,z) \mapsto (y,z+ty^a)$
is a straightener of $X_a$.

Now suppose that $X_a$ is $\nsQ\times\{0\}$-translatable on $\valring^2$ for some $a \in \nsN \setminus \bbN$.
Then we can find a $1$-dimensional definable subset $Y \subseteq X_a$ whose projection to the first coordinate is
equal to $\valring$. (Choose a straightener $\phi\colon \valring^2 \to \valring^2$ and set $Y:= \phi(\valring \times \{z_0\})$ for a suitable $z_0 \in \valring$.)
Set $Y' := \{(y,z) \in \valring^2 \mid (y, tz) \in Y\}$.
Since $\dim Y' = 1$, there exists a polynomial $f \in \valring[y,z] \setminus \maxid[y,z]$ vanishing on $Y'$
(e.g.\ by Lemma~\ref{lem:alg-dim}), and thus $\bar f := \res(f)$ is a non-zero polynomial vanishing on
$\res(Y') = \{(y, y^a) \mid y \in \nsQ\}$.
However, for any $d \in \bbN$, we have
\[
\bbQ \models \forall k \in \bbN, k > d\colon \text{no polynomial of degree $\le d$ vanishes on }
\{(y, y^k) \mid y \in \bbQ\}.
\]
Since this also holds in $\nsQ$ and since 
$a > d$ for all $d \in \bbN$, we obtain a contradiction.
\end{ex}

A property of t-stratifications which is important
for inductive arguments is that on an affine subspace of a ball which is transversal to the translatability space on that ball, 
they again induce
t-stratifications. This is the statement of the following lemma.
(It is formulated for a t-stratification reflecting a
map, but of course, we can apply it to the trivial map if we
are interested in a ``pure'' t-stratification.)

\begin{lem}\label{lem:t-strat-ball}
Let $(S_i)_i$ be a t-stratification of $B_0 \subseteq K^n$ reflecting
a definable map $\chi\colon B_0 \to \RV\eq$, and let $B \subseteq B_0$ be a ball.
Let $\pi\colon B \to K^d$ be an exhibition of $\tsp_B((S_i)_i)$ and
suppose that $F = \pi\1(x)$ is a $\pi$-fiber (for some $x \in \pi(B)$).
Set $T_i := S_{i+d} \cap F$ for $0 \le i \le n - d$.
Then $(T_i)_{i \le n - d}$ is a t-stratification of $F \cap B$
reflecting $\chi\auf{F \cap B}$.
\end{lem}

\begin{proof}
By Lemma~\ref{lem:I&D}, $\dim S_{i+d} \le i + d$ implies $\dim T_i \le i$,
so it remains to show the translatability condition.
Consider a ball $B' \subseteq B$ with $B' \cap F \ne \emptyset$ and suppose that
$j$ is minimal with $B' \cap F \subseteq T_{\ge j}$. We have to show $j$-translatability on
$B' \cap F$.

Set $V := \tsp_B((S_i)_i)$ and $V' := \tsp_{B'}((S_i)_i)$.
By $V$-translatability of $S_{\ge j + d}$, $B' \cap F \subseteq S_{\ge j + d}$ implies
$B' \subseteq S_{\ge j + d}$, so $\dim V' = j+d$.
Since $V \subseteq V'$, we have $\pibar(V') = k^d$, so Lemma~\ref{lem:trans&fiber}
implies $(V' \cap \ker\pibar)$-translatability of $((S_i)_i, \chi)$ on $B' \cap F$.
Now we are done since $\dim (V' \cap \ker\pibar) = j$.
\end{proof}

Here are some ``global'' properties of t-stratifications.

\begin{lem}\label{lem:t-strat-glob}
Let $(S_i)_i$ be a t-stratification of $B_0 \subseteq K^n$. Then the following holds.
\begin{enumerate}
\item\label{it:max-open}
For each $d$ and each $x \in S_{\ge d+1}$, there exists a maximal ball
$B$ containing $x$ such that $B \cap S_{\le d} = \emptyset$. Moreover,
if $B \ne B_0$ then $B$ is open.
In particular, the sets $S_{\le d}$ are topologically closed.
\item\label{it:loc-dim}
$S_d$ has dimension exactly $d$ locally at each point $x \in S_d$.
In particular, either $\dim S_d = d$ or $S_d = \emptyset$.
\end{enumerate}
\end{lem}

\begin{proof}
(\ref{it:max-open})
For $d = 0$, this is clear since $S_0$ is finite; now suppose $d > 0$.
By induction, there is a maximal ball $B$ containing $x$
with $B \cap S_{\le d - 1} = \emptyset$. If $B \cap S_d = \emptyset$, then
$B$ is the ball we are looking for, so suppose now that $B \cap S_d \ne \emptyset$.
Then $V := \tsp_B((S_i)_i)$ is $d$-dimensional; let $\pi\colon B \to K^d$ be an exhibition
of $V$ and let $F \subseteq B$ be the $\pi$-fiber containing $x$.
Since $F \cap S_d$ is finite and non-empty, we find a maximal open ball $B' \subseteq B$
such that $B' \cap F \cap S_d = \emptyset$.
Now $V$-translatability implies $B' \cap S_d = \emptyset$, so $B'$ is the ball we were looking for.

(\ref{it:loc-dim})
Let $x \in S_d$ be given. By (\ref{it:max-open}), there exists a ball $B$ containing $x$ with $B \subseteq S_{\ge d}$, hence on any sub-ball $B' \subseteq B$,
we have $d$-translatability. Now $\dim (S_d \cap B') < d$ would contradict
Lemma~\ref{lem:I&D}.
\end{proof}

\subsection{Families of sets up to risometry}
\label{subsect:sak}

Given a definable family $(\chi_q)_{q \in Q}$ of maps $K^n \to \RV\eq$,
whether two maps $\chi_q$, $\chi_{q'}$ are definably risometric defines an equivalence relation on $Q$.
This equivalence relation is in general not definable. (For the family
$(X_{a})_{a \in \nss\bbN}$ from Example~\ref{ex:tsp-ndef},
$\bbN$ is one of the equivalence classes.)
The main result of this subsection is that it does become definable if we equip each $\chi_q$ with a t-stratification. Moreover, for each equivalence class, we
can find a definable family of risometries which are compatible under composition.

Under an additional assumption, we can even get some more information about
these risometries; this will be needed in the proof of the main theorem.
To formulate that assumption we temporarily introduce
the following, very weak variant of translatability; this definition will only be used
in this subsection and in the application in the proof of Lemma~\ref{lem:ZP}.

\begin{defn}\label{defn:pt-trans}
Suppose that $B \subseteq K^n$ is a definable subset (usually a ball), $\chi\colon B \to \RV\eq$ is a definable map,
and $V \subseteq k^n$ is a vector space exhibited by $\pi\colon B \to K^d$.
We say that $\chi$ is \emph{pointwise translatable} on $B$ in direction $V$
with respect to $\pi$ (or simply \emph{$V$-$\pi$-pointwise translatable})
if for any $y \in B$ and any $x' \in \pi(B)$,
there exists an $y' \in \pi\1(x')$ with
$\chi(y') = \chi(y)$ and
$\dir(y - y') \in V$.
We also use Convention~\ref{conv:trans}.
\end{defn}

Notice the similarity of pointwise translatability to Condition~(\ref{it:dir}) in Lemma~\ref{lem:translater} (the definition of translater).

\begin{prop}\label{prop:sak}
Suppose that $Q$ is a $\emptyset$-definable set (in any sort),
$(S_{i})_{i \le n}$ is a $\emptyset$-definable partition of $Q \times K^n$
and $\chi\colon Q \times K^n \to \RV\eq$ is a $\emptyset$-definable map. Write $\pi$ for
the projection $Q \times K^n \surject Q$.
For $q \in Q$, set $S_{i,q} := S_i \cap \pi\1(q)$ and
$\chi_q := \chi\auf{\pi\1(q)}$. Then we have the following.
\begin{enumerate}
\item
The set $Q' \subseteq Q$ of those $q$ for which 
$(S_{i,q})_{i \le n}$ is a t-stratification of $\{q\} \times K^n$ reflecting $\chi_q$ is $\emptyset$-definable.
\item
There exists a $\emptyset$-definable map $\chi'\colon Q' \to \RV\eq$ such that
for all $q, q' \in Q'$,
$\chi'(q) = \chi'(q')$ if and only if there exists a
definable risometry $\phi\colon \{q\} \times K^n \to \{q'\} \times K^n$ respecting $((S_{i})_{i}, \chi)$.
\item
For each $\chi'$-fiber $C \subseteq Q'$, there exists a
compatible $\code{C}$-definable family $\alpha_{q,q'}\colon ((S_{i,q})_{i}, \chi_{q}) \to ((S_{i,q'})_{i}, \chi_{q'})$ of risometries, where $q,q'$ run through $C$.
Compatible means: $\alpha_{q',q''} \circ \alpha_{q,q'} = \alpha_{q,q''}$ for $q,q', q'' \in C$.
\end{enumerate}
Suppose now that $Q \subseteq K^m$ and that $V \subseteq k^{m+n}$ is exhibited
by the projection $\pi\colon Q \times K^n \surject Q$. Then we also have the following variant of (3):
\begin{enumerate}
\item[(3')]
If, for some $\chi'$-fiber $C \subseteq Q'$,
$(S_i)_i$ is $V$-$\pi$-pointwise translatable on $C \times K^n$
and moreover $S_0 \cap (C \times K^n) \ne \emptyset$, then
the family $\alpha_{q,q'}$ can be chosen such that additionally,
$\dir(\alpha_{q,q'}(x) - x) \in V$ for all $q, q' \in C$ and all $x \in \{q\} \times K^n$.
\end{enumerate}
All of the above works uniformly for all models $K$ of our theory $\Tx$, i.e.,
given formulas defining $Q$, $S_i$ and $\chi$, we can find formulas defining
$Q'$, $\chi'$ and $\alpha_{q,q'}$ not depending on $K$.
\end{prop}

\begin{rem}
By taking $Q := \{0, 1\}$, in particular we obtain:
if there exists a definable risometry between two $\emptyset$-definable t-stratifications, then there already exists a $\emptyset$-definable one.
\end{rem}

Before we prove the proposition, let us consider the following corollary, which
shows how statement (3') can be used to deduce translatability.

\begin{cor}\label{cor:sak}
Suppose that $B \subseteq K^n$ is a ball, $\chi\colon B \to \RV\eq$ is a definable map,
$V \subseteq k^n$ is a sub-space exhibited by $\pi\colon B \to K^d$, and $(S_i)_{0 \le i \le n-d}$ is a definable partition of $B$ such that for each $\pi$-fiber $F \subseteq B$,
$(S_i \cap F)_i$ is a t-stratification reflecting $\chi\auf{F}$.
Suppose moreover that $S_0$ is non-empty, that for any two $\pi$-fibers $F, F'$
there exists a definable risometry $\phi\colon F \to F'$ respecting $((S_i)_i, \chi)$, and that
$(S_i)_i$ is $V$-$\pi$-pointwise translatable on $B$. Then $((S_i)_i, \chi)$ is $V$-translatable on $B$.
\end{cor}

\begin{proof}
Without loss, $\pi$ is the projection to the first $d$ coordinates;
set $Q := \pi(B)$. We extend the domains of $\pi$, $\chi$, and $(S_i)_i$
from $B$ to $Q \times K^{n-d}$ as follows. For $\chi$, we send all of $(Q \times K^{n-d}) \setminus B$
to a single new element in $\RV\eq \setminus \chi(B)$,
and for $(S_i)_i$, we simply enlarge $S_{n-d}$ (and keep all $S_i$ for $i < n - d$).
Then for each $q \in Q$, $(S_i \cap \pi\1(q))_i$ is a t-stratification of $\pi\1(q)$ reflecting
$\chi\auf{\pi\1(q)}$ (this uses $S_0 \cap \pi\1(q) \ne \emptyset$; cf.\ Remark~\ref{rem:ball->all}).
When applying Proposition~\ref{prop:sak} to this data, the map $\chi'$ we obtain is constant on $Q$
and $Q$ satisfies the prerequisites of (3'), hence we obtain a
family of risometries $\alpha_{q_1, q_2}\colon \pi\1(q_1) \to \pi\1(q_2)$
as in (3') for $q_1, q_2 \in Q$. Define a family $(\beta_q)_{q \in Q-Q}$ of maps $B \to B$ by 
$\beta_q(x) := \alpha_{\pi(x), \pi(x)+q}(x)$.
We claim that this family is a translater
proving $V$-translatability of $((S_{i})_{i}, \chi)$ on $B$.

It is clear that these $\beta_q$ satisfy Conditions (1)~--~(4) of Lemma~\ref{lem:translater}
(by definition of $\beta_q$ and by the properties of $\alpha_{q_1,q_2}$),
so it remains to check that each $\beta_q$ is a risometry.
To see this, choose $x_1, x_2 \in B$ and set $q_i := \pi(x_i)$
and $x_3 := \alpha_{q_2,q_1}(x_2)$. Then
$\beta_{q}$ preserves both, $\rv(x_1 - x_3)$
(since $\alpha_{q_1,q_1+q}$ is a risometry) and $\rv(x_3 - x_2)$
(by Lemma~\ref{lem:vf}~(\ref{it:dir-pi}), since $\pi(x_3 - x_2)$ and $\dir(x_3 - x_2) = V$ are preserved), and
these two values together determine $\rv(x_1 - x_2)$ by Lemma~\ref{lem:vf}~(\ref{it:rv-sum}).
\end{proof}

\begin{proof}[Proof of Proposition~\ref{prop:sak}] 
The whole proof is by induction on $n$, i.e. we assume that the proposition
holds for smaller $n$.

\textbf{(1)}
Here is an informal formula defining $Q'$ (where $q$ is the variable):
\begin{align*}
&\bigwedge_{d = 0}^n \dim S_{d,q} \le d\\
&\llap{$\wedge\,\,\,$}\bigwedge_{d = 1}^n\,\forall\text{ balls $B\subseteq S_{\ge d,q}$ with }B \cap S_{d,q} \ne \emptyset:\\
&\quad\bigvee_{\substack{\rho\colon B \to \{q\} \times K^d\\\text{coordinate}\\\text{projection}}}\hskip-3ex
\exists \,V \subseteq k^n \text{ sub-space}:\\
&\qquad\qquad\qquad\rho \text{ exhibits } V \quad\\
&\qquad\qquad\qquad\llap{$\wedge\,\,\,$}(S_{i,q})_{i} \text{ is $V$-$\rho$-pointwise translatable on $B$}\\
&\text{{\small\emph{[For $x \in \rho(B)$, set $T_{i,x} := S_{i + d,q} \cap \rho\1(x)$ and $\chi_x := \chi_q\auf{\rho\1(x)}$]}}}\\
&\qquad\qquad\qquad\llap{$\wedge\,\,\,$}\text{$(T_{i,x})_{i \le n - d}$ is a t-stratification reflecting $\chi_x$ for all $x \in \rho(B)$}\\
&\qquad\qquad\qquad\llap{$\wedge\,\,\,$}\text{all $((T_{i,x})_i, \chi_x)$ are definably risometric for } x \in \rho(B)
\end{align*}
This is first order: in the first line, we use that 
dimension is definable; in the last two lines, we use
(1), (2) of the induction hypothesis.

If $(S_{i,q})_{i}$ is a t-stratification reflecting $\chi_q$, then it is clear the formula holds.
For the other direction, note that by the induction hypothesis and Corollary~\ref{cor:sak}, the last four lines of the formula
together with $B \cap S_{d,q} \ne \emptyset$
imply that $((S_{i,q})_{i}, \chi_q)$ is $V$-translatable on $B$.

\textbf{(2) and (3)}
Without loss, $Q = Q'$. Moreover, if we have a definable map $\chi' \colon Q \to \RV\eq$
such that the existence of a risometry $\{q\} \times K^n \to \{q'\} \times K^n$ respecting $((S_{i})_{i}, \chi)$ implies $\chi'(q) = \chi'(q')$, we can consider each $\chi'$-fiber
separately for the remainder of the proof (adding the image of the fiber to the language).
We will do this several times; at the end, we will obtain a
definable compatible family of risometries 
on the whole of $Q$, thus proving both (2) and (3).

By Lemma~\ref{lem:fin&iso}~(1), there is at most one risometry sending
$S_{0,q}$ to $S_{0,q'}$. Whether such a risometry exists can definably be
tested by choosing an enumeration $(x_\mu)_\mu$ of $S_{0,q}$
and comparing the matrix $(\rv(x_\mu - x_\nu))_{\mu,\nu}$ to a corresponding matrix for $S_{0,q'}$.
(Note that the cardinality $|S_{0,q}|$ is bounded.)
Thus we can suppose that
for each $q, q' \in Q$, a risometry $\beta_{q,q'}\colon S_{0,q} \to S_{0,q'}$ exists
and that $\beta_{q,q'}$ respects $\chi\auf{S_0}$.
Moreover (again by uniqueness of this risometry), $\beta_{q,q'}$ is $(q,q')$-definable
and the family $(\beta_{q,q'})_{q,q'}$ is compatible with composition (as required in (3)).

Consider a set $R \subseteq \RV^{(n)}$ such that
$B_{R,q} := \{x \in \{q\} \times K^n \mid \rv(x - S_{0,q}) = R\}$ is non-empty.
By Lemma~\ref{lem:fin&iso}, this non-emptiness condition does not depend on $q$, $B_{R,q}$
is a maximal ball not intersecting $S_{0,q}$ (possibly equal to $\{q\} \times K^n$), and any risometry
$\{q\} \times K^n \to \{q'\} \times K^n$ respecting $S_{0}$ sends
$B_{R,q}$ to $B_{R,q'}$.
This means that we can treat each family $(B_{R,q})_{q}$ separately as follows.
For each $R$ as above, we will construct an $\code{R}$-definable map $\chi'_R \colon Q \to \RV\eq$ and an
$\code{R}$-definable compatible family of risometries $\alpha_{R,q,q'}\colon B_{R,q} \to B_{R,q'}$ 
such that (2) and (3) hold for the restricted ($\code{R}$-definable) family
$((S_{i,q} \cap B_{R,q})_i, \chi_q\auf{B_{R,q}})_{q \in Q}$. By compactness,
we can assume that the definitions of $\chi'_R$ and $\alpha_{R,q,q'}$ are uniform in $\code{R}$.
Using stable embeddedness of $\RV$ (Hypothesis~\ref{hyp:general} (\ref{it:st-emb})),
we define the aggregate map $\chi'\colon Q \to \RV\eq, q \mapsto \code{\code{R} \mapsto \chi'_R(q)}$.
For $q, q' \in Q$ with $\chi'(q) = \chi'(q')$, the risometries $\beta_{q,q'}$ and
$\alpha_{R,q,q'}$ can be assembled to a $\emptyset$-definable map $\alpha_{q,q'}\colon K^n \to K^n$;
this map is a risometry by Lemma~\ref{lem:fin&iso}, and for varying $q, q'$
we have the required compatibility.

Thus from now on, we fix $R$, we add $\code{R}$ to the language, and to simplify notation, we write $B_{q}$ instead of $B_{R,q}$. Moreover, we set $B := \bigcup_{q \in Q} B_q \subseteq Q \times K^n$.

For some $q \in Q$, set $W := \tsp_{B_q}((S_{i,q})_{i})$. By Lemma~\ref{lem:trans-def}, $W$ is $q$-definable, so
using the map $q \mapsto \code{W} \in \RV\eq$, we may assume that $W$ does not depend on $q$.
Set $d := \dim W$, choose an exhibition $\rho'\colon K^n \surject K^d$ of $W$,
and set $\rho\colon B \to Q \times K^d, (q,y) \mapsto (q,\rho'(y))$.
Note that $d \ge 1$, since $B_q \cap S_{0,q} = \emptyset$.
By Lemma~\ref{lem:t-strat-ball}, we get a family of t-stratifications
of the fibers of $\rho$, parametrized by $\tilde{Q} := \rho(B)$.
Applying induction (2) to this family (after using Remark~\ref{rem:ball->all}) yields a map $\tilde{\chi} \colon \tilde{Q} \to \RV\eq$ depending only on
the $Q$-coordinate. By Lemma~\ref{lem:tr-faser}, if there exists a risometry $B_{q} \to B_{q'}$ respecting
$((S_i)_i, \chi)$, then there also exist such risometries between any two $\rho$-fibers contained in $B_{q}$ or $B_{q'}$, so we can assume
that $\tilde{\chi}$ is constant on $\tilde{Q}$. Now induction (3) yields a definable compatible family of risometries $\tilde{\alpha}_{\tilde{q},\tilde{q}'}\colon
\rho\1(\tilde q) \to \rho\1(\tilde q')$ for $\tilde q, \tilde q' \in \tilde Q$.

To finish the construction of a definable compatible family of risometries
$\alpha_{q,q'}\colon B_q \to B_{q'}$, it remains to find a definable compatible family
of risometries $\gamma_{q,q'}\colon \rho(B_q) \to \rho(B_{q'})$ (which does not need to respect anything); after that, we can set
$\alpha_{q,q'}(y) := \tilde{\alpha}_{\tilde{q}, \tilde{q}'}(y)$, where $\tilde{q} := \rho(y)$ and $\tilde{q}' := \gamma_{q,q'}(\tilde{q})$.

If $S_0$ is empty, then $B_q = \{q\} \times K^n$ and we can set $\gamma_{q,q'}(q,z) := (q',z)$ for
every $q,q' \in Q$, $z \in K^n$,
so suppose now that $S_0$ is non-empty.
For each $q$, let $N_q$ consist of those elements of $S_{0,q}$ that are closest to $B_q$,
i.e., $N_q = S_{0,q} \cap \bar{B}_q$, where $\bar{B}_q$ is the unique closed ball containing $B_q$ with $\radc(\bar{B}_q) = \rado(B_q)$.
Define $c_q := \frac1{|N_q|}\sum_{s\in N_q} s$ to be the barycenter of $N_q$.
The translation $\tilde{\gamma}_{q,q'}\colon y \mapsto y - c_q + c_{q'}$ sends $B_q$ to $B_{q'}$,
so we can define $\gamma_{q,q'}\colon \rho(B_q) \to \rho(B_{q'})$ to be
the induced translation on the projections.

\textbf{(3')}
Let us say that a map $\phi$ between two subsets of $Q \times K^n$ ``moves in direction $V$''
if $\dir(\phi(x) - x) \in V$ for all $x$ in the domain. We claim that under the additional
assumptions of (3'), the maps $\alpha_{q,q'}$ constructed in the proof of (3) do already
move in direction $V$; so let us go through the construction of $\alpha_{q,q'}$.

First, we have to check that the risometries $\beta_{q,q'}\colon S_{0,q} \to S_{0,q'}$
move in direction $V$. Set $\delta := \val(q - q')$, and let us say that
$T \subseteq S_{0,q'}$ is a set of $\delta$-representatives (of $S_{0,q'}$)
if for each $s \in S_{0,q'}$ there exists exactly one $t \in T$
with $\val(s - t) > \delta$. Choose any set of
$\delta$-representatives $T \subseteq S_{0,q'}$.
For each $t \in T$, using pointwise translatability of $S_0$,
we can choose an element $\phi(t) \in S_{0,q}$ with $\dir(\phi(t) - t) \in V$.
Using that $\val(t - t') \le \val(q - q')$ for any two different $t, t' \in T$,
we get that $\phi\colon T \to \phi(T)$ is a risometry.
Composing with $\beta_{q,q'}$ yields a risometry from $T$ to
$T' := \beta_{q,q'}(\phi(T))$,
and $T'$ is also a set of $\delta$-representatives of $S_{0,q'}$.
The bijection from $T$ to $T'$ sending $t$ to the unique $t' \in T'$
with $\val(t - t') > \delta$ is also a risometry, so by Lemma~\ref{lem:fin&iso},
it is equal to
$\beta_{q,q'}\circ\phi$. This implies that $\dir(y - \beta_{q,q'}(y)) \in V$,
first for $y \in \phi(T)$ and then also for all other $y \in S_{0,q}$.

To get that the maps $\alpha_{q,q'}\colon B_q \to B_{q'}$ move in direction $V$, it remains to check that both, the maps $\tilde{\alpha}_{\tilde{q},\tilde{q}'}$
and the maps $\tilde{\gamma}_{q,q'}$ move in direction $V$.
Let us first consider the maps $\tilde{\gamma}_{q,q'}$.
By assumption, $S_0 \ne \emptyset$, so $\tilde{\gamma}_{q,q'}(y) = y - c_q + c_{q'}$,
which moves
in direction $V$ since $\beta_{q,q'}(N_q) = N_{q'}$ and $\beta_{q,q'}$ moves in direction $V$.

To obtain that the maps $\tilde{\alpha}_{\tilde{q}, \tilde{q}'}$ move in direction $V$,
it suffices to check that we can apply (3') instead of (3) in the induction. For this, we take
$\tilde{V} := V + (\{0\}^m \times W) \subseteq k^{m+n}$. Since $B_{q} \cap S_{d,q} \ne \emptyset$,
the 0-dimensional stratum of the induction is non-empty, and it remains to check
pointwise translatability: for given
$(q,x), (q',x') \in \tilde{Q}$, $i \le n$, and $y \in \rho\1(q, x) \cap S_i$, we need to 
find an element
$y' \in \rho\1(q',x') \cap S_i$ such that
$\dir(y - y') \in \tilde{V}$.

Set $\delta := \val(q - q')$.
If $\delta \le \rado B_{q}$, then any $y' \in B_{q'}$ satisfies 
$\dir(y - y') \in V$, since $\tilde{\gamma}_{q,q'}$ moves in direction $V$
and sends $B_{q}$ to $B_{q'}$, so the risometry $B_{q} \to B_{q'}$ from
the proof of (3) yields an $y'$ with the desired properties.

If $\delta > \rado B_{q}$, then let $y'' \in \pi\1(q') \cap S_i$ be a point
obtained from the $V$-pointwise translatability
in the assumptions.
Using again that $\tilde{\gamma}_{q,q'}$ moves in direction $V$, we get
$\val(y'' - \tilde{\gamma}_{q,q'}(y)) > \delta$ and thus $y'' \in B_{q'}$.
Now we use $W$-translatability of $B_{q'}$ to move $y''$ to
the fiber $\rho\1(q',x')$.
\end{proof} 

A priori, being a t-stratification is not first order, since there might be
no bound on how complicated the straighteners in a single t-stratification are. (Recall that the straighteners are the risometries making things translation invariant; see Definition~\ref{defn:trans}.)
However, Proposition~\ref{prop:sak}~(1) says that after all, being a t-stratification is
first order; from this, we can deduce a posteriori that
all straighteners appearing in a single t-stratification can be defined uniformly.
(In fact, these uniformly defined straighteners can also directly be extracted from the proof of Proposition~\ref{prop:sak}.)

\begin{cor}\label{cor:unif-straight}
If $(S_i)_i$ is a t-stratification reflecting a definable map $\chi\colon B_0 \to \RV\eq$,
then the straighteners on all balls can be defined uniformly, i.e.,
there is a formula $\eta(x, x', y)$, where
$x, x'$ are $n$-tuples of valued field variables and $y$ is an arbitrary tuple of variables, such that
for any ball $B \subseteq S_{\ge d}$, there exists an element $b$ such that
$\eta(x, x', b)$ defines the graph of a straightener of $((S_i)_i, \chi)$ on $B$ witnessing
$d$-translatability. If $(S_i)_i$ and $\chi$ are given by formulas, then $\eta(x, x', y)$ can
be chosen to work in all $K \models \Tx$ where $(S_i)_i$ reflects $\chi$.
\end{cor}

\begin{proof}
For any formula $\eta(x, x', y)$, let $\operatorname{str}_\eta(b, \code{B})$ be a formula expressing that $\eta(x, x', b)$ defines a straightener
which witnesses $d$-translatability of $((S_i)_i, \chi)$ on the ball $B$, where $d$ is minimal with $B \cap S_d \ne \emptyset$.
Applying Proposition~\ref{prop:sak}~(1) to $((S_i)_i, \chi)$, where $Q$ is a one-point-set, yields a sentence $\psi$
such that a model $K \models \Tx$ satisfies $\psi$ if and only if
for every ball $B\subseteq K^n$, a straightener exists.
In particular, for every $K \models (\Tx,\psi)$ and every $B$, there exists an $\eta(x, x', y)$
such that $K \models \exists b\, \operatorname{str}_\eta(b, \code{B})$.
By compactness, there is a single $\eta(x, x', y)$ such that
$\Tx \cup \{\psi\}$ implies $\forall B\, \exists b\, \operatorname{str}_\eta(b, \code{B})$;
this $\eta$ uniformly defines all straighteners.
\end{proof}

\section{Proof of the existence of t-stratifications}
\label{sect:bew}

We now come to the proof of the main theorem about existence of t-stratifications:
Theorem~\ref{thm:main}.
If not specified otherwise, we will only assume Hypothesis~\ref{hyp:general}$_0$ everywhere.
In fact, the only places where we will need more than that are Proposition~\ref{prop:Stwf} and the main theorem itself.

Here is a very rough sketch of the proof
(omitting many technicalities).
Suppose we have a definable map $\chi \colon K^n \to \RV\eq$.
The overall idea is to construct the sets $S_d$ one after the other, starting with $S_n$.
Suppose that $S_n, \dots, S_{d+1}$ are already constructed and let
$X := K^n \setminus S_{\ge d+1}$ be the remainder, which we suppose to be of dimension at most $d$.
To obtain $S_d$,
we only have to find a set $X' \subseteq X$
which is at most $(d-1)$-dimensional such that on any ball not intersecting $X'$,
we have (at least) $d$-translatability; then we can set $S_d := X \setminus X'$.
However, to be able to obtain such an $X'$ in a definable way, we have to drop the
condition $X' \subseteq X$. This is not a problem; we simply shrink the sets $S_i$ we already constructed before (removing $X'$ from them).

To prove $d$-translatability on many balls $B$ (where ``many'' means: outside of a ($d-1$)-dimensional set), we roughly proceed as follows.
First note that we can always ``refine'' $\chi$, i.e., we can replace it by a map $\chi'$ such that each $\chi'$-fiber
is entirely contained in a $\chi$-fiber. We use
Proposition~\ref{prop:Stwf} (which in turn uses the Jacobian property) to refine $\chi$ in such a way that each $\chi$-fiber
$C$ separately becomes $(\dim C)$-translatable on suitable balls. The biggest difficulty then consists in
showing that on many balls $B$, these individual translatabilities fit together well enough to yield translatability of the whole map $\chi$; this is done using Lemma~\ref{lem:ZP}.
A main prerequisite for that lemma is that there exists a coordinate projection $\pi\colon B \to K^d$
such that between any two $\pi$-fibers, there exists a risometry respecting $\chi$.

One strategy to obtain the required risometries between $\pi$-fibers is as follows. Consider a fiber $F \subseteq K^n$ of a coordinate projection $\rho\colon K^n \surject K^{d-1}$.
If $d \ge 2$, then $\dim F < n$, so by induction, we may assume that there exists a t-stratification
reflecting $\chi\auf{F}$. This implies that in many cases, risometries between the $\pi$-fibers inside $F$ exist.
(``In many cases'' means: for many balls $B$ and many projections $\pi$.)
By doing this for different $\rho$, we finally obtain that in many cases, we have risometries between
any two $\pi$-fibers (Lemma~\ref{lem:IIK}). More precisely, we indeed find a 
$(d-1)$-dimensional ``bad'' set $X'$ such that for any ball $B$ not intersecting $X'$,
there exists a $\pi\colon B\to K^d$ such that the risometries exist.

In the case $d=1$, the method from the previous paragraph does not work, since we can not apply induction.
In that case, we apply Proposition~\ref{prop:sak} to the family of all $\pi$-fibers,
which implies that between many $\pi$-fibers, we find risometries respecting $\chi$.
More precisely, by doing this for all $\pi$ (and using $d = 1$), we can confine the bad set $X'$ to a 
union of finitely many ($n-1$)-dimensional hyperplanes. This argument can be repeated and each time, the dimension
of the hyperplanes drops by one. In this way, we finally obtain that $X'$ is contained in a finite set;
then we set $S_0 := X'$ and we are done.

\subsection{Sub-affine pieces}
\label{subsect:sub-aff}

\begin{defn}\label{defn:refine}
For two maps $\chi, \chi'\colon K^n \to \RV\eq$, we say that $\chi'$ \emph{refines} $\chi$ (or: $\chi'$
is a \emph{refinement} of $\chi$) if $\chi = f \circ \chi'$ for a suitable map $f$
(or, equivalently, if the partition of $K^n$ given by the fibers of $\chi'$ refines one for $\chi$).
Using Convention~\ref{conv:trans}, we also speak of maps $\chi'$ refining tuples of maps and sets.
\end{defn}

To get translatability of a definable map $\chi\colon K^n \to \RV\eq$ on certain balls, the first step is to
refine it such that each fiber $C$
is, up to a risometry, a subset of a sub-vector space of $K^n$ of the same dimension as $C$.
More precisely, we will require each $C$ to have the following property.

\begin{defn}\label{defn:sub-aff}
For any subset $C \subseteq K^n$, we define the \emph{affine direction space} of $C$
to be the sub-space $\affdir(C) \subseteq k^n$ generated by $\dir(x - x')$, where $x, x'$ run through $C$ (and $x \ne x'$). We call $C$ \emph{sub-affine} (in direction $\affdir(C)$) if for every $x \in C$, $\dim_x(C) = \dim (\affdir(C))$.
\end{defn}

We use $\dim_x(C)$ and not $\dim(C)$ in the definition to ensure
that the intersection of a sub-affine set with a ball is again sub-affine.
Note that by (\ref{it:sub-graph}) of the next lemma, $\dim (\affdir(C))$
can not be less than $\dim(C)$.

\begin{lem}\label{lem:sub-prop}
Suppose that $C \subseteq K^n$ is a definable set and that $\pi\colon K^n \surject K^{d}$ exhibits
$V = \affdir(C)$. Then we have the following.
\begin{enumerate}
\item\label{it:sub-graph} Each $\pi$-fiber contains at most one point of $C$; in particular, $C$ is the graph of a map $c\colon \pi(C) \to K^{n-d}$
and $\dim C \le d = \dim V$.
\item\label{it:sub-liso}
Suppose that $c$ is as in (\ref{it:sub-graph}) and to simplify notation, suppose that $\pi$ is the projection to the first $d$ coordinates.
Then the map
\[\phi\colon\pi(C) \times K^{n-d} \to \pi(C) \times K^{n-d},\quad
(x, y) \mapsto (x, y + c(x))\]
can be written as $\phi = \psi \circ M$, where $\psi$ is a risometry, $M \in \GL_n(\valring)$, and $\pi \circ \psi = \pi \circ M = \pi$.
\item\label{it:sub-pi}
For any coordinate projection $\rho\colon K^n \surject K^{d'}$, we have $\rhobar(\affdir(C)) \subseteq \affdir(\rho(C))$.
\end{enumerate}
\end{lem}

\begin{proof}
(1)
Two different elements $x, x'$ in the same $\pi$-fiber would have $\pibar(\dir(x - x')) = 0$,
contradicting that $\pi$ exhibits $V$.

(2)
Choose any lift $\tilde{V} \subseteq K^n$ of $V$ and define $M$ to be the linear map
sending $(x, y) \in K^d \times K^{n-d}$ to $(x, y + b(x))$, where $(x, b(x)) \in \tilde{V}$.
Since $\pi$ exhibits $V$, we have $\val(b(x)) \ge \val(x)$ and hence $M \in \GL_n(\valring)$; it remains
to check that the map $\psi = \phi \circ M\1$, which sends $(x, y)$ to $(x, y - b(x) + c(x))$, is a risometry.
For $x, x' \in \pi(C)$, Lemma~\ref{lem:vf} (\ref{it:dir-pi}) implies $\rv(x - x', c(x) - c(x')) = \rv(x - x', b(x) - b(x'))$; from this, the claim follows using Lemma~\ref{lem:vf} (\ref{it:rv-sum}).

(3)
For any pair of points $x, x' \in C$ we have to check that $\rhobar(\dir(x - x')) \in \affdir(\rho(C))$.
If $\rhobar(\dir(x - x')) = 0$, there is nothing to prove; otherwise Lemma~\ref{lem:vf} (\ref{it:dir-pi})
implies $\rhobar(\dir(x - x')) = \dir(\rhobar(x) - \rhobar(x'))$.
\end{proof}

Being sub-affine is closely related to translatability.

\begin{lem}\label{lem:sub-aff}
Let $B \subseteq K^n$ be a ball and $C \subseteq B$ a definable subset.
\begin{enumerate}
\item\label{it:trans-in-sub}
If $C$ is $V$-translatable on $B$ for some $V\subseteq k^n$, then $V \subseteq \affdir(C)$.
\item\label{it:sub=>trans}
If there is an exhibition $\pi\colon B \to K^d$ of $V := \affdir(C)$
with $\pi(C) = \pi(B)$, then $C$ is $V$-translatable on $B$.
\end{enumerate}
\end{lem}
\begin{proof}
(1) Clear.

(2) Assume without loss $0 \in B$. Then the risometry $\psi$ obtained from Lemma~\ref{lem:sub-prop} (\ref{it:sub-liso})
sends $B$ to itself and it is a straightener.
\end{proof}

Now we can formulate the main result of this subsection.
Its proof is the only place in the proof of Theorem~\ref{thm:main} where the Jacobian property is needed.

\begin{prop}\label{prop:Stwf}
Assume Hypothesis~\ref{hyp:general}$_{n-1}$ and
let $\chi\colon K^n \to \RV\eq$ be a $\emptyset$-definable map.
Then there exists a $\emptyset$-definable refinement $\chi'$ of $\chi$ such that
each $\chi'$-fiber $C' \subseteq K^n$ is sub-affine.
\end{prop}

In the proof, we will need the following lemma.

\begin{lem}\label{lem:jac2subaff}
Suppose that $X \subseteq K^d$ is a definable set of dimension $d$ and that
$Z \subseteq K^{d+1}$ is the graph of a definable function $f\colon X \to K$
that has the Jacobian property (Definition~\ref{defn:jac}). Then $Z$ is sub-affine.
\end{lem}
\begin{proof}[Proof of Lemma~\ref{lem:jac2subaff}]
We may assume that $f$ is not constant; let $z \in K^n \setminus \{0\}$ be as in Definition~\ref{defn:jac}
and write $\pi\colon K^{d+1} \surject K^d$ for the projection to the first $d$ coordinates.

Suppose that $Z$ is not sub-affine. Since $\dim Z = d$, this implies $\affdir(Z) = k^{d+1}$,
so by definition, there exist
$d+1$ pairs of points $x'_i, x''_i \in X$ with $x'_i \ne x''_i$ such that
$\big(\dir(x'_i - x''_i, f(x'_i) - f(x''_i))\big)_i$ is a basis of $k^{d+1}$.
Set $x_i := x'_i - x''_i$ and $y_i := f(x'_i) - f(x''_i)$.
Using these notations, the inequality in Definition~\ref{defn:jac} becomes
\[
\tag{$+$}
\vali(y_i - \langle z, x_i\rangle) > \val(z) + \val(x_i)
.
\]

Suppose first that $\val(z) < 0$. Choose $i$ with
$\langle \dir(z), \pibar(\dir(x_i, y_i)) \rangle \ne 0$. Then
by Lemma~\ref{lem:vf} (\ref{it:dir-pi}) we have $\pibar(\dir(x_i, y_i)) = \dir(x_i)$ and
$\vali(y_i) \ge \val(x_i) > \val(z) + \val(x_i)$; moreover,
Lemma~\ref{lem:vf} (\ref{it:dir-scal}) implies $\vali(\langle z, x_i\rangle) = \val(z) + \val(x_i)$.
Thus we have $\vali(y_i - \langle z, x_i\rangle) = \val(z) + \val(x_i)$,
contradicting ($+$).

Now suppose $\val(z) \ge 0$. Then ($+$) implies
$\val((x_i, y_i) - (x_i, \langle z, x_i\rangle)) > \val(x_i, y_i)$ and hence
$\dir (x_i, y_i) = \dir(x_i, \langle z, x_i\rangle)$.
Now $\tilde{V} := \{(x, \langle z, x\rangle) \mid x \in K^d\}$ is a $d$-dimensional subspace of $K^{d+1}$,
so each $\dir(x_i, \langle z, x_i\rangle)$ lies in the $d$-dimensional subspace $\res(\tilde V) \subseteq k^{d+1}$
(by Lemma~\ref{lem:vf} (\ref{it:dir-res})), contradicting that the $\dir(x_i, y_i)$ form a basis of $k^{d+1}$.
\end{proof}

\begin{proof}[Proof of Proposition~\ref{prop:Stwf}]
We will prove the following claim: For any $\emptyset$-definable set $C \subseteq K^n$ of dimension $d$, there is a
$\emptyset$-definable map $\tilde{\chi}\colon C \to \RV\eq$ such that each
$\tilde\chi$-fiber $\tilde{C} \subseteq C$ of dimension $d$ satisfies $\dim(\affdir(\tilde{C})) = d$.

Once we have this, we can finish the proof of the proposition as follows. We do an induction
over the maximum of the dimensions of $\chi$-fibers that are not
sub-affine. Denote this maximum by $d$. On each $\chi$-fiber $C$
of dimension $d$ which is not sub-affine, we refine $\chi$ as follows.

First, we apply the claim to $C$ (with $\code{C}$ added to the language), which yields a $\code{C}$-definable map $\tilde{\chi}\colon C \to \RV\eq$.
Now consider a $\tilde{\chi}$-fiber $\tilde{C}$ of dimension $d$.
By Lemma~\ref{lem:loc-dim}, the set
$\tilde{D} := \{x \in \tilde{C} \mid \dim_x \tilde{C} < d\}$ has dimension less than $d$,
so in particular $\tilde{C} \setminus \tilde{D}$ is sub-affine. Refine $\tilde{\chi}$ such that each
$\tilde{C} \setminus \tilde{D}$ becomes a separate fiber. The result is a
refinement of $\chi$ whose fibers of dimension $d$ are sub-affine,
so then we are done by the induction hypothesis.

It remains to prove the claim. For $d = n$, there is nothing to show, so we assume $d < n$. We define $\tilde{\chi}(x)$ (for $x \in C$)
to be the tuple $(\chi_\pi(x))_\pi$, where $\pi$ runs over all coordinate projections $\pi\colon K^n \surject K^d$ and where
$\chi_\pi\colon C \to \RV\eq$ is defined as follows:
\begin{enumerate}
\item
Let $C_0$ be the union of all $\pi$-fibers $F = \pi\1(y) \cap C$ (for $y \in \pi(C)$) satisfying $\dim F = 0$
and let $C_1 := C \setminus C_0$ be the remainder. We set $\chi_\pi(x) = 1_\pi$ for all $x \in C_1$, where
$1_\pi \in \RV\eq$ is any element that is not used again (i.e., $C_1$ is one fiber of $\chi_\pi$).
\item
By Lemma~\ref{lem:fin&iso}~(2) and using stable embeddedness of $\RV$,
we get a $\emptyset$-definable map $\chi_0\colon C_0 \to \RV\eq$
that is injective on each $\pi$-fiber (namely $\chi_0(x) = \code{\rv\!\big(x - \big(\pi\1(\pi(x)) \cap C_0\big)\big)}$).
\item
A $\chi_0$-fiber $C'=\chi_0\1(\sigma) \subseteq C_0$ can be seen
as the graph of a $\sigma$-definable function
$f\colon \pi(C') \to K^{n-d}$. Since the theory has the Jacobian property up to dimension $n - 1 \ge d$, we find a $\sigma$-definable map $\chi'_\sigma \colon \pi(C') \to \RV\eq$ such that on each $\chi'_\sigma$-fiber of dimension $d$, each coordinate of $f$ has the Jacobian property
(Definition~\ref{defn:jac}).
For $x \in C_0$, set
$\chi_\pi(x) := (\chi_0(x), \chi'_{\chi_0(x)}(\pi(x)))$.
\end{enumerate}

We now have to check that if $\tilde{C}$ is a $d$-dimensional $\tilde{\chi}$-fiber,
then $\dim(\affdir(\tilde{C})) = d$, so assume for contradiction that 
$d' := \dim(\affdir(\tilde{C})) > d$. Choose an exhibition $\rho\colon K^n \surject K^{d'}$ of $\affdir(\tilde{C})$. Then $\dim \rho(\tilde{C}) = d$ since otherwise,
there would be $\rho$-fibers containing several points of $\tilde{C}$, contradicting
that $\rho$ exhibits $\affdir(\tilde{C})$. Next choose a coordinate projection
$\pi'\colon K^{d'} \surject K^d$ such that for $\pi := \pi' \circ \rho$, we still
have $\dim \pi(\tilde{C}) = d$, and choose an arbitrary decomposition of $\pi'$ into two coordinate projections $\rho'\colon K^{d'} \surject K^{d+1}$ and $\pi''\colon K^{d+1} \to K^d$.
Note that by the choice of $\rho$, we have $\rhobar'(\rhobar(\affdir(\tilde{C}))) = k^{d+1}$.

Since $\dim(\pi(\tilde{C})) = \dim(\tilde{C})$, we have $\chi_\pi(\tilde{C}) \ne 1_\pi$, so by definition of $\chi_\pi$,
$Z := \rho'(\rho(\tilde{C}))$ is the graph of a function $f\colon \pi(\tilde{C}) \to K$ that
has the Jacobian property; hence $Z$ is sub-affine by Lemma~\ref{lem:jac2subaff}.
However, this contradicts that by Lemma~\ref{lem:sub-prop} (\ref{it:sub-pi}) we have $\affdir(Z) \supseteq \rhobar'(\rhobar(\affdir(\tilde{C}))) = k^{d+1}$.
\end{proof}

\subsection{Merging translatability}

In the previous subsection, we obtained some first translatability separately for
each fiber of a definable map $\chi\colon K^n\to\RV\eq$. Now we will show how this can be merged
to translatability of the whole map $\chi$ (under a lot of technical assumptions).
We start with a lemma that allows us to relate affine direction spaces of different
fibers.

\begin{lem}\label{lem:sub-GlI}
Let $B \subseteq K^n$ be a ball and let $C, C' \subseteq B$ be non-empty definable subsets that are
sub-affine in directions $V$ and $V'$, respectively. Suppose that
$\pi\colon B \to K^d$ exhibits $V$ and that for any two elements
$y_1, y_2 \in \pi(B)$, there exists a risometry
$\pi\1(y_1) \to \pi\1(y_2)$ respecting $(C, C')$.
Then $V \subseteq V'$.
\end{lem}

\begin{proof}
It suffices to find $x'_1, x'_2 \in C'$ with $\dir(x'_1 - x'_2) = v$ for any
given $v \in V \setminus \{0\}$, so let such a $v$ be given.

Choose any $x'_1 \in C'$ and set $y_1 := \pi(x'_1)$.
Any fiber of $\pi$ contains exactly one element of $C$;
let $x_1$ be this unique element of $C \cap \pi\1(y_1)$. Choose $y_2 \in \pi(B)$ such that
$\dir(y_1 - y_2) = \pibar(v)$ and $\val(y_1 - y_2) = \val(x'_1 - x_1)$.
Now let $x_2$ and $x_2'$ be the images of $x_1$ and $x_1'$ under
a risometry $\phi\colon \pi\1(y_1) \iso \pi\1(y_2)$.

We have $x_2 \in C$, so $\dir(x_1 - x_2) = v$. Since $\phi$ is a risometry,
we have $\rv(x'_1 - x_1) = \rv(x'_2 - x_2)$, so
$\val((x'_1 - x'_2) - (x_1 - x_2)) > \val(x'_1 - x_1) = \val(x_1 - x_2)$,
which implies $\rv(x'_1 - x'_2) = \rv(x_1 - x_2)$ and thus $\dir(x'_1 - x'_2) = v$.
\end{proof}

The following lemma is the main tool to prove $V$-translatability of a map
$\chi\colon B \to \RV\eq$, where $B \subseteq K^n$ is a ball and $V \subseteq k^n$ is a $d$-dimensional vector space.
Let $\pi\colon B \surject K^d$ exhibit $V$.
The prerequisites are (i) that there exist risometries between the $\pi$-fibers respecting $\chi$,
(ii) that the fibers of $\chi$ are sub-affine,
(iii) that there exists a $\chi$-fiber $C$ with $\affdir C = V$,
and (iv) that we already have a partial t-stratification which works outside 
of a $d$-dimensional set.
However, in applications of the lemma, we will not be able to ensure (iv) simultaneously with
(i) -- (iii); therefore, we allow (i) -- (iii) to apply to a refinement $\chi'$
of $\chi$, which is enough to get the result.

\begin{lem}\label{lem:ZP}
Suppose that we have definable
sets $C \subseteq B \subseteq K^n$, definable maps $\chi, \chi' \colon B \to \RV\eq$,
an integer $d \in \{1,\dots, n\}$, a definable partition
$(S_i)_{d \le i \le n}$ of $B$ and a coordinate projection $\pi\colon B \to K^d$,
with the following properties:
\begin{itemize}
\item
$B$ is a ball;
\item
$\dim S_i \le i$;
\item
for any ball $B'\subseteq B \setminus S_d$,
$((S_i)_i, \chi)$ is $j$-translatable on $B'$, where $j$ is minimal with $B' \cap S_{j} \ne \emptyset$;
\item
$\chi'$ is a refinement of $((S_i)_i, \chi)$ (in the sense of Definition~\ref{defn:refine} and Convention~\ref{conv:trans});
\item
for each pair of points $x, x' \in \pi(B)$, there exists a
definable risometry $\pi\1(x) \to \pi\1(x')$ respecting $\chi'$;
\item
$C$ is a $\chi'$-fiber whose
affine direction space $V := \affdir(C)$ is exhibited by $\pi$
(in particular, $\dim C = \dim V = d$).
\end{itemize}
Then $((S_i)_i, \chi)$ is $V$-translatable.
\end{lem}

\begin{proof}
By Lemma~\ref{lem:sub-GlI}, for any $\chi'$-fiber $C'$ we have $V \subseteq \affdir(C')$.
In particular,
if $C' \subseteq S_d$, then $\affdir(C') = V$ and
$C'$ is $V$-translatable on $B$ by Lemma~\ref{lem:sub-aff} (\ref{it:sub=>trans}).

\smallskip

\textbf{Claim 1.} If $B' \subseteq B$ is a ball with $B' \cap S_d = \emptyset$,
then $V \subseteq W := \tsp_{B'}((S_i)_i)$.

\emph{Proof of Claim 1.} Set $d' := \dim W$ and
let $\pi'\colon B' \to K^{d'}$ be an exhibition of $W$.
The $\pi'$-fibers of $S_{d'}$ are finite but non-empty. Choose a subball $B'' \subseteq B'$
such that $B'' \cap (\pi')\1(x) \cap S_{d'}$ is a singleton for each $x \in \pi'(B'')$.
Then $\affdir(S_{d'} \cap B'') = \tsp_{B''}((S_i)_i) = W$.
Now choose any $\chi'$-fiber $C' \subseteq S_{d'}$ with $\dim (C' \cap B'') = d'$.
Then $W \subseteq \affdir C'$ and $\dim (\affdir C') = \dim C' = d'$ together imply
$W = \affdir C' \supseteq V$.

\smallskip

If $S_d = \emptyset$, then we are done using $B' = B$, so from now on suppose $S_d \ne \emptyset$.

\smallskip

\textbf{Claim 2.} Fix $x \in \pi(B)$, let 
$F = \pi\1(x)$ be the fiber over $x$, and set $T_i := S_{i + d} \cap F$
for $i \le n - d$. Then $(T_i)_{i \le n - d}$ is a t-stratification of $F$
reflecting $\chi\auf{F}$.

\emph{Proof of Claim 2.}
Using that the risometries between the $\pi$-fibers respecting $\chi'$ in particular
respect $S_{j+d}$, we obtain $\dim T_j \le j$. Consider a ball $B' \subseteq B$ intersecting $F$ non-trivially and set
$B'_x := B' \cap F$. We have to show that if $B'_x \subseteq T_{\ge j}$, then
$((T_i)_i, \chi)$ is $j$-translatable on $B_x'$.
For $j = 0$, there is nothing to do, so suppose $j \ge 1$.
Then $B'_x \cap S_d = \emptyset$, and using that any
$\chi'$-fiber $C'$ contained in $S_d$ is
$V$-translatable, we get $B' \cap S_d = \emptyset$.
Now Claim~1 together with Lemma~\ref{lem:trans&fiber} implies
$j$-translatability on $B'_x$. 

\smallskip

\textbf{Claim 3.} $(S_i)_i$ is $V$-$\pi$-pointwise translatable (see Definition~\ref{defn:pt-trans}).

\emph{Proof of Claim 3.}
Let $x, x' \in \pi(B)$ and $y \in \pi\1(x)$ be given; we need
to find $y' \in \pi\1(x')$ with $\dir (y - y') \in V$ such that
$y$ and $y'$ are elements of the same set $S_i$.
Set $\delta := \val(x - x')$.
If $\ball{y, \ge \delta} \cap S_d = \emptyset$, then
$y'$ is obtained using Claim~1.
Otherwise, let $C' \subseteq S_d$ be a $\chi'$-fiber
intersecting $\ball{y, \ge \delta}$ non-trivially,
let $z, z'$ be the unique elements of $C' \cap \pi\1(x)$
and $C' \cap \pi\1(x')$, respectively, and
let $y'$ be the image of $y$ under a risometry
$\pi\1(x) \to \pi\1(x')$ respecting $\chi'$.
Now $\val(y - z) \ge \delta$ and $\rv(y - z) = \rv(y' - z')$ together imply
$\val((y - y') - (z - z')) > \delta$, and thus
$\dir(y-y') = \dir(z-z') \in V$.

\smallskip

Now, Claims~2 and 3 (together with $S_d \ne \emptyset$) are all we need to
apply Corollary~\ref{cor:sak}, which yields the desired
$V$-translatability.
\end{proof} 

The next lemma will be useful to prove the prerequisites of the previous lemma.
More precisely, given $\pi$ and a $\chi'$ as above, it yields a way to prove
that there exist risometries respecting $\chi'$ between any to $\pi$-fibers.


\begin{lem}\label{lem:IIK}
Let the following be given:
\begin{itemize}
\item a definable map $\chi\colon B \to \RV\eq$, where $B \subseteq K^n$ is a ball;
\item a vector space $V \subseteq k^n$ of dimension $d \ge 1$ exhibited by $\pi\colon B \to K^d$;
\item a $\chi$-fiber $C \subseteq B$ with $\affdir(C) \subseteq V$.
\end{itemize}
Suppose that the for each coordinate projection $\rho\colon K^d \surject K^{d-1}$ and each
$y \in \rho(\pi(B))$, $\chi$ is $1$-translatable
on the fiber $\pi\1(\rho\1(y))$.

Then for any $x_1, x_2 \in \pi(B)$, there exists a definable
risometry $\pi\1(x_1) \to \pi\1(x_2)$ respecting $\chi$.
\end{lem}

\begin{rem}
A posteriori, this implies $\dim C = d$, so $\affdir(C) = V$ and $C$ is sub-affine.
\end{rem}

\begin{proof}[Proof of Lemma~\ref{lem:IIK}]
It is enough to find such risometries $\pi\1(x_1) \to \pi\1(x_2)$ under the assumptions that $x_1$ and $x_2$ differ
in only one coordinate and that $\pi\1(x_1) \cap C \ne \emptyset$. Indeed,
the existence of such a risometry implies that we also have $\pi\1(x_2) \cap C \ne \emptyset$,
so by repeatedly applying this (starting with a $\pi$-fiber intersecting $C$ non-trivially and
modifying coordinates one by one), we first get that every $\pi$-fiber intersects $C$ non-trivially, and
then we obtain risometries between any two fibers by composition.

So suppose now that $x_1$ and $x_2$ differ only in one coordinate and
let $\rho \colon K^d \surject K^{d-1}$ be the coordinate projection
satisfying $\rho(x_1) = \rho(x_2) =: y$; let $F := \pi\1(\rho\1(y)) \subseteq B$ be the corresponding fiber.
By assumption, there exists a one-dimensional
$W \subseteq \ker (\rhobar \circ \pibar) \subseteq k^n$ such that
$\chi$ is $W$-translatable on $F$.
In particular, the non-empty set $C \cap F$ is $W$-translatable,
so $W \subseteq V$ by Lemma~\ref{lem:sub-aff}~(\ref{it:trans-in-sub}).
Since $\dim (\ker (\rhobar \circ \pibar) \cap V )= 1$,
$W$ is equal to this intersection, so $\pi\auf{F}$ exhibits $W$.
From a translater $(\alpha_x)_{x \in \pi(F - F)}$
of $\chi\auf F$ with respect to $\pi$ (see Definition~\ref{defn:translater}), we obtain a risometry
$\phi\colon \pi\1(x_1) \to \pi\1(x_2)$ by restricting $\alpha_{x_2 - x_1}$ to $\pi\1(x_1)$.
\end{proof}

\subsection{The big induction}

This subsection contains the actual proof of the main theorem.
We first prove the case $n = 1$ separately.

\begin{lem}\label{lem:n=1}
For every $\emptyset$-definable map $\chi\colon K \surject Q \subseteq \RV\eq$, there exists a finite
$\emptyset$-definable set $T_0 \subseteq K$ such that $\chi$ is constant on each ball $B \subseteq K \setminus T_0$.
\end{lem}

\begin{proof}
By Hypothesis~\ref{hyp:general} (\ref{it:strat1}), for each $q \in Q$,
there exists a finite set $S_q$ such that each ball $B \subseteq K \setminus S_q$
is either disjoint from $\chi\1(q)$ or contained in $\chi\1(q)$.
By Lemma~\ref{lem:dim} (\ref{it:dim-u}), the union $T_0 := \bigcup_{q} S_q$
is $0$-dimensional, so by Remark~\ref{rem:0fin}, it is finite.
The construction of $T_0$ ensures that $\chi$ is constant on each ball $B \subseteq K \setminus T_0$.
\end{proof}

Now we are ready for the main theorem. For readers who jumped directly to this point,
let us recall: the basic notation is fixed in Subsection~\ref{subsect:notn},
the language $\LHen$ is introduced in Definition~\ref{defn:LHen},
and t-stratifications are introduced in Definition~\ref{defn:t-strat}.
Recall also Remark~\ref{rem:param}, by which we can replace each occurrence of ``$\emptyset$-definable'' in the theorem by ``$A$-definable'' for some fixed parameter set $A$ (which, using compactness, yields that the theorem works uniformly in families).

\begin{thm}\label{thm:main}
Fix $n \in \bbN$.
Let $\Lx$ be an expansion of the valued field language $\LHen$ and let $K$ be
an $\Lx$-structure whose theory satisfies Hypothesis~\ref{hyp:general}$_{n-1}$.
(In particular, $K$ is a Henselian valued field of equi-characteristic $0$.)
Then, for every $\emptyset$-definable ball $B_0 \subseteq K^n$ and every $\emptyset$-definable map
$\chi\colon B_0 \to \RV\eq$, there exists
a $\emptyset$-definable t-stratification $(S_i)_{i \le n}$ of $B_0$ reflecting $\chi$.
\end{thm}

\begin{proof}
The case $n = 1$ is exactly Lemma~\ref{lem:n=1}. For $n \ge 2$,
we do a big induction on $n$, i.e., we assume that the theorem holds for all smaller $n$.

By extending $\chi$ trivially outside of $B_0$, we may suppose $B_0 = K^n$.

By decreasing induction on $d$, we prove the following.
\labelledclaim{$\star_d$}{%
There exists a $\emptyset$-definable partition $(S_i)_{d \le i \le n}$ of $K^n$
with $\dim S_i \le i$ such that
for any ball $B \subseteq S_{> d}$,
$((S_i)_{i}, \chi)$ is $j$-translatable on $B$, where $j$ is minimal with $B \cap S_{j} \ne \emptyset$.
}
Note that $(\star_0)$ implies the theorem. (For balls intersecting $S_0$, there is nothing to prove.)

The start of induction $(\star_n)$ is trivial (set $S_n = K^n$).
Now suppose that $(S_i)_{d \le i \le n}$ is given such that $(\star_d)$ holds (for some $d \ge 1$). It suffices to find a set $S_{d - 1}$
of dimension at most $d - 1$ such that on any ball $B \subseteq K^n \setminus S_{d-1}$,
$((S_i)_i, \chi)$ is $d$-translatable; after that, we obtain $(\star_{d-1})$ using the partition
$S_{d-1}, (S_i \setminus S_{d-1})_{i \ge d}$. Moreover, it is enough to check $d$-translatability on balls $B$
with $B \cap S_{d} \ne \emptyset$.

We have to do the case $d = 1$ separately.

\smallskip

\textbf{The case $d \ge 2$:}

First, we choose a refinement $\chi'$ of
$((S_i)_i, \chi)$ whose fibers
are sub-affine (using Proposition~\ref{prop:Stwf}).
Now consider a coordinate projection $\pi \colon K^n \surject K^{d-1}$.
By induction on $n$ (and using $d \ge 2$),
we can find t-stratifications of the fibers of $\pi$
reflecting $\chi'$ on the fibers.
Taking the union of corresponding strata of different fibers yields a 
$\emptyset$-definable partition
$(T_{i})_{i \le n - d + 1}$ of $K^n$ with $\dim T_i \le i + d - 1$.
Define $S_{d-1}$ to be the union of the sets $T_0$ for all coordinate projections
$\pi \colon K^n \surject K^{d-1}$.

Now let a ball $B \subseteq K^n \setminus S_{d-1}$ with $B \cap S_d \ne \emptyset$
be given;
we have to prove that $((S_i)_i, \chi)$ is $d$-translatable on $B$.
We will do this by applying Lemma~\ref{lem:ZP} to $B$, $\chi$, $\chi'$,
and $(S_i)_i$; so let us produce the remaining ingredients.

Let $C \subseteq S_d$ be any $\chi'$-fiber intersecting $B$ non-trivially,
let $V \subseteq k^n$ be $d$-dimensional such that $\affdir(C) \subseteq V$
(which exists since $\dim C \le d$), and choose
an exhibition $\pi\colon B \to K^d$ of $V$. The only missing ingredient for
Lemma~\ref{lem:ZP} is now that between any two $\pi$-fibers, there exists a definable risometry respecting $\chi'$;
this then also implies $\dim C = d$ and hence $\affdir(C) = V$.

To get the risometries between the fibers, we apply Lemma~\ref{lem:IIK}
to $\chi'$, $V$, and $C$. Suppose that
$\rho\colon K^d \surject K^{d-1}$ is a coordinate projection and $F$ is a fiber
of $\pi' := \rho \circ \pi$.
Consider the partition $(T_i)_{i \le n-d+1}$ of $K^n$
obtained from t-stratifications of the fibers of $\pi'$
in the above definition of $S_{d-1}$. Since $T_0 \subseteq S_{d-1}$
we have $B \cap T_0 = \emptyset$, so in particular,
$\chi'\auf{F}$ is 1-translatable
on $B \cap F$, which is what we need for Lemma~\ref{lem:IIK}.

\smallskip

\textbf{The case $d = 1$:}

Recall that we do already have a partition $(S_i)_{i \ge 1}$ which is good
outside of $S_1$. We will now carry out an additional induction, during which
the ``bad set'' will become ``more and more 0-dimensional''. More precisely,
consider the following statement for $e \in \{0, \dots, n\}$.

\labelledclaim{$\star\star_e$}{%
There exists a family of definable sets $X_\rho$ parametrized by the
coordinate projections $\rho\colon K^n \surject K^e$, such that
$\rho(X_\rho)$ is finite, $\dim X_\rho \le 1$, and
for any ball $B \subseteq K^n \setminus \bigcup_\rho X_\rho$,
$((S_i)_{i}, \chi)$ is $1$-translatable on $B$.
}

Write $X := \bigcup_\rho X_\rho$ for the union.
The statement $(\star\star_0)$ follows from $(\star_{1})$, since we can
take $X = X_\rho = S_1$ (where $\rho\colon K^n \surject K^0$).
The statement $(\star\star_n)$ is what we want to prove; in that case,
$\rho = \id_{K^n}$ implies that $X$ itself is finite, so we can set $S_0 = X$
(and replace $S_i$ by $S_i \setminus S_0$ for $i \ge 1$);
then $(\star\star_n)$ implies $1$-translatability on balls $B \subseteq S_{\ge 1}$
and $(\star_1)$ implies $d$-translatability on balls $B \subseteq S_{\ge d}$ for $d \ge 2$.

Thus it remains to prove ``$(\star\star_{e}) \Rightarrow (\star\star_{e+1})$'' for $0 \le e < n$.
Let $X = \bigcup_\rho X_\rho$ be given for $e$, and let us construct a set
$X'$ for $e + 1$.
We start by choosing a refinement $\chi'$ of $((S_i)_i, \chi, (X_\rho)_\rho)$ whose
fibers are sub-affine.

Let $\rho\colon K^n \surject K^e$ and $\pi\colon K^n \surject K$ be coordinate projections
``projecting to different coordinates'', i.e., such that $(\rho, \pi)\colon K^n \to K^{e} \times K$ is surjective.

By the main induction on $n$, we can find t-stratifications of the fibers of $\pi$
reflecting $\chi'$ on the fibers.
By Proposition~\ref{prop:sak} (2), there exists a definable map
$\chi_0\colon K \to \RV\eq$ such that for any $\chi_0$-fiber
$C_0 \subseteq K$ and any $x, x' \in C_0$, we have a definable risometry $\pi\1(x) \to \pi\1(x')$ respecting $\chi'$.
Lemma~\ref{lem:n=1} yields a finite subset
$T_0 \subseteq K$ such that $\chi_0$ is constant on each ball $B' \subseteq K \setminus T_0$.
Recall that $\picomp\colon K^n \surject K^{n-1}$ denotes the ``complement'' of $\pi$ and define the set $X_{\rho,\pi}$ as follows:
\[
X_{\rho,\pi} := \{x \in \pi\1(T_0) \mid \picomp(x) \in \picomp(X_\rho)\}
.
\]
We define $X'$ to be the union of all such $X_{\rho,\pi}$
(for all $\rho,\pi$ as above).

Since $T_0$ is finite, $\dim X_{\rho,\pi} \le \dim X_\rho \le 1$ and $(\rho,\pi)(X_{\rho,\pi})$ is finite,
so it remains to check that on a ball $B \subseteq K^n \setminus X'$,
$((S_i)_{i}, \chi)$ is $1$-translatable. If $B \cap X = \emptyset$, then we know this by induction on $e$, so suppose that
$B \cap X_\rho \ne \emptyset$ for some $\rho\colon K^n \surject K^e$.

Let $C \subseteq X_\rho$ be a $\chi'$-fiber with
$C \cap B \ne \emptyset$; note that $\dim C \le \dim X_\rho \le 1$. If $\dim (C \cap B) = 0$, then let $\pi\colon K^n \surject K$ be any coordinate projection projecting to a different coordinate than $\rho$.
Otherwise, set $V := \affdir(C)$ and let $\pi\colon K^n \surject K$ be an exhibition of $V$.
Since $\rho(C)$ is finite, we have $V \subseteq \ker \rhobar$, so in this case too,
$\rho$ and $\pi$ project to different coordinates.

Let $\chi_0$, $T_0$ be as in the construction of $X_{\rho,\pi}$.
Then $\pi(B) \cap T_0 = \emptyset$, since otherwise, for $x \in \pi(B) \cap T_0$
and $y \in B \cap X_\rho$, the point $y' \in K^n$
with $\pi(y') = x$ and $\picomp(y') = \picomp(y)$ lies both in $B$
and in $X_{\rho,\pi}$,
contradicting $B \cap X_{\rho,\pi} = \emptyset$.
By our choice of $T_0$, this implies that $\chi_0$ is constant on $\pi(B)$ and
thus there are risometries respecting $\chi'$ between any two fibers $\pi\1(x)$ (for $x \in \pi(B)$).
In particular, $C$ intersects every fiber non-trivially and thus $\dim (C \cap B) = 1$.

Now we can apply Lemma~\ref{lem:ZP} to
$C \cap B$, $B$,
$\chi$, $\chi'$, $\pi$, and the partition $(S_1 \cup X, (S_i \setminus X)_{i \ge 2})$ (restricted to $B$); this yields
that $((S_i)_i, \chi)$ is $V$-translatable on $B$, which is what
we had to show.
\end{proof}

\subsection{Corollaries}
\label{subsect:cor}

Using compactness, we can deduce a version of the main theorem which
works uniformly for all models of
a theory $\Tx$ satisfying Hypothesis~\ref{hyp:general}, and also for all models of a finite subset of $\Tx$,
provided that the notion of t-stratification makes sense.
In particular, we get t-stratifications
in all Henselian valued fields of sufficiently big residue characteristic
(both, in the equi-characteristic and the mixed characteristic case).
Note that in equi-characteristic, there is no good notion of dimension of
a definable set; there, ``$\dim S_i = i$'' means that we stupidly
apply Definition~\ref{defn:dim}. However, in the case of the pure
valued field language, this problem will be solved by
Corollary~\ref{cor:alg}, which says that we can choose the t-stratification such that each set $S_{\le i}$
is Zariski closed and has dimension $i$ in the algebraic sense.

\begin{cor}\label{cor:pos-char}
Suppose $\Tx$ is an $\Lx$-theory satisfying Hypothesis~\ref{hyp:general}.
Let $\chi$ be an $\Lx$-formula defining a map $\chi_K \colon K^n \to \RV\eq$ (for any model $K \models \Tx$).
Then there exist $\Lx$-formulas $\psi_0,\dots,\psi_n$
and a finite subset $\Tx_0 \subseteq \Tx$ such that for each model $K$ of
$\Tx_0$, $(\psi_i(K))_i$ is a t-stratification of $K^n$ reflecting $\chi_K$.
(For this to make sense, we assume that $\Tx_0$ in particular says that $K$ is
a valued field.) 
\end{cor}
\begin{proof}
By Theorem~\ref{thm:main}, we find formulas $(\psi_i)_i$ defining a t-stratification
for any fixed model $K \models \Tx$.
Moreover, by Corollary~\ref{cor:unif-straight}, we also find a formula $\eta$ (depending on
$(\psi_i)_i$) defining the corresponding straighteners on all balls $B \subseteq K^n$ (using parameters).
This allows us to formulate a first order sentence which holds in an $\Lx$-structure $K'$
iff $(\psi_i(K'))_i$ is a t-stratification reflecting $\chi_{K'}$, namely:
\labelledclaim{$\bigtriangleup$}{%
For each $i$, $\psi_i(K')$ is either empty or has dimension $i$ in the sense of
Definition~\ref{defn:dim}, and\\
for each ball $B \subseteq (K')^n$, there exists a parameter $b$
such that $\eta(K', b)$ defines a straightener on $B$ which witnesses that
$((\psi_i(K'))_i, \chi_{K'})$ is $j$-translatable on $B$, where $j$ is minimal with $B \cap S_{j} \ne \emptyset$.
}
By compactness, $\psi_i$ and $\eta$ can be chosen such that ($\bigtriangleup$)
holds in all models of $\Tx$. Moreover,
($\bigtriangleup$) then follows already from a finite subset of $\Tx$.
\end{proof}

The next corollary says that in some sense, the risometry type of a definable
subset of $K^n$, or, more generally, of a definable map $K^n \to \RV\eq$,
can be encoded using only $\RV$-data.

\begin{cor}\label{cor:isotyp}
We assume Hypothesis~\ref{hyp:general}.
Let $\chi_q \colon K^n \to \RV\eq$ be a $\emptyset$-definable family of maps,
parametrized by $q \in Q$ (for some definable set $Q$ in any sort).
Then there exists a $\emptyset$-definable map $\chi'\colon Q \to \RV\eq$ such that
$\chi'(q_1) = \chi'(q_2)$ implies that there exists a
$(q_1,q_2)$-definable risometry $\phi\colon K^n \to K^n$ with
$\chi_{q_1} \circ \phi = \chi_{q_2}$.
This also works uniformly for all models $K$ of a finite subset of $\Tx$.
\end{cor}
\begin{proof}
If we add a constant symbol for $q$ to the language, then
Corollary~\ref{cor:pos-char} yields uniformly defined
t-stratifications reflecting $\chi_q$ in each model of a finite subset of $\Tx$
and for each $q \in Q$.
Now $\chi'$ is obtained from Proposition~\ref{prop:sak}~(2).
\end{proof}

In Section~\ref{sect:alg} (where we prove an algebraic version of the main result),
we will give an algebraic version of this corollary (Corollary~\ref{cor:alg-isotyp}).
(That version follows directly from Corollary~\ref{cor:isotyp}, but thematically,
it fits better into Section~\ref{sect:alg}.)

\subsection{Characterizations of reflection}
\label{subsect:reflect}

For applications of the main theorem, it will be useful to understand more precisely what it means
that a t-stratification $(S_i)_i$ reflects a map. Proposition~\ref{prop:reflect} gives
different equivalent conditions for this; in
particular, there exists a finest map reflected by $(S_i)_i$---the ``rainbow'' of $(S_i)_i$.
(If one thinks of the fibers of this map as having different colors, then it indeed looks a bit like a rainbow,
in particular near $S_1$.)

A simple consequence of Proposition~\ref{prop:reflect} (which was not clear from the definition of reflection) is that if $(S_i)_i$ reflects both $\chi$ and $\chi'$, then it also reflects $(\chi, \chi')$ (Remark~\ref{rem:reflect}). Other consequences are Corollary~\ref{cor:refine} and Lemma~\ref{lem:small-changes}, which together will allow us to ``enhance'' t-stratifications in the following sense.
Given $(S_i)_i$, we will find $(S'_i)_i$ which reflects at least as much as $(S_i)_i$ and which has additional good properties.

The last results of this subsection give some more information about the fibers $C$ of a rainbow. These will be needed in Lemma~\ref{lem:an-t-strat}, which, as a side result, yields an even
more precise description of these $C$; see Remark~\ref{rem:convex}.

\begin{defn}\label{defn:rainbow}
Let $(S_i)_i$ be a t-stratification of a ball $B_0 \subseteq K^n$.
We define the \emph{rainbow} of $(S_i)_i$ to be the map $\rho\colon B_0 \to \RV\eq, x \mapsto \code{(\rv(x - S_i))_{i \le n}}$, where $\rv(x - S_i) = \{\rv(x - y) \mid y \in S_i\}$.
(Recall that such a code exists in $\RV\eq$ by stable embeddedness of $\RV$.)
\end{defn}

The rainbow is not uniquely determined, since we have to choose a code; however, this choice will never matter; therefore, we take the freedom to speak of ``the'' rainbow of a t-stratification.

\begin{rem}\label{rem:rainbow}
It is clear that the rainbow $\rho$ of $(S_i)_i$ refines $(S_i)_i$ in the sense of Definition~\ref{defn:refine} (and
Convention~\ref{conv:trans}), so any risometry $\phi\colon B_0 \to B_0$
respecting $\rho$ also respects $(S_i)_i$.
Vice versa, if $\phi$ is a risometry
respecting $(S_i)_i$, then for any $x\in B_0$ we have $\rv(x - S_i) = \rv(\phi(x) - S_i)$ and
hence $\phi$ respects $\rho$.
\end{rem}

\begin{prop}\label{prop:reflect}
Let $(S_i)_i$ be a t-stratification of $B_0 \subseteq K^n$ and let $\chi\colon B_0 \to \RV\eq$ be a definable map. Then the following are equivalent.
\begin{enumerate}
\item
$(S_i)_i$ reflects $\chi$.
\item
The rainbow of $(S_i)_i$ is a refinement of $\chi$.
\item
Any definable risometry $\phi\colon B_0 \to B_0$ respecting $(S_i)_i$ also respects $\chi$.
\end{enumerate}
\end{prop}

\begin{proof}
(2) $\Rightarrow$ (3) follows from Remark~\ref{rem:rainbow}.

(3) $\Rightarrow$ (1):
For every ball $B \subseteq B_0$, we have to show that $\tsp_B((S_i)_i, \chi) = \tsp_B((S_i)_i)$.
Let $(\alpha_x)_{x}$ be a translater of $(S_i)_i$ on $B$, with respect to any exhibition
of $\tsp_B((S_i)_i)$ (see Definition~\ref{defn:translater}).
Extending each $\alpha_x$ by the identity
on $B_0 \setminus B$ yields risometries $\alpha_x\colon B_0 \to B_0$ respecting $(S_i)_i$.
By (3), these risometries also respect $\chi$, hence $(\alpha_x)_{x}$ is also a translater
for $((S_i)_i, \chi)$.

(1) $\Rightarrow$ (2):
Let $\rho$ be the rainbow of $(S_i)_i$ and
suppose that for two points $y_1, y_2 \in B_0$, we have $\rho(y_1) = \rho(y_2)$ but
$\chi(y_1) \ne \chi(y_2)$.
Let $B :=\ball{y_1, \ge \val(y_1 - y_2)}$ be the smallest ball 
containing $y_1$ and $y_2$ and set $V := \tsp_B((S_i)_i)$. We may assume that
$y_1, y_2$ have been chosen such that $d := \dim V$ is maximal.

Choose an exhibition $\pi\colon B \to K^d$ of $V$ and a corresponding translater
$(\alpha_x)_{x \in \pi(B - B)}$ of $((S_i)_i, \chi)$.
Set $x_j := \pi(y_j)$ and let
$F_j := \pi\1(x_j)$ be the fiber containing $y_j$. Then for $y_1' := \alpha_{x_2-x_1}(y_1) \in F_2$,
we have $\chi(y_1') = \chi(y_1) \ne \chi(y_2)$. Moreover, since $\alpha_{x_2-x_1}$ respects
$(S_i)_i$, it also respects its rainbow (by Remark~\ref{rem:rainbow}),
i.e., $\rho(y_1') = \rho(y_1) = \rho(y_2)$.

Now set $B' := \ball{y'_1, \ge \val(y'_1 - y_2)} \subseteq B$. It remains to show
that $B' \cap S_d = \emptyset$ to get a contradiction to the maximality of $d$.
The set $T := S_d \cap F_2$ is finite but non-empty.
For any $y \in F_2$, we have $\rv(y - T) = \rv(y - S_d) \cap \rv(F_2 - F_2)$
by $V$-translatability on $B$,
thus $\rho(y_1')= \rho(y_2)$ implies
$\rv(y_1' - T) = \rv(y_2 - T)$. Now, Lemma~\ref{lem:fin&iso}
implies $(B' \cap F_2) \cap T = \emptyset$, which in turn implies
$B' \cap S_d = \emptyset$.
\end{proof}

\begin{rem}\label{rem:reflect}
The equivalence (1) $\iff$ (2) implies that for any two definable maps $\chi_1, \chi_2\colon B_0 \to \RV\eq$,
$(S_i)_i$ reflects the product $(\chi_1, \chi_2)$ if and only if it reflects $\chi_1$ and $\chi_2$ separately.
In particular, in the remainder of the article, we will use these two statements interchangeably.
\end{rem}

As mentioned at the beginning of this subsection, the following two results will be useful to ``enhance'' a given t-stratification; see the proofs of Lemma~\ref{lem:an-t-strat} and Proposition~\ref{prop:unif-cl} for applications.

\begin{cor}\label{cor:refine}
Let $(S_i)_i$ and $(S'_i)_i$ be two t-stratifications. Then the following are equivalent.
\begin{enumerate}
\item The rainbow of $(S'_i)_i$ refines the rainbow of $(S_i)_i$;
\item any definable map into $\RV\eq$ reflected by $(S_i)_i$ is also reflected by $(S'_i)_i$;
\item $(S'_i)_i$ reflects $S_j$ for each $j \le n$.
\end{enumerate}
\end{cor}

\begin{proof}
(1) $\iff$ (2) follows from Proposition~\ref{prop:reflect} (1) $\iff$ (2).

(3) $\Rightarrow$ (2) is obtained by using
Proposition~\ref{prop:reflect} (1) $\iff$ (3) to translate everything into statements about
which risometries respect what.

For (1) $\Rightarrow$ (3), note that (1) implies that
the rainbow of $(S'_i)_i$ refines $(S_i)_i$, so (3) follows from
Proposition~\ref{prop:reflect} (2) $\Rightarrow$ (1).
\end{proof}

If the above conditions hold, then $S_{\le i} \subseteq S'_{\le i}$ for all $i$
by Lemma~\ref{lem:I&D}. However, requiring $S_{\le i} \subseteq S'_{\le i}$ for all $i$
is not enough to imply the conditions of the corollary.

\begin{lem}\label{lem:small-changes}
Suppose that $X \subseteq K^n$ is a definable set of dimension $d$, that $\chi\colon X \to \RV\eq$ is a definable map, and that $(S_i)_i$, $(T_i)_i$ are
two t-stratifications of $K^n$, where
$(T_i)_i$ reflects $(S_i)_i$ and $\chi$. (Here, we extend the domain of $\chi$ to $K^n$ by sending
$K^n \setminus X$ to a single new element.) Then the following defines a
t-stratification which reflects $(S_i)_i$ and $\chi$ and which agree with $(S_i)_i$ outside of $X \cup T_{\le d - 1}$:
\[
S'_{\le i} := \begin{cases}
                T_{\le i}& \text{for } i < d\\
                S_{\le i} \cup X \cup T_{\le d - 1}& \text{for } i \ge d.
              \end{cases}
\]
\end{lem}


\begin{proof}
It is clear that
$\dim S'_{\le i} \le i$, so now consider a ball $B \subseteq S'_{\ge j}$;
we have to show that $((S'_i)_i, (S_i)_i, \chi)$ is $j$-translatable on $B$.
If $j \le d$, then we have $j$-translatability since $B \subseteq S'_{\ge j} = T_{\ge j}$
and $(T_i)_i$ reflects $((S_i)_i, \chi)$;
if $j \ge d + 1$, then $B \cap (X \cup T_{\le d - 1}) = \emptyset$ and
$S'_i \cap B = S_i \cap B$ for every $i$, so $j$-translatability follows from
$j$-translatability of $(S_i)_i$.
\end{proof}

To finish this subsection, we give some more properties of the fibers of a rainbow.

\begin{lem}\label{lem:rbp}   
Let $C$ be a fiber of the rainbow of a t-stratification $(S_i)_i$. Then 
either $C$ consists of a single element of $S_0$ or it is entirely contained in a ball $B \subseteq S_{\ge 1}$ and moreover, for any exhibition $\pi\colon B \to K^d$ of $\tsp_B((S_i)_i)$ and any fiber $F = \pi\1(y)$ (where $y \in \pi(B)$), $C \cap F$ is
exactly one fiber of the rainbow of the induced t-stratification $(S_{i + d} \cap F)_{i \le n - d}$ of $F$.
\end{lem}
\begin{proof}
Suppose that $x, x' \in C$. Since
$\rv(x - S_0) = \rv(x' - S_0)$, Lemma~\ref{lem:fin&iso} implies that there exists a risometry $K^n \to K^n$ fixing $S_0$ pointwise and sending $x$ to $x'$.
If $x \in S_0$, we obtain $x' = x$, so $C = \{x\} \subseteq S_0$; otherwise
by Lemma~\ref{lem:fin&iso} (2) (b), there exists a ball $B \subseteq S_{\ge1}$ containing both $x$ and $x'$. By doing this for all pairs $x, x' \in C$, we obtain a single ball $B \subseteq S_{\ge1}$ with $C \subseteq B$.

Let $\pi\colon B \to K^d$ be an exhibition of $V := \tsp_B((S_i)_i)$ and let $F = \pi\1(y)$ be a fiber (where $y \in \pi(B)$);
denote the induced t-stratification of $F$ by
$(S'_i)_{i \le n - d}$ and consider $x_1, x_2 \in F$. We have to check that $x_1$ and $x_2$ have the same image under the rainbow of $(S_i)_i$ iff they have the same image under the rainbow of $(S'_i)_i$.
It is clear that for any $i$, we have $\rv(x_1 - (S_i \setminus B)) = \rv(x_2 - (S_i \setminus B))$. This implies
$\rv(x_1 - S_i) = \rv(x_2 - S_i)$ for $i < d$ and it remains to verify that for $i \ge d$, we have
$\rv(x_1 - (S_i \cap B)) = \rv(x_2 - (S_i \cap B))$ 
iff $\rv(x_1 - S'_{i - d}) = \rv(x_2 - S'_{i - d})$.
This equivalence follows from $V$-translatability of $S_{i}$ on $B$. Indeed, translatability on the one hand implies
$\rv(x_j - S'_{i - d}) = \rv(x_j - S_{i}) \cap \rv(F - F)$ and hence ``$\Rightarrow$'';
on the other hand, we obtain $\rv(x_j - (S_{i} \cap B)) = \rv(x_j - S'_{i - d}) + \rv(\tilde V \cap (B - B))$ for any lift $\tilde V \subseteq K^n$ of $V$ (the sum on the right hand side is well-defined by Lemma~\ref{lem:vf}~(\ref{it:rv-sum})), and this implies ``$\Leftarrow$''.
\end{proof}

\begin{lem}\label{lem:rbp-sub}
Suppose that $(S_i)_i$ is a t-stratification and that $C \subseteq S_d$ is a fiber of the rainbow of $(S_i)_i$.
Then $C$ is sub-affine (cf.\ Definition~\ref{defn:sub-aff}) with $\affdir(C) = \tsp_B((S_i)_i)$ for any $B \subseteq S_{\ge d}$ with $B \cap C \ne \emptyset$.
\end{lem}

\begin{proof}
Lemma~\ref{lem:I&D} implies $\dim C = d$, hence it suffices to prove that $\affdir(C) \subseteq \tsp_B((S_i)_i)$.
We do an induction on $d$. If $d = 0$, then $C$ consists of a single element by Lemma~\ref{lem:rbp} and the claim is clear.
Otherwise, choose a ball $B' \subseteq S_{\ge 1}$ containing $C$ and choose an exhibition $\pi\colon B' \to K^{d'}$ of $V' := \tsp_{B'}((S_i)_i)$. By induction, for each $\pi$-fiber $F \subseteq B'$ and each suitable $B'' \subseteq F$, we have
$V'' := \affdir(C \cap F) = \tsp_{B''}((S_i \cap F)_i)$. By $V'$-translatability, $V''$ does not depend on the choice of $F$,
and by choosing $F$ and $B''$ such that $B'' = B \cap F \ne \emptyset$, we obtain $\tsp_B((S_i)_i) = V' + V''$.
Now choose any $x, x' \in C$, denote the corresponding $\pi$-fibers by $F$ and $F'$, respectively, and let $x''$
be the image of $x'$ in $F$ under a translater sending $F'$ to $F$. Then we have $\dir(x - x'') \in V''$
and $\dir(x'' - x') \in V'$ and hence $\dir(x - x') \in V' + V''$ by Lemma~\ref{lem:vf} (\ref{it:dir-sum}).
\end{proof}

\section{Fields with analytic structure satisfy Hypothesis~\ref{hyp:general}}
\label{sect:an->hyp}

Now it is time to prove that there do exist theories $\Tx$ satisfying Hypothesis~\ref{hyp:general}.
In \cite{CL.analyt}, Cluckers and Lipshitz introduce a notion of ``Henselian valued field with analytic structure'', which generalizes many older notions of analytic structures.
We will prove that any Henselian valued field of residue characteristic $0$
with analytic structure in that sense satisfies Hypothesis~\ref{hyp:general}.

Hypothesis~\ref{hyp:general}$_0$ follows directly from the results of \cite{CL.analyt}; we will give the arguments in Subsection~\ref{subsect:an->hyp0}.
Proving the Jacobian property needs more work. In dimension one, this is also done in \cite{CL.analyt} (up to some more minor differences in the definitions), but one essential ingredient to that proof is that any definable map $K \to \RV\eq$ can be refined in such a way that each fiber becomes a ball or a point. In higher dimensions, we will instead use our main theorem inductively to refine any definable map $K^n \to \RV\eq$ to the rainbow $\rho$ of a t-stratification $(S_i)_i$. We will see that each $\rho$-fiber is, up to a linear map and a risometry, a product of balls, and if
we are careful with the choice of $(S_i)_i$, we can moreover ensure that these risometries are analytic (Lemma~\ref{lem:an-t-strat}). This will allow the arguments from \cite{CL.analyt} proving the Jacobian property in dimension one to go through (Proposition~\ref{prop:an->jac}).

Since $\Tx = \THen$ is a special case of a theory of fields with analytic structure
(see Example~\ref{ex:an-pure} below), the proofs of this section in particular apply to
that case, i.e., $\THen$ satisfies Hypothesis~\ref{hyp:general}.

\subsection{The setting}
\label{subsect:an-set}

Let us fix the setting for the whole of Section~\ref{sect:an->hyp}. Concerning the notion of fields with analytic structure,
I will only repeat those
properties from \cite{CL.analyt} that are relevant to us, since the complete definition is somewhat technical.

Given a ``separated Weierstra\ss{} system'' $\mathcal A$ \cite[Definition~4.1.5]{CL.analyt},
one obtains a language $\LHenA$ \cite[beginning of Section~6.2]{CL.analyt} and
there is a notion of a ``separated analytic $\mathcal A$-structure'' on a valued field $K$ \cite[Definition~4.1.6]{CL.analyt} which turns $K$ into an $\LHenA$-structure.
The sorts of $\LHenA$ are $K$ and $\RV$. (More precisely, there are several $\RV$-like sorts, but these are all the same when the residue characteristic is zero.)
We let $\Lx$ be the union of $\LHenA$ and the remaining sorts of $\RV\eq$ (together with the canonical maps), and we let 
$\Tx$ be the $\Lx$-theory of all Henselian valued fields of equi-characteristic $0$ with separated analytic $\mathcal A$-structure.
(This is the same as the $\LHenA$-theory $\THenA$ defined at the beginning of \cite[Section~6.2]{CL.analyt},
except for the additional sorts in $\Lx$ and for the fact that $\THenA$ does not require the residue characteristic to be $0$.)

As always, $K$ will denote a model of $\Tx$.

\begin{ex}\label{ex:an-pure}
By \cite[Section~4.4, Example~(13)]{CL.analyt}, on any Henselian valued field
there exists an analytic structure whose definable sets are exactly the $\LHen$-definable sets.
Since we care about the language only up to interdefinability, this implies that the results of this section apply to $\LHen$ and $\THen$.
\end{ex}

Here is another concrete example of fields with analytic structure (a special case of \cite[Section~4.4~(1)]{CL.analyt}); a lot of more general examples are given in \cite[Section~4.4]{CL.analyt}.

\begin{ex}\label{ex:an}
Let $A := \bbZ[[t]]$ be equipped with the $t$-adic valuation (which we denote by $\vali$),
let
\[
T_{m} := A\langle \xi_1, \dots, \xi_m \rangle = \{\sum_{\nu \in \bbN^m} c_\nu\xi^\nu \mid c_\nu \in A, \lim_{|\nu| \to \infty} \vali(c_\nu) = \infty\}
\]
be the algebra of restricted power series (here, we use multi-index notation), and set
\[
S_{m,n} := T_{m}[[\rho_1, \dots, \rho_n]]
.
\]
As a language, take $\Lx := \LHen \dcup \bigdcup_{m,n} S_{m,n}$, where each element of $S_{m,n}$ is a symbol for an
$(m+n)$-ary function.

Now suppose that $K$ is a complete valued field of rank one and of residue characteristic $0$ (e.g., $K = \bbC((t))$).
Suppose moreover that $K$ extends $A$ as a valued ring; in particular, we identify $t \in A$ with an element of $\maxid$. Then each element of $S_{m,n}$ naturally defines a function $\valring^m \times \maxid^n \to \valring$.
This turns $K$ into an $\Lx$-structure (after extending these functions trivially to $K^{m+n}$)
and as such, $K$ is a valued field with analytic structure in the sense of \cite{CL.analyt}.
\end{ex}

\subsection{Fields with analytic structure satisfy Hypothesis~\ref{hyp:general}$_0$}
\label{subsect:an->hyp0}

We recall those definitions and results of \cite{CL.analyt} that we will need.

The following definition of ``b-minimality with centers'' is a combination of
\cite[Definitions~6.3.1 and 6.3.2]{CL.analyt} (slightly simplified, since we work in
equi-characteristic $0$). It differs from the 
notion of b-minimality from Definition~\ref{defn:b-min-eq} in two aspects.
First, ``with centers'' is a strengthening of Condition~(1),
and second, it only uses the sorts $K$ and $\RV$, whereas in Subsection~\ref{subsect:hyp}, we also allowed
the other sorts of $\RV\eq$.

\begin{defn}\label{defn:b-min}
An expansion $\Tx$ of $\THen$ is \emph{b-minimal with centers over $\RV$} if for every model $K \models \Tx$
and every set $A \subseteq K \cup \RV$, the following holds.
\begin{enumerate}
\item
For every $A$-definable set $X \subseteq K$, there exists an $A$-definable
map $\chi\colon X \surject Q \subseteq \RV^\ell$ (for some $\ell$)
and an $A$-definable map $c\colon Q \to K$ such that
for each $q \in Q$, the fiber $\chi\1(q)$ is of the form
$c(q) + \rvi\1(\xi)$ for some $\xi \in \RV$ (depending on $q$).
\item
There exists no surjective definable map from a subset of $\RV^\ell$
to an open ball $B \subseteq K$ (for any $\ell$).
\item
For every $A$-definable $X\subseteq K$ and $f\colon X \to K$, there exists
an $A$-definable map $\chi\colon X \surject Q \subseteq \RV^\ell$ (for some $\ell$) such that for each $q \in Q$,
$f\auf{\chi\1(q)}$ is either injective or constant.
\end{enumerate}
\end{defn}

\begin{lem}[{\cite[Theorem~6.3.7]{CL.analyt}}]\label{lem:an-qe}
Every $\Lx$-formula is, modulo $\Tx$, equivalent to an $\Lx$-formula without valued field quantifiers.
Moreover, $\Tx$ is b-minimal with centers over $\RV$.
\end{lem}

In \cite{CL.analyt}, the first statement of that lemma is formulated with $\LHenA$ instead of $\Lx$,
but it is easy to deduce the $\Lx$-version from the $\LHenA$-version.

The next result uses a definitory expansion
$\LHenA^{*}$ of the language
$\LHenA$ by certain functions $h_{m,n}$ which yield zeros of polynomials (cf.\ \cite[Definition~6.1.7]{CL.analyt});
we set $\Lx^* := \Lx \cup \LHenA^*$. A property of $\LHenA^*$ (and hence $\Lx^*$) we shall use is that for any function symbol
defining a function $f\colon K^m \times \RV^\ell \to K$ and any fixed $a \in K^m$, the function $f(a, \cdot)\colon \RV^\ell \to K$
has finite image.

\begin{lem}[{\cite[Theorem~6.3.8]{CL.analyt}}]\label{lem:an-term}
For any parameter set $A \subseteq K \cup \RV\eq$,
any $A$-definable function $K^n \to K$ can be written as $t(x, g(x))$ where $t$ is an $\Lx^*(A)$-term and $g$ is an $A$-definable function from $K^n$ to $\RV^\ell$ for some $\ell$.
\end{lem}

Again, \cite{CL.analyt} proves this for $\LHenA^*$ (and in particular $A \subseteq K \cup \RV$),
but the proof can easily be adapted to $\Lx^*$.

\begin{prop}\label{prop:an->hyp0}
Hypothesis~\ref{hyp:general}$_0$ holds in the setting described in Subsection~\ref{subsect:an-set}.
\end{prop}
\begin{proof}
Hypothesis~\ref{hyp:general} (\ref{it:st-emb}) follows from Lemma~\ref{lem:an-qe}
using that in $\Lx$, the only connection between $K$ and $\RV\eq$ is the map $\rvi$.
Hypothesis~\ref{hyp:general} (\ref{it:fin}) follows from Lemma~\ref{lem:an-term} for $n = 0$,
using the property of $\Lx^*$ described above that lemma.

Concerning Hypothesis~\ref{hyp:general} (\ref{it:strat1}),
let $K \models \Tx$, $A \subseteq K \cup \RV\eq$ and $X \subseteq K$ be given and set $A_0 := A \cap K$.
Then $X$ is $(A_0 \cup \{b\})$-definable for some tuple $b \in \RV^m$.
Applying Definition~\ref{defn:b-min}~(1) yields $(A_0 \cup \{b\})$-definable maps $\chi\colon X \surject Q \subseteq \RV^\ell$ and $c\colon Q \to K$. By writing $c(q) = c'(q, b)$ for some $A_0$-definable map $c'$ and
applying Hypothesis~\ref{hyp:general}~(\ref{it:fin}) to $c'$, we find that the image of $c$ is contained in a finite, $A_0$-definable set $S_0$.
Now suppose that for some ball
$B \subseteq K \setminus S_0$, we have $B \cap X \ne \emptyset$, say $x_0 \in B \cap X$;
we have to verify that $B \subseteq X$.
Using that $c(\chi(x_0)) \in S_0$, we obtain that
the map $x \mapsto \rvi(x - c(\chi(x_0)))$ is constant on $B$. By the description of the fibers of $\chi$,
this implies $B \subseteq \chi\1(\chi(x_0)) \subseteq X$.

To obtain Hypothesis~\ref{hyp:general}$_0$~($4''$), we only have to
generalize Definition~\ref{defn:b-min}~(3) from parameter sets in $K \cup \RV$ to parameter sets in $K \cup \RV\eq$.
Thus let $K$, $A$, $X$ and $f$ be given (where $A \subseteq K \cup \RV\eq$).
As before, set $A_0 := A \cap K$ and
choose $b \in \RV^m$ such that
$X$ and $f$ are $(A_0 \cup \{b\})$-definable. We find an $(A_0 \cup \{b\})$-definable map $\chi'\colon X \to \RV\eq$ that is as
desired, and it remains to modify it to make it $A$-definable. Again, we have
$\chi'(x) = \chi''(x, b)$ for some $A_0$-definable map $\chi''\colon X \times \RV^m \to \RV\eq$;
we define $\chi(x) := \code{y \mapsto \chi''(x, y)}$. This map is $A$-definable,
we may suppose that its range lies in $\RV\eq$ by Hypothesis~\ref{hyp:general} (\ref{it:st-emb}),
and since $\chi$ refines $\chi'$, $f$ has the desired property on each $\chi$-fiber.
\end{proof}

\subsection{Fields with analytic structure have the higher-dimensional Jacobian property}

As described at the beginning of this section, the strategy to obtain the Jacobian property
is to use our main result inductively and to describe fibers of rainbows more precisely.
Before we will do that, we will adapt two lemmas of \cite{CL.analyt}.
For all this, we will use the abstract notion of analytic functions introduced in \cite{CL.analyt}.
This notion works more smoothly when the field $K$ is algebraically closed, so in most of this subsection,
we will restrict to that case. We will see at the end that this is enough to obtain the Jacobian property
also for non-algebraically closed $K$.

So from now on assume that $K$ is algebraically closed. Then we have the following notions and results from \cite{CL.analyt}.
\begin{itemize}
\item 
A \emph{domain} is a particular kind of definable subset of $K^n$ \cite[Definition~5.2.2]{CL.analyt}; in particular, products of balls $B_i \subsetneq K$ are domains.\footnote{Strictly speaking, according to \cite[Definition~5.2.2]{CL.analyt}, only subsets of $\valring^n$ can be domains; however, via scaling one easily generalizes the notion to include bigger domains;
cf.\ also \cite[Remark~5.2.15]{CL.analyt}.} There are also notions of \emph{open} and \emph{closed} domains.
(A closed domain is essentially an abstract version of a rational domain in the sense of rigid geometry.)
\item
If $X$ is either an open domain or a closed domain, then
we have a well-defined \emph{ring of analytic functions} $\anfct(X)$ consisting of certain definable maps from $X$ to $K$ \cite[Definition~5.2.2 and Corollary~5.2.14]{CL.analyt}.
(In \cite{CL.analyt}, the ring $\anfct(X)$ is denoted by $\mathcal{O}^{\sigma}_{K}(X)$.)
\end{itemize}

We start by proving a higher-dimensional version of \cite[Lemma~6.3.15]{CL.analyt};
the idea of the proof is the same, with balls replaced by subsets of domains.

\begin{lem}\label{lem:Stw-an}
Suppose that $K$ is algebraically closed and let
$f\colon K^n \to K$ be an $A$-definable function (for some set of parameters $A \subseteq K \cup \RV\eq$).
Then there exists an $A$-definable map
$\chi\colon K^n \to RV\eq$ such that for each $\chi$-fiber $C \subseteq K^n$,
there exists an open domain $X \subseteq K^n$ containing $C$ such that $f\auf{C}$ is equal to $\tilde f\auf{C}$ for some analytic function $\tilde f\in \anfct(X)$.
\end{lem}

\begin{proof}
By Lemma~\ref{lem:an-term}, $f$ can be written as $f(x) = t(x, g(x))$, where
$t$ is an $\Lx^*(A)$-term and $g$ is an $A$-definable function with range in $\RV\eq$. By compactness, it therefore suffices to prove the lemma for the maps $x \mapsto t(x, \xi)$, where $\xi$ is a fixed element of the image of $g$.
We may assume that all subterms of $t$ are $K$-valued and we will omit $\xi$ from the notation.

The cases $t = x_{i_0}$ (for $i_0 \le n$) and $t = a$ (for $a \in A$) are trivial
(we only need to set $\chi(x_1, \dots, x_n) := (\rvi(x_1), \dots, \rvi(x_n))$ to ensure that
each $\chi$-fiber is contained in a domain),
so suppose that $t = h(t'_1, \dots, t'_\ell)$ where $h$ is a function symbol in
$\Lx^*$. Induction yields definable maps $\chi'_j\colon K^n \to \RV\eq$ for the terms $t'_j$.
We set $\chi(x) := \big(\chi'_1(x), \dots, \chi'_\ell(x), \rvi(t'_1(x)), \dots, \rvi(t'_\ell(x))\big)$.

Fix $q_j \in \im \chi'_j$ and $\xi_j \in \RV$ for each $j$ and consider
$C := \chi\1(q_1, \dots, q_\ell, \xi_1, \dots, \xi_\ell)$.
For each $j$, we have $C \subseteq C'_j := (\chi'_j)\1(q_j)$ and
induction yields a domain $X'_j \supseteq C'_j$
and an analytic function $\tilde t'_j \in \anfct(X'_j)$ with $t'_j\auf{C'_j} = \tilde t'_j\auf{C'_j}$.

We define $X := \bigcap_j X''_j$, where $X''_j = \{x \in X'_j \mid \rvi(\tilde t'_j(x)) = \xi_j\}$ if
$\xi_j \ne 0$ and $X''_j = X'_j$ if $\xi_j = 0$. This is an open domain and it contains $C$.
Next, we define $\tilde t := h(\tilde t''_1, \dots, \tilde t''_\ell)$, where
$\tilde t''_j := \tilde t'_j$ if $\xi_j \ne 0$ and $\tilde t''_j := 0$ if $\xi_j = 0$.
With this definition, we have
$t\auf{C} = \tilde t\auf{C}$, and it remains to check that $\tilde t$ is analytic on $X$.
If $h$ is $y_1 + y_2$ or $y_1 \cdot y_2$, then this is clear; otherwise, this follows
from \cite[Lemma~6.3.11]{CL.analyt}, since
our construction ensures that $\rvi(\tilde t''_j)$ is constant on $X$ for each $j$.
\end{proof}

The next lemma is a variant of \cite[Lemma~6.3.9]{CL.analyt}.

\begin{lem}\label{lem:639}
Suppose that $K$ is algebraically closed and that $g \in \anfct(\valring)$ is an analytic function
such that $g'(x) \in \valring$ for every $x \in \valring$ and $\res(g'(x))$ is constant.
Then, for any $x_1, x_2 \in \valring$ with $x_1 \ne x_2$, we have
\[
\vali(g(x_1) - g(x_2) - g'(0)\cdot(x_1 - x_2)) > \vali(x_1-x_2)
\]
(where $g'$ denotes the derivative of $g$).
\end{lem}

\begin{proof}
First of all, note that we may suppose that $g'(0) = 0$; otherwise, replace $g(x)$ by $g(x) - g'(0)x$. Thus the assumption becomes that $g'(x) \in \maxid$ for all $x$ and we want to prove that $\vali(g(x_1) - g(x_2)) > \vali(x_1 - x_2)$.

Recall how analytic functions on $\valring^m$ are defined in \cite{CL.analyt}: an analytic structure on the field $K$ yields a
ring of abstract power series $A^\dag_{m,0}(K) \subseteq \valring[[\xi_1, \dots, \xi_m]] \otimes_{\valring} K$ and a ring homomorphism $f \mapsto f^\sigma$ from $A^\dag_{m,0}(K)$
to the ring of maps from $\valring^m$ to $K$; $\anfct(\valring^m)$ is the image of $A^\dag_{m,0}(K)$ under $\sigma$. Let $f = \sum_i a_i \xi^i \in A^\dag_{1,0}(K)$ be a preimage of $g$ under $\sigma$, i.e., $g = f^\sigma$.

The assumption $g'(\valring) \subseteq \maxid$ implies $a_i \in \maxid$ for all $i \ge 1$.
Indeed, otherwise, choose any $r \in K$ with $\vali(r) = -\min_{i \ge 1} \vali(a_i)$ (the minimum exists by \cite[Remark~4.1.10]{CL.analyt}) and consider the series $r\cdot f' \in \valring[[\xi]]$. By \cite[Definition~4.1.2 (iv)]{CL.analyt}, $\res(r\cdot f')$ is a polynomial and
since $g'(\valring) \subseteq \maxid$ and $\vali(r) \ge 0$, it is $0$ everywhere on $k$. However, this implies that all its coefficients are $0$ and hence $\vali(ra_i) > 0$ for all $i \ge 1$, contradicting the choice of $r$.

Using Weierstra\ss{} division (\cite[Definition~4.1.3]{CL.analyt}), we find a power series $h(x, y) \in A_{2,0}(K)$
such that $f(x) - f(x + y) = h(x, y)\cdot y$. Since the constant term of the left hand side cancels, all its coefficients lie in $\maxid$. This then also holds for $h$ and hence
(using  \cite[Remark~4.1.10]{CL.analyt} again) the range of $h^\sigma$ is contained in $\maxid$. This implies the lemma, since $g(x_1) - g(x_2) = h^\sigma(x_1, x_2-x_1) \cdot (x_2 - x_1)$.
\end{proof}

Recall that if a t-stratification $(S'_i)_i$ reflects another t-stratification $(S_i)_i$, then it also reflects any definable map into $\RV\eq$ reflected by $(S_i)_i$ (Corollary~\ref{cor:refine}).
The following lemma morally consists of two separate statements, namely (1) given
any $(S_i)_i$, one can find $(S'_i)_i$ that reflects $(S_i)_i$ and that is ``piecewise analytic'';
and (2) any rainbow of any t-stratification has fibers of a simple form.
In the lemma, the piecewise analyticity is formulated in terms of the fibers of the rainbow of $(S'_i)_i$ (which is what we really need in the application). In Remark~\ref{rem:convex}, we will give a separate formulation of statement (2).

We call a map from an open or closed domain to $K^d$ analytic if each of its coordinates is analytic.

\begin{lem}\label{lem:an-t-strat}
Suppose that we are in the setting of Subsection~\ref{subsect:an-set} and that $\Tx$ additionally satisfies Hypothesis~\ref{hyp:general}$_{n-1}$.
Suppose moreover that $K$ is algebraically closed and that $(S_i)_i$
is an $A$-definable t-stratification of $K^n$ for some set of parameters $A \subseteq K \cup \RV\eq$. Then there exists an $A$-definable t-stratification $(S'_i)_i$ reflecting $(S_i)_i$ such that for any fiber $C$
of the rainbow of $(S'_i)_i$ with $C \subseteq S'_n$, we have the following.
\begin{enumerate}
\item $C$ is an open domain.
\item There exist open balls $B_1, \dots, B_{n} \subseteq K$ and an
analytic bijection $\phi \colon B_1 \times \dots \times B_{n} \to C$ that can be written as the composition of a risometry and a matrix from $\GL_n(\valring)$.
\end{enumerate}
\end{lem}

\begin{rem}\label{rem:liso}
By Remark~\ref{rem:risoGLO}, the order of the composition in (2) doesn't matter,
and the composition of several such maps can again be written as a composition of a single risometry and a single matrix from $\GL_n(\valring)$.
\end{rem}

\begin{proof}[Proof of Lemma~\ref{lem:an-t-strat}]
The proof consists of two parts. First, we will construct $(S'_i)_i$ in such a way that each fiber
of the rainbow is (essentially) the graph of an analytic function. (This will be needed for all fibers and not just for the ones contained in $S'_n$). In the second part, we will show that this is enough to imply the lemma.

\smallskip

\textbf{Part 1:}

Recall that if $C$ is a fiber of the rainbow of a t-stratification, then $C$ is sub-affine 
by Lemma~\ref{lem:rbp-sub}, and by Lemma~\ref{lem:sub-prop}, if $\pi \colon K^n \surject K^j$ exhibits $\affdir(C)$
(with $j = \dim C$), then $C$ is the graph of a function $c \colon \pi(C) \to K^{n-j}$.

By downwards induction on $d$, we will prove that we can find an $A$-definable t-stratification $(S'_i)_i$ reflecting $(S_i)_i$ with the following property. If $C \subseteq S'_{j}$ is a fiber of the rainbow of $(S'_i)_i$ for some $j \ge d$ and $\pi, c$ are as above, then there exists an open domain $X \subseteq K^{j}$ containing $\pi(C)$ such that $c$ is the restriction
to $\pi(C)$ of an analytic function $X \to K^{n-j}$.

By induction, we may assume that $(S_i)_i$ itself has the above property for $d+1$. To obtain a t-stratification $(S'_i)_i$ which has the property for $d$, we proceed as follows. For any rainbow fiber $C \subseteq S_{d}$ and for $\pi$, $c$ as above, we apply Lemma~\ref{lem:Stw-an} to find a $\code{C}$-definable map $\tilde{\chi} \colon \pi(C) \to \RV\eq$ such
that each $\tilde{\chi}$-fiber $F \subseteq \pi(C)$ is contained in an open domain $X \subseteq K^d$ such that $c\auf{F}$ is the restriction of an analytic function on $X$. Composing $\tilde{\chi}$ with the projection $C \to \pi(C)$ yields a definable map $C \to \RV\eq$; we take the product of these maps
for all exhibitions $\pi$ of $\affdir(C)$ to obtain a single map $C \to \RV\eq$ and we do this for all $C \subseteq S_d$;
using compactness, this yields an $A$-definable map $\chi \colon S_{d} \to \RV\eq$.

Let $(T_i)_i$ be a t-stratification reflecting $(S_i)_i$ and $\chi$
(this exists by Theorem~\ref{thm:main}, since we are assuming Hypothesis~\ref{hyp:general}$_{n-1}$).
Applying Lemma~\ref{lem:small-changes} to $(S_i)_i$, $(T_i)_i$, $X := S_d$, and $\chi$,
we obtain a t-stratification $(S'_i)_i$ reflecting both $(S_i)_i$ and $\chi$ with
$S'_{i} \subseteq S_{i}$ for $i \ge d$. We claim that this $(S'_i)_i$ has the desired properties.

Let $C' \subseteq S'_{j}$ be a fiber of the rainbow of $(S'_i)_i$, where $j \ge d$. By Corollary~\ref{cor:refine}, there exists a fiber $C$ of the rainbow of $(S_i)_i$ entirely containing $C'$ and since $S'_{j} \subseteq S_{j}$, we have $C \subseteq S_{j}$ and hence (using Lemma~\ref{lem:rbp-sub}) $\affdir(C') = \affdir(C)$. Let $\pi\colon K^n \surject K^{j}$ be an exhibition of $\affdir(C)$ and let $c\colon \pi(C) \to K^{n-j}$ and $c'\colon \pi(C') \to K^{n-j}$ be the functions whose graphs are $C$ and $C'$, respectively; note that $c'$ is simply the restriction of $c$ to $\pi(C')$.

If $j > d$, then by induction $c$ is the restriction of an analytic function to $\pi(C)$, so $c'$ is the restriction of the same function to $\pi(C')$. If $j = d$, then the fact that $(S'_i)_i$ reflects the map $\chi$ implies that $\chi$ is constant on $C'$ and hence that the corresponding map $\tilde{\chi} \colon \pi(C) \to \RV\eq$ is constant on $\pi(C')$. By construction of $\tilde{\chi}$, this implies the claim.

\smallskip
\textbf{Part 2:}

To avoid some special cases, we assume that $S'_0 \ne \emptyset$ (if $S'_0 = \emptyset$, we can replace it by $\{0\}$).
Fix a fiber $C \subseteq S'_n$ of the rainbow of $(S'_i)_i$.
For $d \le n$, consider the following statement. There exists a coordinate projection $\pi\colon K^n \surject K^d$ and a $\lambda \in \Gamma$ such that the following conditions hold:
\begin{enumerate}
\item  
  For each $q \in \pi(C)$, $C \cap \pi\1(q)$ is contained in an open ball $B_q \subseteq \pi\1(q)$ of radius $\lambda$;
\item
  for each $q_1, q_2 \in \pi(C)$, there exists a risometry $B_q \to B_{q'}$ respecting 
  the rainbow of $(S'_i)_i$;
\item
  for each $q \in \pi(C)$, $\tsp_B((S'_i)_i)$ is exhibited by $\pi$, where $B$ is the open subball of $K^n$ of radius $\lambda$ containing $B_q$ (in particular  $\dim \tsp_B((S'_i)_i) = d$ and hence $B_q \cap S'_d \ne \emptyset$);
\item
  $\pi(C)$ is an open domain;
\item
 there exists a set $R = B_1 \times \dots \times B_d$, where each $B_i \subseteq K$ is an open ball,
and an analytic bijection
$\phi\colon R \to \pi(C)$ which can be written as the composition of a risometry and a matrix from $\GL_d(\valring)$.
\end{enumerate}
In the case $d = n$, (4) and (5) together imply the lemma.
(In that case, we are not interested in (1) -- (3), which are a bit pathological since $B_q$ is supposed to be a zero-dimensional ball.
However, everything makes sense even in that case, using that a ball of any radius in $K^0$ is equal to (the one-point set) $K^0$ itself.)

For $d = 0$, this statement follows from $S'_0 \ne \emptyset$. Indeed, the latter implies that
$C$ is contained in a ball $B \subsetneq K^n$ (which is needed for (1)) and that for $B$ sufficiently big,
$\tsp_B((S'_i)_i) = \{0\}$ (which is needed for (3)).

To finish the proof, we will show that if the statement holds for some $d < n$, then there also exists a $d'>d$ for which it holds.

To simplify notation, we assume that $\pi$ projects to the first $d$ coordinates.
By (3), $(S'_i)_i$ induces a t-stratification of $B_q$ for each $q \in \pi(C)$ and by Lemma~\ref{lem:rbp},
$B_q \cap C$ is a fiber of the rainbow of that t-stratification and 
it is contained in a ball $B'_q \subseteq B_q \setminus S'_d$.
We may assume this $B'_q$ to be a maximal open ball in $B_q \setminus S'_d$.
Since $B_q \cap S'_d \ne \emptyset$ (which follows from (3)), we have $B'_q = s_q + (\{0\}^d \times \rv\1(\xi_q))$ for some $s_q \in S'_d \cap B_q$ and some $\xi_q \in \RV^{(n-d)}$. We first fix choices of $s_q$ and $\xi_q$ for one single $q \in \pi(C)$ and then use risometries from (2) to obtain analogous elements $s_{q'}$ and $\xi_{q'}$ for all other $q' \in \pi(C)$; this ensures that
$\xi := \xi_q$ does not depend on $q$ and that all the $s_q$ are contained in a single fiber $\tilde{C} \subseteq S'_d$ of the rainbow of $(S'_i)_i$. Set $\lambda' := \vRV(\xi)$.

By (2), $V := \tsp_{B'_q}((S'_i \cap B_q)_i) \subseteq k^{n-d}$ does not depend on $q$. Since
$B'_q \cap S'_d = \emptyset$, we have $\dim V \ge 1$.
Set $d' := d + \dim V$, let $\rho\colon K^{n-d} \surject K^{d' - d}$ be an exhibition of $V$, assume
that $\rho$ is the projection to the first $d' - d$ coordinates of $K^{n-d}$, and
let $\pi' := \rho \circ \pi\colon K^{n} \surject K^{d'}$ be the projection to the first $d'$ coordinates. We will now verify that the above statement holds for these choices of $d', \lambda', \pi'$; we denote the corresponding conditions by (1') to (5').

(1'), (2'), and (3') are clear. By (3) and Lemmas~\ref{lem:rbp-sub} and \ref{lem:sub-prop}, $\tilde{C}$ is the graph of a function
$c\colon \pi(\tilde{C}) \to K^{n-d}$ which, by part 1 of the proof, is equal to the restriction of an analytic function; from this, we will now easily deduce (4') and (5'). First note that since $\pi(C)$ is an open domain, also $c\auf{\pi(C)}$ is analytic. (In fact,
one could also check that $\pi(\tilde{C}) = \pi(C)$.) Define $e\colon \pi(C) \to K^{d' - d}, q \mapsto \rho(c(q))$.

By $V$-translatability of $C \cap B_q$ on $B'_q$, we have $\pi'(C \cap B_q) = \pi'(B'_q)$, which in turn is equal to $\{q\} \times (e(q) + \tilde{B})$, where $\tilde{B} := \rho(\rv\1(\xi)) \subseteq K^{d'-d}$. This implies that $\pi'(C)$ is an open domain, and for $R,\phi$ as in (5) and
$R' := R \times \tilde{B}$, we obtain an analytic bijection
$\phi'\colon R' \to \pi'(C), (r, x) \mapsto (\phi(r), e(\phi(r)) + x)$. By Lemma~\ref{lem:sub-prop},
the map $\psi\colon \pi(C) \times \tilde{B} \to \pi'(C), (q, x) \mapsto (q, e(q) + x)$ is the composition of a risometry and an element of $\GL_{d'}(\valring)$, so using Remark~\ref{rem:liso}, the same holds for $\phi' = \psi \circ (\phi \times \id_{\tilde{B}})$, which finishes the proof.
\end{proof}

\begin{rem}\label{rem:convex}
By skipping the whole first part of this proof, we obtain that for any fiber $C \subseteq S_n$ of the rainbow of any t-stratification $(S_i)_i$ (in any valued field satisfying Hypothesis~\ref{hyp:general}$_0$), there exists a definable bijection $\phi \colon B_1 \times \dots \times B_{n} \to C$ that can be written as the composition of a risometry and a matrix in $\GL_n(\valring)$ (and where each $B_i \subseteq K$ is an open ball). Moreover, it is not difficult to modify the proof to obtain a similar statement for fibers $C \subseteq S_d$ with $d < n$.
\end{rem}

Now we can finally prove that fields with analytic structure have the Jacobian property in any dimension.

\begin{prop}\label{prop:an->jac}
In the setting of Subsection~\ref{subsect:an-set},
$\Tx$ has the Jacobian property in the sense of Definition~\ref{defn:jac};
in particular, Hypothesis~\ref{hyp:general} holds.
\end{prop}

\begin{proof}
We first reduce the general case to the case where $K$ is algebraically closed.

Let $K$ be a model of $\Tx$ (as usual) and let $f\colon K^n \to K$ be $A$-definable.
As in the proof of Lemma~\ref{lem:Stw-an}, using Lemma~\ref{lem:an-term} we may assume that
$f(x) = t(x, \xi)$, where
$t$ is an $\Lx^*(A)$-term and $\xi \in \RV\eq$ is $A$-definable.
By \cite[Theorem~4.5.11 (i)]{CL.analyt}, the analytic structure on $K$ has a (unique)
extension to an analytic structure on the algebraic closure $\Kacl$ of $K$ (in the same language). In particular, our term $t(x, \xi)$ also defines a function $\Kacl^n \to \Kacl$. Assuming that the proposition holds in $\Kacl$, we find an $A$-definable map $\tilde{\chi}\colon \Kacl^n \to \RVacl\eq$ such that $t$ has the Jacobian property on each $n$-dimensional $\tilde{\chi}$-fiber; here $\RVacl\eq$ denotes the $\RV\eq$-sorts corresponding to $\Kacl$. By elimination of valued field quantifiers in the language $\Lx$ (Lemma~\ref{lem:an-qe}), we may assume that $\tilde{\chi}(x) = (\rvi(t_1(x)), \dots, \rvi(t_k(x)))$ for some $\Lx(A)$-terms $t_i$. Indeed, after elimination of the quantifiers, the map $x \mapsto (\rvi(t_i(x)))_i$ is a refinement of $\tilde{\chi}$, where $t_i$ runs over all $\Kacl$-valued $\Lx(A)$-terms appearing in $\tilde{\chi}$.

Now $\tilde{\chi}$ is given by a tuple of terms, so it restricts to a map $\chi\colon K^n \to \RV\eq$. If a fiber $\chi\1(q)$ has dimension $n$, then so has $\tilde\chi\1(q)$, so 
$t$ has the Jacobian property on $\tilde\chi\1(q)$ and in particular on $\chi\1(q)$.

\medskip

We now have to prove the proposition in the case where $K$ is algebraically closed.
As announced, we prove the ``Jacobian property up to dimension $n$''
by induction over $n$. For $n = 0$, there is nothing to show,
so now suppose $n \ge 1$ and suppose that $f\colon K^n \to K$ is an
$A$-definable map.

Lemma~\ref{lem:Stw-an} yields an $A$-definable map $\chi\colon K^n \to \RV\eq$
such that for each $\chi$-fiber $C$, there exists an open domain $X \supseteq C$ and an analytic function $g \in \anfct(X)$ such that $f\auf{C} = g\auf{C}$. Recall that such an analytic function $g$ has 
a well-defined \emph{Jacobian} at every $x \in X$: $(\Jac g)(x) = (\del g/\del x_1, \dots, \del g/\del x_n) \in K^n$.
Suppose that $\dim C = n$.
By refining $\chi$ (and using Lemma~\ref{lem:loc-dim}), we may assume that $\dim_x C = n$ for every $x \in C$. This implies that for $x \in C$, $(\Jac g)(x)$ is determined by $f\auf{C}$, so we can further refine $\chi$
in such a way that $\rv((\Jac f\auf{C})(x)) \in \RV^{(n)}$ is constant on each $\chi$-fiber $C$ of dimension $n$.

We apply Theorem~\ref{thm:main} to $\chi$ and obtain an $A$-definable t-stratification $(S_i)_i$ of $K^n$
reflecting $\chi$ (the prerequisites of Theorem~\ref{thm:main} hold by Proposition~\ref{prop:an->hyp0} and the induction hypothesis). Then we apply Lemma~\ref{lem:an-t-strat} (using the induction hypothesis again); the resulting t-stratification $(S'_i)_i$ still reflects $\chi$ by Corollary~\ref{cor:refine}.
Let $\rho$ be the rainbow of $(S'_i)_i$ and let $C$ be an $n$-dimensional $\rho$-fiber.
Then $C$ is an open domain, $f$ is analytic on $C$ (since $f\auf{C}$ is the restriction of an analytic function on a larger domain), $\rv((\Jac f\auf{C})(x))$ is constant on $C$, and there exists an analytic bijection $\phi = \psi \circ M \colon R \to C$ where $R = B_1 \times \dots \times B_n$ for some open balls $B_i \subseteq K$, where $\psi$ is a risometry, and where $M \in \GL_n(\valring)$.

To finish the proof, we will check that $f$ has the Jacobian property on $C$, i.e, either $f\auf{C}$ is constant or for any $x, x' \in C$ with $x \ne x'$, we have
\[\tag{$*$}
\vali(f(x) - f(x') - \langle z, x - x'\rangle) > \val(z) + \val(x - x')
\]
for some $z \in K^n$ not depending on $x, x'$. In fact, we will prove ($*$) for
$z = (\Jac f\auf{C})(x')$, which is enough by Remark~\ref{rem:jac} and since $\rv(\Jac f\auf{C})$ is constant.

We may assume that $\rv((\Jac f\auf{C})) \ne 0$, since otherwise, $f\auf{C}$ is constant and we are done. Let $x, x' \in C$, $x \ne x'$ be given and
define $\eta\colon \valring \to R$, $y \mapsto y\cdot  \phi\1(x) + (1 - y)\cdot \phi\1(x')$. The remainder of the proof consists in pulling back the problem on $C$ to a problem on $\valring$ using 
$\theta := \phi \circ \eta\colon \valring \to C$; the problem on $\valring$ then follows from Lemma~\ref{lem:639}. All this is straightforward, but we give some details.

The map $\theta$ is obviously analytic and an easy computation shows that
$\rv(\frac{\theta(y) - \theta(y')}{y-y'})\in \RV^{(n)}$ is constant for $y, y' \in \valring$, $y \ne y'$. Indeed, we have $\rv(\theta(y) - \theta(y')) = \rv(M(\eta(y)) - M(\eta(y')))$
since $\psi$ is a risometry, and $\rv(\frac{M(\eta(y)) - M(\eta(y'))}{y-y'})$ is constant since $M \circ \eta$ is linear.
Moreover, constantness of $\rv(\frac{\theta(y) - \theta(y')}{y-y'})$ implies that $\rv(\Jac\theta)$ is also constant and that these two values are equal.
In particular, by plugging $y = 1$ and $y' = 0$ (the preimages of $x$ and $x'$) into this equality, we obtain 
\[\tag{$**$}
\val(x - x' - (\Jac \theta)(0)) > \val(x - x')
.
\]

Set $g := f \circ \theta\colon \valring \to K$; it is clear that $g$ is analytic.
Using the chain rule $g' = \langle \Jac f\auf{C}, \Jac \theta\rangle$ and that $\rv(\Jac f\auf{C})$ and $\rv(\Jac \theta)$ are constant, we obtain $\vali(g'(y)) \ge \val(\Jac f\auf{C}) + \val(\Jac \theta)$ and $\vali(g'(y) - g'(y')) > \val(\Jac f\auf{C}) + \val(\Jac \theta)$ 
for any $y, y' \in \valring$ with $y \ne y'$.
Therefore, we can apply Lemma~\ref{lem:639} to the function $y \mapsto \frac{g(y)}{r}$, where $r \in K$ is any element with $\vali(r) = \val(\Jac f\auf{C}) + \val(\Jac \theta)$; this yields
\[
\vali(g(y) - g(y') - g'(0)\cdot (y-y')) > \val(\Jac f\auf{C}) + \val(\Jac \theta) + \vali(y - y').
\]
For $y = 1$, $y' = 0$, and $z = (\Jac f\auf{C})(x')$, this becomes
\[\tag{$*\mathord{*}*$}
\vali(f(x) - f(x') - g'(0)) > \val(z) + \val(\Jac \theta)
.
\]
Now ($**$) already implies that ($*$) and ($*\mathord{*}*$) have the same right hand side,
and it remains to verify that
\[
\vali(\langle z, x - x'\rangle - g'(0)) > \val(z) + \val(x - x').
\]
Indeed, ($**$) implies
\[
\vali(\langle z, x - x'\rangle - \langle z, (\Jac \theta)(0) \rangle) > \val(z) + \val(x - x'),
\]
which is what we want, since $g'(0) = \langle z, (\Jac \theta)(0)\rangle$.
\end{proof}

\section{Algebraic results}
\label{sect:alg}

Up to now, for a t-stratification $(S_i)_i$ we know that the sets $S_{\le i}$
are closed in the valued field topology. However, in a purely algebraic setting, it
would be natural to require the
sets $S_{\le i}$ to be Zariski closed. We will now show that indeed this can be achieved (Corollary~\ref{cor:alg}).
In fact, we will first prove a more general result (Proposition~\ref{prop:unif-cl}) in arbitrary theories
satisfying Hypothesis~\ref{hyp:general}.

\subsection{Getting closed sets \texorpdfstring{$S_{\le i}$}{S\_\042\144i}}

In this subsection we assume that $\Tx$ satisfies Hypothesis~\ref{hyp:general}.

In Proposition~\ref{prop:unif-cl}, we
introduce a set $\Delta$ of formulas which can be thought of as defining the closed sets of a topology
(although the conditions on $\Delta$ will be weaker) and we prove that any t-stratification $(S_{i})_i$ can be enhanced
to a t-stratification $(S'_{i})_i$ in such a
way that each $S'_{\le i}$ is closed in this sense. For this to be possible, we only need that taking the closure of a definable set doesn't increase its dimension. Here, ``enhancing'' means that any
map into $\RV\eq$ reflected by $(S_{i})_i$ is also reflected by $(S'_{i})_i$. By Corollary~\ref{cor:refine},
this is equivalent to: $(S'_{i})_i$ reflects $S_j$ for each $j$.

To be able to work uniformly for all models of $\Tx$, we introduce a uniform notion of dimension.

\begin{defn}\label{defn:unif-dim}
For an $\Lx$-formula $\phi$ whose free variables live in the valued field sort,
set $\dim\phi := \max_{K \models \Tx} \dim\phi(K)$.
\end{defn}

\begin{prop}\label{prop:unif-cl}
In the following, all $\Lx$-formulas have $n$ free valued field variables,
and ``$\phi \rightarrow \psi$'' means ``$\Tx \vdash \forall x\,(\phi(x) \rightarrow \psi(x))$''.

Suppose that we have a family $\Delta$ of $\Lx$-formulas with the following properties:
\begin{enumerate}
\item $\Delta$ is closed under disjunctions and contains $\bot$.
\item
For each $\Lx$-formula $\phi$, there exists a minimal formula $\phi^\Delta \in \Delta$ with $\phi \rightarrow \phi^\Delta$;
minimal means: for any other $\psi \in \Delta$ with
$\phi \rightarrow \psi$, we have $\phi^\Delta \rightarrow \psi$.
\item
For each $\Lx$-formula $\phi$, we have $\dim \phi^\Delta = \dim\phi$.
\end{enumerate}
Suppose moreover that $(\phi_i)_{0 \le i \le n}$ is a
tuple of formulas defining a t-stratification in every model of $\Tx$.
Then we can find a tuple of formulas
$(\phi'_i)_i$ which, in every model, defines a t-stratification reflecting the sets
defined by the formulas $\phi_i$ and such that for each $i$,
$\phi'_0 \vee \dots \vee \phi'_i$ is equivalent to a formula in $\Delta$.
\end{prop}

\begin{proof}
For any formula $\phi$, we set $\del \phi := \phi^\Delta \wedge \neg\phi$. Note that
using $\phi^\Delta \vee \psi^\Delta \in \Delta$, one obtains $(\phi \vee \psi)^\Delta \to (\phi^\Delta \vee \psi^\Delta)$
and hence $\del (\phi \vee \psi) \to (\del \phi \vee \del\psi)$.

We write $\phi_{\le i}$ for $\phi_0 \vee \dots \vee \phi_i$, and similarly for $\phi'_{i}$ and $\psi_{i}$ (which will be introduced below).

Suppose that for some given $d \in \{0, \dots, n\}$, $(\phi_i)_i$ satisfies $\dim \del \phi_{\le i} \le d$
for all $i$. From this, we construct a 
t-stratification $(\phi'_{i})_i$ reflecting $(\phi_i)_i$ and satisfying $\dim \del \phi'_{\le i} \le d - 1$. Applying this repeatedly yields the proposition (where $\dim \phi \le -1$
will mean $\Tx \vdash \neg\exists x\phi(x)$). So let $d$ be given as above.

For $i$ from $n$ to $0$, recursively define $\delta_i := \del(\phi_{\le i} \vee \delta_{i+1} \vee \dots \vee \delta_n)$. Inductively, we get $\dim \delta_i \le d$.
Set $\delta := \bigvee_{i=0}^n \delta_i$, choose any t-stratification $(\psi_i)_i$ reflecting $((\phi_i)_i, \delta)$,
and apply Lemma~\ref{lem:small-changes} to $X = \delta(K)$, a constant map $\chi\colon X \to \RV\eq$, $(S_i)_i = (\phi_i(K))_i$ and $(T_i)_i = (\psi_i(K))_i$.
We claim that for the resulting t-stratification $(\phi'_i(K))_i := (S'_i)_i$, we have
$\del \phi'_{\le i}  \to (\phi'_{\le d-1})^\Delta$ (where $\phi'_{\le -1} = \bot$); in particular, we obtain
$\dim \del \phi'_{\le i} \le d - 1$.

For $i \le d - 1$, the claim is clear, so suppose $i \ge d$.
The formula $\phi'_{\le i} = \phi_{\le i} \vee \delta \vee \psi_{\le d - 1}$
is equivalent to $\bigvee_{j \le i} (\phi_{\le j} \vee \delta_j \vee \dots \vee \delta_n) \vee \psi_{\le d - 1}$. Each formula $\phi_{\le j} \vee \delta_j \vee \dots \vee \delta_n$ is equivalent to a formula in $\Delta$
by definition of $\delta_j$, so we get $\del \phi'_{\le i} \to \del \psi_{\le d - 1}$, and $\del\psi_{\le d - 1}$ is equal to $\del \phi'_{\le d - 1}$.
\end{proof}

\subsection{Algebraic strata}
\label{subsect:alg}

In the pure valued field language $\LHen$, we can apply Proposition~\ref{prop:unif-cl}
to the family of formulas that are conjunctions of polynomial equations and in this way obtain
t-stratifications where the sets $S_{\le i}$ are Zariski closed.
This yields a version of the main theorem that can almost be formulated in a purely
algebraic language---only
almost, since in the definition of t-stratification, we require the straighteners
to be definable. (Of course, this condition can simply be omitted, but this weakens
the result.) Nevertheless, we take the opportunity to present a setting
that is as algebraic as possible.

Fix a Noetherian integral domain $A$ of characteristic $0$. We set
$\Lx := \LHen(A)$ and $\Tx := \THen \cup \{$positive atomic diagram of $A$ in $\Lring\}$;
in other words, models $K$ of $\Tx$ are Henselian valued fields of equi-characteristic $0$
together with a ring homomorphism $A \to K$.
We fix $n \in \bbN$ and let $\Delta$ be the set of conjunctions of polynomial equations in $n$ variables
with coefficients in $A$. For any model $K \models \Tx$, by considering the sets $\phi(K), \phi \in \Delta$
as closed, we obtain the Zariski topology (``over $A$'') on $K^n$. More precisely,
formulas in $\Delta$ correspond to Zariski closed subsets of the scheme $\bbA_A^n$, and
our topology on $K^n$ is the one which the Zariski topology on $\bbA_A^n$ induces on the $K$-valued points $\bbA_A^n(K)$.

Note that for $\phi \in \Delta$, we have two notions of dimension: the one given in
Definition~\ref{defn:unif-dim} and the algebraic one, where we consider $\phi$ as
a variety over $A$. However, by considering an algebraically closed model of $\Tx$, we see
that the two notions of dimension coincide.

Given an arbitrary $\Lx$-formula $\phi$, it is clear that there exists a well-defined
``Zariski closure'' of $\phi$, i.e., a minimal formula
$\phi^\Delta \in \Delta$ implied by $\phi$. To be able to apply Proposition~\ref{prop:unif-cl}, it remains to check that
$\phi^\Delta$ has the same dimension as $\phi$.
This has been proven in \cite{Dri.dimDef} or \cite{CL.mot} for example, but in slightly
different contexts than ours, so let us quickly repeat the proof from \cite{Dri.dimDef}. We first work
in a fixed model $K$.

\begin{lem}\label{lem:alg-dim}
For every $\emptyset$-definable set $X \subseteq K^n$, there exists a formula $\psi \in \Delta$ such that $X \subseteq \psi(K)$ and $\dim \psi = \dim X$.
\end{lem}

\begin{proof}
In this proof, we will write $\rvi^\ell\colon K^\ell \to \RV^\ell$ for the cartesian power of $\rvi$
(in contrast to the map $\rv\colon K^\ell \to \RV^{(\ell)}$ mainly used in the remainder
of the article).

By quantifier elimination (see e.g.\ \cite{CL.analyt}, Theorem~6.3.7), $X$ is of the form
\[
X = \{x \in K^n \mid \big(\rvi(f_1(x)), \dots, \rvi(f_\ell(x))\big) \in \Xi\}
= f\1((\rvi^\ell)\1(\Xi))
\]
where $f = (f_1, \dots, f_\ell)$ is an $\ell$-tuple of polynomials with coefficients in $A$ and $\Xi \subseteq \RV^\ell$ is $\emptyset$-definable.
The statement of the lemma is preserved by finite unions, so
we can do a case distinction on whether $f_i(x) = 0$ or not for each $i$;
in other words, $X$ is of the form
\[
X = \psi(K) \cap f\1((\rvi^\ell)\1(\Xi))
\]
for $\psi \in \Delta$, $f$ as above and $\Xi \subseteq (\RV \setminus \{0\})^\ell$.

Write $\Kacl$ for the algebraic closure of $K$.
We may assume that $\psi(\Kacl)$ is the Zariski closure of $X$ in $\Kacl^n$.
In particular, $X$ contains a regular point $x$ of $\psi(\Kacl)$,
i.e., on a Zariski-neighborhood of $x$, $\psi(\Kacl)$ is defined by
$n-\dim \psi$ polynomials and the Jacobian matrix at $x$ of this tuple of polynomials
has maximal rank.

Now $(\rvi^\ell)\1(\Xi)$ is open in the valuation topology, so in that topology,
there is a neighborhood $U \subseteq \psi(K)$ of $x$ which is contained in $X$.
Using the implicit function theorem
and regularity at $x$, we find a coordinate projection $\pi\colon K^n \surject K^{\dim \psi}$ such that $\pi(U)$ contains a ball in $K^{\dim \psi}$. This implies $\dim X = \dim \psi$.
\end{proof}

Now we make the result uniform for all models of $\Tx$.

\begin{lem}\label{lem:alg-dim-univ}
For every $\Lx$-formula $\phi$ in $n$ valued field variables,
there exists a formula $\psi \in \Delta$
with $\dim \psi = \dim \phi$ and $\phi(K) \subseteq \psi(K)$ for all models $K \models \Tx$.
\end{lem}

\begin{proof}
For each $K$ separately, Lemma~\ref{lem:alg-dim} yields a formula $\psi_K \in \Delta$ with
$\phi(K) \subseteq \psi_K(K)$ and $\dim \psi_K = \dim \phi(K)$. By compactness,
there exists a finite disjunction $\psi$ of some of the $\psi_K$
such that $\phi(K) \subseteq \psi(K)$ for all $K$. Since
\[\dim \psi \le \max_K \dim \psi_K = \max_K \dim \phi(K) = \dim \phi
,
\]
we are done.
\end{proof}

Now Proposition~\ref{prop:unif-cl} can be applied to the Zariski topology and we get
t-stratifications such that each set $S_{\le i}$ is defined by a
conjunction of polynomials (uniformly for all models).
Moreover, using that being a t-stratification is first order
in the sense of Corollary~\ref{cor:pos-char},
the same t-stratification also works in models of a finite subset of $\Tx$.
Here is the precise result.

\begin{cor}\label{cor:alg}
Let $A$ be a Noetherian integral domain of characteristic $0$, $\Lx = \LHen(A)$, and $\Tx$
the theory of Henselian valued fields $K$ of equi-characteristic $0$ together with a
ring homomorphism $A \to K$. 
For every $\Lx$-formula $\chi$ defining a map $\chi_K\colon K^n \to \RV\eq$ (for any $K \models \Tx$),
there exists a finite subset $\Tx_0 \subseteq \Tx$ and
formulas $(\phi_i)_i$ such that:
\begin{itemize}
\item
Either $\dim \phi_i = i$ (in the sense of Definition~\ref{defn:unif-dim}) or
$\phi_i = \bot$.
\item
The disjunction $\phi_0 \vee \dots \vee \phi_d$ is equivalent to a conjunction of polynomial equations
with coefficients in $A$.
\item
For every model $K \models \Tx_0$, $(\phi_i(K))_i$ is a t-stratification reflecting
$\chi_K$.
\end{itemize}
\end{cor}

(As in Corollary~\ref{cor:pos-char}, we assume that models of $\Tx_0$ are valued fields
for the statements to make sense.)

Here is a an algebraic formulation of Corollary~\ref{cor:alg}; by a ``subvariety of $\bbA^n_A$'', we simply mean a reduced (not necessarily irreducible) subscheme.
Since the notion of a
definable map to $\RV\eq$ is not so algebraic, we instead formulate the theorem for a finite
family $(X_\nu)_\nu$ of subvarieties of $\bbA^n_A$.

\begin{thm}\label{thm:alg}
Let $A$ be a Noetherian integral domain of characteristic $0$ and let $X_\nu$ be subvarieties of $\bbA^n_A$ for $\nu = 1, \dots, \ell$. Then there exists an integer $N \in \bbN$ and
a partition of $\bbA^n_A$ into
subvarieties $S_i$ with the following properties:
\begin{enumerate}
\item
$\dim S_i = i$ or $S_i = \emptyset$
\item
Each $S_{\le i}$ is a closed subvariety of $\bbA^n_A$.
\item
For every Henselian valued field $K$ over $A$ of residue characteristic either $0$ or
at least $N$,
$(S_i(K))_i$ is a t-stratification of $K^n$ reflecting the family of sets
$(X_\nu(K))_\nu$, i.e., for every $d \le n$ and every ball $B \subseteq S_{\ge d}(K)$,
the tuple $\big(S_d(K), \dots, S_n(K), X_1(K), \dots, X_\ell(K)\big)$ is $d$-translatable on $B$
(see Definition~\ref{defn:trans} or \ref{defn:t-strat} and Convention~\ref{conv:trans}).
\end{enumerate}
\end{thm}

In the next section, we will prove that a t-stratification $(S_i)_i$ as in Theorem~\ref{thm:alg} induces a classical Whitney stratification $(S_i(\bbC))_i$ of $\bbC^n$ (for any ring homomorphism $A \to \bbC$), and similarly for $\bbC$ replaced by $\bbR$ (see Theorem~\ref{thm:t->whit}). In particular, this implies that each $S_i$ is smooth over the fraction field of $A$.

We conclude this section with an algebraic formulation of Corollary~\ref{cor:isotyp} about how the risometry type can vary in a uniform family.

\begin{cor}\label{cor:alg-isotyp}
Let $A$ be a Noetherian integral domain of characteristic $0$,
let $Q$ be any affine variety over $A$,
and let $X_\nu$ be subvarieties of $\bbA^n_Q$
for $\nu = 1, \dots, \ell$.
Then there exists an integer $N \in \bbN$ and
algebraic maps $f_1,\dots, f_m\colon Q \to \bbA^1_A$
such that 
for every Henselian valued field $K$ over $A$ of residue characteristic either $0$ or
at least $N$, we have the following.

Given $q \in Q(K)$, write $X_{\nu,q} =  X_\nu \times_Q \operatorname{spec} K$ for the fiber of $X_\nu$ over $q$
and consider $X_{\nu,q}(K)$ as a subset of $K^n$.
If two elements $q, q' \in Q(K)$ satisfy $\rvi(f_\mu(q)) = \rvi(f_\mu(q'))$ for all $\mu$,
then there exists a risometry $K^n \to K^n$ sending $X_{\nu,q}(K)$ to $X_{\nu,q'}(K)$
for each $\nu$.
\end{cor}

\begin{proof}
We fix an embedding $Q \inject \bbA^\ell_A$.
Applying Corollary~\ref{cor:isotyp} yields an integer $N$ and a formula
$\phi$ such that for every $K$ as above, $\phi$ defines a
map $\phi_K \colon Q(K) \to \RV\eq$ such that $\phi_K(q) = \phi_K(q')$ implies existence of
a risometry as above for $q, q'\in Q(K)$.
By quantifier elimination \cite[Theorem~6.3.7]{CL.analyt},
we may refine the map defined by $\phi$ to a map of the form
$q \mapsto (\rvi(f_1(q)), \dots, \rvi(f_m(q)))$ for some polynomials $f_i$;
this implies the claim.
\end{proof}

\section{Obtaining classical Whitney stratifications}
\label{sect:whit}

The main result of this section is that the existence of t-stratifications implies
the existence of classical Whitney stratifications. More precisely, a non-standard
model of $\bbR$ or $\bbC$ can be considered as a valued field, and we will see that any definable partition
in the standard model that induces a t-stratification in the non-standard model
is already a Whitney stratification. We will start by proving that t-stratifications
satisfy a kind of analogue of Whitney's Condition~(b) (Corollary~\ref{cor:whit-b}).
This needs the following additional natural Hypothesis on the language $\Lx$ (and the theory $\Tx$).

\begin{hyp}\label{hyp:orth}
In this section, we require that the residue field is orthogonal to the value group,
i.e., in any model of $\Tx$, any definable set $X \subseteq k^n \times \Gamma^m$ is a finite union of
sets of the form $Y_i \times Z_i$, for some definable sets $Y_i \subseteq k^n$ and $Z_i \subseteq \Gamma^m$.
\end{hyp}

\begin{prop}\label{prop:an->orth}
The theory of any Henselian valued field $K$ with analytic structure in the sense
of \cite{CL.analyt} satisfies Hypothesis~\ref{hyp:orth}.
\end{prop}
\begin{proof}
We work with the language $\Lx$ introduced in Subsection~\ref{subsect:an-set}.
By quantifier elimination (Lemma~\ref{lem:an-qe}), any definable subset
of $\RV^n$ can be defined
in the restriction to $\RV$ of $\Lx$.
To that restricted language, add the sorts $k$ and $\Gamma$ and a splitting $\RV \setminus \{0\} \to k\mult$ of the sequence
$k\mult \inject \RV \setminus \{0\} \surject \Gamma$ (such a splitting corresponds to an angular component map $K \to k$ of the valued field); then it becomes interdefinable with the language $\Lx'$ consisting of $k$ with the ring language and $\Gamma$ with the language $\{0, +, -, <\}$ of ordered abelian groups (where $\RV$ is identified with $k \times \Gamma$).
In particular, any set $X \subseteq k^n \times \Gamma^m$ definable in our original language
is also $\Lx'$-definable. Since $\Lx'$ contains no connection between $k$ and $\Gamma$, the proposition follows.
\end{proof}

At some point, we will use the following easy consequence of the above hypothesis.

\begin{rem}\label{rem:dclGamma}
For any parameter set $A \subseteq k$, we have $\acl(A) \cap \Gamma = \acl(\emptyset) \cap \Gamma$ (where $\acl$ is
the algebraic closure in the model theoretic sense). Using the order on $\Gamma$, we
get the same with $\acl$ replaced by the definable closure.
In particular, if $\Lx = \LHen$ and $K$ is either real closed or algebraically closed,
then $\Gamma$ is a pure divisible ordered abelian group and the only finite,
$A$-definable subsets of $\Gamma$ are $\emptyset$ and $\{0\}$.
\end{rem}

\subsection{An analogue of Whitney's Condition~(b)}

Our main theorem about the existence of t-stratifications only speaks about the dimension
of translatability spaces. The following theorem additionally (partially) specifies their
direction. The analogue of Whitney's Condition~(b) will then be a corollary.
(Recall from Definition~\ref{defn:dir} that $\dirRV \colon \RV^{(n)} \to \bbP^{n-1} k$ is the map induced
by $\dir\colon K^n \to \bbP^{n-1} k$.)

\begin{thm}\label{thm:kegel}
Suppose that the language $\Lx$ and the $\Lx$-theory $\Tx$ satisfy Hypotheses~\ref{hyp:general} and \ref{hyp:orth} and that $K$ is a model of $\Tx$.
Let $\chi \colon K^n \to \RV\eq$ be a definable map and let $x \in K^n$
be any point. Let $\Xi \subseteq \RV^{(n)} \setminus \{0\}$ be the set of those $\xi$ such that
$\chi$ is not $\dirRV(\xi)$-translatable on the ball $B := x + \rv\1(\xi)$.
Then $\vRV(\Xi)$ is finite.
\end{thm}


\begin{rem}\label{rem:Xi-def}
In general, $\Xi$ is not definable. However, we can choose a t-stratification
$(S_i)_i$ reflecting $\chi$ and refine $\chi$ to $((S_i)_i, \chi)$; after this modification,
$\Xi$ is definable (over $x$ and the parameters of the original $\chi$) by Lemma~\ref{lem:trans-def}.
\end{rem}

\begin{proof}[Proof of Theorem~\ref{thm:kegel}]
By Remark~\ref{rem:Xi-def}, we may assume that $\Xi$ is definable.
Without loss, fix $x = 0$ and suppose for contradiction that $\vRV(\Xi)$ is infinite.
By orthogonality of the value group and the residue field, there exists
a one-dimensional $V \subseteq k^n$ such that
the subset $\Xi_0 := \{\xi \in \Xi \mid \dirRV(\xi) \in V\}$ is already infinite.
Choose a lift $\tilde{V} \subseteq K^n$ of $V$ and consider the map $\chi'$ obtained from $(\chi, \tilde{V})$
via Convention~\ref{conv:trans}.
For $\xi \in \Xi_0$, $\chi$ is not $V$-translatable on the ball $B :=\rv\1(\xi)$;
on the other hand, $\tilde{V} \cap B \ne \emptyset$, so $\tsp_{B}(\tilde{V}) = V$, which implies that $\chi'$ is not translatable at all on $\rv\1(\xi)$.
In particular, if $(S_i)_i$ is a t-stratification reflecting $\chi'$ (which exists by Theorem~\ref{thm:main}), then
we have $\rv\1(\xi) \cap S_0 \ne \emptyset$ for all $\xi \in \Xi_0$, which contradicts $S_0$ being finite.
\end{proof}

Note that the only way we used Hypothesis~\ref{hyp:general} in this proof is to apply Theorem~\ref{thm:main}
to $\chi$ and $\chi'$.

In the classical version of Whitney's Condition~(b), one has two sequences of points in two different strata $S_d$ and $S_j$ with $d < j$, and both sequences converge
to the same point in $S_d$. In the valued field version,
each sequence is replaced by a single point, and ``converging to the same point in $S_d$'' is replaced
by ``lying in a common ball $B \subseteq S_{\ge d}$''. In the proof of Proposition~\ref{prop:t->whit},
we will see how this implies the classical Condition~(b) via non-standard analysis.

\begin{cor}\label{cor:whit-b}
Assume Hypotheses~\ref{hyp:general} and \ref{hyp:orth}.
Let $(S_i)_i$ be a $\emptyset$-definable t-stratification of a $\emptyset$-definable ball $B_0 \subseteq K^n$,
let $B \subseteq B_0$ be a sub-ball, and let $d$ be maximal with $B \subseteq S_{\ge d}$. Then 
there exists a finite $\code{B}$-definable set $M \subseteq \Gamma$ such that the following holds.
For any $j > d$, any $x' \in B \cap S_d$, and any $y' \in B \cap S_j$, if $\val(x' - y') \notin M$, then
$\dir(x' - y') \in \tsp_{B'}((S_i)_i)$, where $B' \subseteq S_{\ge j}$ is a ball containing $y'$.
\end{cor}

\begin{proof}
Let $\pi\colon B \to K^d$ be an exhibition of $\tsp_B((S_i)_i)$.
Choose any $z \in \pi(B)$ and any $x \in \pi\1(z) \cap S_d$, and apply
Theorem~\ref{thm:kegel} to $(S_i)_i$ and $x$.
This yields a finite set $\vRV(\Xi) \subseteq \Gamma$, which is $x$-definable by Remark~\ref{rem:Xi-def}. Doing this for all $x \in \pi\1(z) \cap S_d$
and taking the union of the (finitely many) corresponding sets $\vRV(\Xi)$ yields a finite,
$(\code{B}, z)$-definable set which we denote by $M_z$.
For any other $z' \in \pi(B)$, Lemma~\ref{lem:translater} yields a risometry $\alpha\colon B \to B$ sending
$\pi\1(z)$ to $\pi\1(z')$; extending $\alpha$ by the identity on $B_0 \setminus B$, we get
a risometry which shows that $M_z = M_{z'}$; hence $M := M_z$ is $\code{B}$-definable.

Now let $x' \in B \cap S_d$ and $y' \in B \cap S_j$ be given with $\val(x' - y') \notin M$ and set $z = \pi(x')$. Since $\lambda := \val(y' - x') \notin M_z$, $(S_i)_i$ is $\dir(x' - y')$-translatable
on $B_1 := \ball{y', > \lambda}$ and hence also on $B' \subseteq B_1$.
\end{proof}

\subsection{The classical Whitney conditions}
\label{subsect:defn-whit}

We now recall the definition of Whitney stratifications; see e.g.\ \cite{BCR.realGeom}
for more details.
We will consider Whitney stratifications both over $k = \bbR$ and $k = \bbC$,
in a semi-algebraic resp.\ algebraic setting.
A Whitney stratification is a partition of $k^n$ into certain kinds of manifolds.
In the case $k = \bbR$, we will work with Nash manifolds and also with a weakening
of that notion.

\begin{defn}\label{defn:nash}
A \emph{Nash manifold} is a $C^\infty$-sub-manifold of $\bbR^n$ (for some $n$), which is
$\Lring$-definable (or, equivalently, which is semi-algebraic). By a
\emph{$C^1$-Nash manifolds}, we mean a $C^1$-sub-manifold of $\bbR^n$ which is
$\Lring$-definable.
\end{defn}

Note that by ``$M$ is a sub-manifold of $k^n$'' we mean that also the inclusion
map $M \inject k^n$ is in the corresponding category, i.e., either $C^1$ or $C^\infty$ (but we do not require $M$ to be
closed in $k^n$).
All our manifolds will be sub-manifolds of some $k^n$ in this sense (for $k$ either $\bbR$ or $\bbC$);
this will not always be written explicitly.

In the case $k = \bbC$, we will only have one notion of manifolds, namely algebraic sub-manifolds of $\bbC^n$.
Note that this is in perfect analogy to the case $k = \bbR$; if we simply replace
$\bbR$ by $\bbC$ in Definition~\ref{defn:nash}, then
``definable'' means constructible instead of semi-algebraic; moreover
``differentiable'' should now be read as ``complex differentiable''. Thus
in that case, both kinds of manifolds introduced in Definition~\ref{defn:nash} simply become
algebraic manifolds.

In the remainder of the section, we will treat $k = \bbR$ and $k = \bbC$ simultaneously, and we will write 
``Nash/algebraic manifolds'' or ``$C^1$-Nash/algebraic manifolds'' (depending on the
notion of manifold we want to consider in the case $k = \bbR$).

We will not require our manifolds to be connected, but if they are not, then each
connected component has to have the same dimension.

For a $C^1$-Nash/algebraic manifold $M \subseteq k^n$ and a point $x \in M$,
there is a well-defined notion of tangent space $T_xM \subseteq k^n$ of $M$ at $x$.
Such a space can be seen as an element of the corresponding Grassmanian $\grass_{n,\dim M}(k)$ and as such, it makes sense to speak of limits of sequences of such spaces.

\newcounter{enumsave}

\begin{defn}\label{defn:whit}
Let $k$ be either $\bbR$ or $\bbC$.
A \emph{Whitney stratification} of $k^n$ is a partition of $k^n$
into Nash/algebraic manifolds $(S_i)_{0 \le i \le n}$ with the following properties.
(As always, we write $S_{\le i}$ for $S_0 \cup \dots \cup S_i$.)
\begin{enumerate}
\item\label{it:dimi} For each $i$, either $\dim S_i = i$ or $S_i = \emptyset$.
\item\label{it:closi} Each set $S_{\le i}$ is topologically closed in the analytic topology.
\item\label{it:wai} Each pair $S_i, S_j$ with $i < j$ satisfies \emph{Whitney's Condition~(a)},
i.e., for any element $u \in S_i$ and any sequence $v_\mu \in S_j$ converging to
$u$, if $\lim_{\mu \to \infty} T_{v_\mu}S_j$ exists, then
\[
T_uS_i \subseteq \lim_{\mu \to \infty} T_{v_\mu}S_j
.
\]
\item\label{it:wbi} Each pair $S_i, S_j$ with $i < j$ satisfies \emph{Whitney's Condition~(b)},
i.e., for any two sequences $u_\mu \in S_i$, $v_\mu \in S_j$ both converging to the same element $u \in S_i$, if both $\lim_{\mu \to \infty} T_{v_\mu}S_j$ and $\lim_{\mu\to\infty}k\cdot(u_\mu - v_\mu)$ exist,
then
\[
\lim_{\mu\to\infty}k\cdot(u_\mu - v_\mu) \subseteq \lim_{\mu \to \infty} T_{v_\mu}S_j
.
\]
\setcounter{enumsave}{\theenumi}
\end{enumerate}
We will say that $(S_i)_{i}$ is a \emph{$C^1$-Whitney stratification} if
it is a partition of $k^n$
into $C^1$-Nash/algebraic manifolds satisfying the above conditions (\ref{it:dimi}) -- (\ref{it:wbi}).
\end{defn}

In fact, it is known that anyway (4) implies (3); we will prove (3) separately nevertheless,
since the argument is short and elegant.

Often, one additionally requires that the topological closure of any connected component
of any $S_j$ is the union of some connected components of some of the $S_i$, $i < j$. However, once one knows how to obtain Whitney stratifications in our sense,
it is easy to also obtain this additional condition.

\subsection{Transfer to the Archimedean case}
\label{subsect:t->whit}

Let $k$ be either $\bbR$ or $\bbC$.
We will consider $k$ as a structure in the language
$\Labsring := \Lring \cup \{|\cdot|\}$, where $|\cdot| \colon k \to \bbR_{\ge 0} \subseteq k$
is the absolute value. (Of course, in the case $k = \bbR$, $|\cdot|$ is already $\Lring$-definable.)
Fix $K$ to be a (non-principal) ultra-power of $k$ with index set $\bbN$;
this will be the non-standard model of $k$ we will be working in. (In fact, any
$\aleph_1$-saturated elementary extension of the $\Labsring$-structure $k$ would do.)

The image in $K$ of any $u \in k$ (under the canonical embedding) is denoted by $\nsu$.
Similarly, for any set $X \subseteq k^n$, the ultra-power of $X$,
considered as a subset of $K^n$, will be denoted by $\nsX$. (In particular, $\ns{k} = K$ and $\nsR \subseteq K$.)

Define
\[
\valring := \{x \in K \mid \exists (u \in \bbR)\, |x| < u\}
;
\]
this is a valuation ring, turning $K$ into a valued field which is Henselian and of equi-characteristic $0$. The maximal ideal is
\[
\maxid = \{x \in K \mid \forall (u \in \bbR_{>0})\, |x| < u\}
,
\]
the residue field is $k$, and $\res\colon \valring \surject k$
is simply the standard part map.

Using the absolute value on $k$, we can define the Euclidean norm on $k^n$, which we denote by $|\cdot|_2\colon k^n \to \bbR_{\ge0}$; this also induces an ``Euclidean norm'' $|\cdot|_2\colon K^n \to \nsR_{\ge0}$,
and this Euclidean norm induces a topology on $K^n$, given by the subbase $\{x \in K^n \mid |x - a|_2 < r\}$,
$a \in K^n$, $r \in \nsR_{>0}$. This topology is the same as the 
valuation topology on $K^n$, since for any $a \in K^n$, any $\lambda \in \Gamma$, and any
$r \in \nsR_{>0}$ with $\vali(r) > \lambda$, we have
\[
\ball{a, >\vali(r)} \subseteq \{x \in K^n \mid |x - a|_2 < r\}
\subseteq \ball{a, >\lambda};
\]
note that we continue to use the notations $\ball{a, >\lambda}$, $\ball{a, \ge\lambda}$ for balls in the valuative sense.

Let $X \subseteq k^n$ be any definable set.
Any sequence $(u_\mu)_{\mu \in \bbN}$ with $u_\mu \in X$ and $\lim_{\mu \to \infty} u_\mu = u \in k^n$
represents an element of $\res\1(u) \cap \nsX$ in the ultra-product; 
vice versa, any element of $\res\1(u) \cap \nsX$ can be represented by a sequence in $X$ converging to $u$.
If $(u_\mu)_{\mu}$ is such a converging sequence,
we will write $[u_\mu]$ for the corresponding element of $\res\1(\lim_\mu u_\mu)$.
We will also use this notation with more complicated expressions inside the square brackets;
the index variable of the sequence will always be $\mu$. Note that square brackets
commute with definable maps as follows. If, in addition to $X$ and $u_\mu$ as above, we have
definable $Y \subseteq k^m$ and $f\colon X \to Y$, then $[f(u_\mu)] = \ns{\!f}([u_\mu])$, where
$\ns{\!f}$ denotes the corresponding map $\nsX \to \nsY$.

The following lemma is almost trivial, but it is the main tool which makes
the transfer between $k$ and $K$ work. (Note that we implicitly use that the
two different topologies on $K$ coincide.)

\begin{lem}\label{lem:nsa}
For any definable set $X \subseteq k^n$ and any element $u \in k^n$, the following are equivalent:
\begin{enumerate}
\item
$u$ lies in the topological closure of $X$;
\item
$\nsu$ lies in the topological closure of $\nsX$;
\item
$\nsX \cap \res\1(u)$ is non-empty.
\end{enumerate}
\end{lem}
\begin{proof}
(1) $\iff$ (2) follows from definability of being in the topological closure.
(1) $\iff$ (3): If $u$ lies in the closure of $X$, then any sequence $v_\mu \in X$
converging to $u$ yields an element $[v_\mu] \in X \cap \res\1(u)$.
Vice versa, if $[v_\mu] \in X \cap \res\1(u)$, then we may assume $v_\mu \in X$ for all $\mu$ and thus
$u = \lim_\mu v_\mu$ lies in the closure of $X$.
\end{proof}

Note that the equivalence (2) $\iff$ (3) does not hold if one replaces $\nsX$ by an arbitrary definable subset of $K^n$; the point is that $\nsX$ is $\Labsring$-definable
and using only parameters from the image of $k$ in $K$.

On the other hand, Lemma~\ref{lem:nsa} also applies to definable subsets of varieties (instead of subsets of $k^n$);
in particular, we will apply it in the Grassmanians $\grass_{n,d}$.


Now we can formulate the main proposition of this section.

\begin{prop}\label{prop:t->whit}
Let $k$ be either $\bbR$ or $\bbC$ and let $K$ be as in the beginning of Subsection~\ref{subsect:t->whit}.
Suppose that $(S_i)_i$ are $\Lring$-definable subsets of $k^n$
such that $(\nsS_i)_i$ is a t-stratification of $K^n$. Then $(S_i)_i$
is a $C^1$-Whitney stratification of $k^n$ (in the sense of Definition~\ref{defn:whit}).
\end{prop}

\begin{proof}
In this proof, we will use the letters $u, v$ for elements of $k^n$ and
$x, x'$ for elements of $K^n$.
We have to prove conditions (1) -- (4) of Definition~\ref{defn:whit}
and that each $S_i$ is a $C^1$-Nash/algebraic manifold.

Since dimension is definable, (\ref{it:dimi}) follows from the corresponding
property of $\nsS_i$. (To obtain a definition of dimension which works both in $k^n$ and
in $K^n$, we can replace the valuative ball in Definition~\ref{defn:dim} by an Euclidean ball.) 

Using Lemma~\ref{lem:nsa}, closedness of
$\nsS_{\le i}$ implies (\ref{it:closi}) and moreover that for any $u \in S_d$,
$\res\1(u)$ is a subset of $\nsS_{\ge d}$; in particular, $(\nsS_i)_i$ is
$d$-translatable on $\res\1(u)$.

Fix $u \in S_d$ and set $B := \res\1(u)$ and $V_u := \tsp_B((\nsS_i)_i)$. We claim that
\[\tag{$\diamond$}
\affdir(B \cap \nsS_d) = V_u\]
(cf.\ Definition~\ref{defn:sub-aff}). For dimension reasons, it suffices to verify ``$\subseteq$'',
i.e., for any $x, x' \in B \cap \nsS_d$, we have $\dir(x - x') \in V_u$.
To prove this, choose an exhibition $\pi\colon K^n \surject K^d$ of
$V_u$, set $F := \pi\1(\pi(\nsu)) \cap \nsS_d$, and
apply Lemma~\ref{lem:nsa} to $F \setminus \{\nsu\}$. Since $\nsu$ does not lie in the closure
of $F \setminus \{\nsu\}$, we obtain $F \cap B = \{\nsu\}$, i.e., $\pi$-fibers of $\nsS_d$ in $B$
consist of a single element.
Now $V_u$-translatability of $\nsS_d$ on $B$ implies the claim. 

Next, we prove that $S_d$ is a $C^1$-manifold and that $V_u$ is the tangent space at $u$
for every $u \in S_d$ and $V_u$ as above. First note that each point $u \in S_d$ has a neighborhood $U \subseteq k^n$ such that for a suitable coordinate projection $\pibar\colon k^n \surject k^d$,
$\pibar$ induces a bijection $U \cap S_d \to \pibar(U)$. (Indeed, this is a first order statement
and it holds in $K$.) We will use $\pibar$ as a chart of $S_d$ around $u$. To prove that its inverse
$(\pibar\auf{U \cap S_d})\1$ is $C^1$ and that $V_u$ is the tangent space at $u$ (for every
$u \in S_d$), it suffices to verify the following.
For any $u \in S_d$ and any two sequences $v_\mu, v'_{\mu} \in S_d$ with $\lim_\mu v_\mu = \lim_\mu v'_\mu = u$ and $v_\mu \ne v'_\mu$, if $\lim_\mu k\cdot (v_\mu - v'_\mu)$ exists (in $\grass_{n,1}(k)$),
then $\lim_\mu k\cdot (v_\mu - v'_\mu) \subseteq V_u$.
So suppose that such $u, v_\mu, v_{\mu'}$ are given. Working in $\grass_{n,1}$,
we have
\[
\lim_\mu k\cdot (v_\mu - v'_\mu)
= \res([k\cdot (v_\mu - v'_\mu)])
= \res(K\cdot ([v_\mu] - [v'_\mu]))
.
\]

Now $[v_\mu], [v'_\mu] \in \res\1(u) \cap \nsS_d$ implies $\dir([v_\mu] - [v'_\mu]) \in V_u$
by ($\diamond$)
and hence 
$\res(K\cdot ([v_\mu] - [v'_\mu])) \subseteq V_u$.

In the case $k = \bbC$, we just proved that $S_d$ is $C^1$ in the sense of complex
differentiation, so in that case, we obtain that $S_d$ is an algebraic manifold.

Sending a point $u \in S_d$ to its tangent space $T_uS_d$ is a 
definable map $S_d \to \grass_{n,d}(k)$; transferring this to $K$ yields a notion
of tangent space of $\nsS_d$ at any $x \in \nsS_d$; we denote that tangent space (which is a sub-space of $K^n$) by $T_x\nsS_d$.
Fix $x \in \nsS_d$.
By definition, if $x' \in \nsS_d \setminus \{x\}$ is close to $x$, then $K\cdot(x' - x)$
is close to a space contained in $T_x\nsS_d$. In particular and more precisely,
there exists a ball $B' \subseteq K^n$ containing $x$ such that for any
$x' \in B' \cap \nsS_d \setminus \{x\}$, we have $\res(K\cdot(x' - x)) \subseteq \res(T_x\nsS_d)$.
After possibly further shrinking $B'$, $(\nsS_i)_i$ becomes $d$-translatable on $B'$
and then, any one-dimensional subspace of $\tsp_{B'}((\nsS_i)_i)$ is of the form
$\res(K\cdot(x' - x))$ for some $x' \in B' \cap \nsS_d \setminus \{x\}$. For dimension reasons,
this implies
\[
\res(T_x\nsS_d) = \tsp_{B'}((\nsS_i)_i).
\]

Now consider Whitney's Condition~(a), i.e., suppose we are given a point $u \in S_d$ and a sequence $v_\mu \in S_j$ ($j > d$) as in
Definition~\ref{defn:whit} (\ref{it:wai}). Set $B := \res\1(u)$ and let $B' \subseteq B \cap \nsS_{\ge j}$ be a ball containing $[v_\mu]$. Then
\[
\lim_\mu T_{v_\mu}S_j = \res([T_{v_\mu}S_j]) = \res(T_{[v_\mu]}\nsS_j) = \tsp_{B'}((\nsS_i)_i) \supseteq
\tsp_{B}((\nsS_i)_i) = T_uS_d
.
\]
For Whitney's Condition~(b), suppose we are given $u \in S_d$ and sequences $u_\mu \in S_d$ and $v_\mu \in S_j$ ($j > d$) as in
Definition~\ref{defn:whit} (\ref{it:wbi}). Again set $B := \res\1(u)$. Since $[u_\mu], [v_\mu] \in B \subseteq \nsS_{\ge d}$,
we can apply Corollary~\ref{cor:whit-b} to $[u_\mu] \in \nsS_d$ and $[v_\mu] \in \nsS_j$ and obtain
a finite, $\LHen(\code{B})$-definable set $M \subseteq \Gamma$ such that $\val([u_\mu] - [v_\mu]) \notin M$ implies
$\dir([u_\mu] - [v_\mu]) \in \tsp_{B'}((\nsS_i)_i)$ for a ball $B' \subseteq \nsS_{\ge j}$ containing $[v_\mu]$.
In particular, $M$ is $\LHen(u)$-definable (viewing $u$ as an element of the residue field),
so $M \subseteq \{0\}$ by Remark~\ref{rem:dclGamma} and thus indeed $\val([u_\mu] - [v_\mu]) \notin M$.
Therefore we obtain
\begin{align*}
\lim_\mu k\cdot (u_\mu - v_\mu) &
= \res([k \cdot (u_\mu - v_\mu)]) = \res(K \cdot ([u_\mu] - [v_\mu]))\\
&\subseteq \tsp_{B'}((\nsS_i)_i)
= \res(T_{[v_\mu]}\nsS_j) = \res([T_{v_\mu}S_j]) = \lim_\mu T_{v_\mu}S_j
,
\end{align*}
which finishes the proof.
\end{proof}

Using Proposition~\ref{prop:t->whit}, it is now easy to deduce that t-stratifications
``are'' also classical Whitney stratifications. To be consistent with Subsection~\ref{subsect:alg},
we fix a Noetherian integral domain $A$ of characteristic $0$, we set $\Lx := \LHen(A)$, and we let $\Tx$ be the theory of
Henselian valued fields $K$ of equi-characteristic $0$ with ring homomorphism $A \to K$.

\begin{thm}\label{thm:t->whit}
Let $A$, $\Lx$, and $\Tx$ be as defined right above. Suppose that $\phi_\nu$ ($\nu = 1, \dots, \ell$) and $\psi_i$
($i = 0, \dots, n$) are $\Lring(A)$-formulas in $n$ free variables such that
for any model $K \models \Tx$, $(\psi_i(K))_i$ is a t-stratification of $K^n$
reflecting $(\phi_\nu(K))_\nu$. Suppose moreover that the formulas
$\psi_i$ are quantifier free.
Then for both $k = \bbR$ and $k = \bbC$ and for any ring homomorphism $A \to k$, we have
the following, where $X_\nu := \phi_\nu(k)$ and $S_i := \psi_i(k)$.
\begin{enumerate}
\item
$(S_i)_i$ is a Whitney stratification of $k^n$ (see Definition~\ref{defn:whit})
\item
Each $X_\nu$ is a union of some of the connected components of the
sets $S_i$ (in the analytic topology).
\end{enumerate}
In particular, each $\psi_i$ is an algebraic variety which is smooth over the fraction field of $A$.
\end{thm}

Note that by Corollary~\ref{cor:alg}, for any $(\phi_\nu)_\nu$ as above we can find
$(\psi_i)_i$ defining a t-stratification as above, so indeed we obtain a new proof of the existence of Whitney
stratifications (for $\Lring(A)$-definable subsets of $\bbR^n$ or $\bbC^n$).

\begin{proof}[Proof of Theorem~\ref{thm:t->whit}]
Let $K$ be the non-standard model of $k$ used in Proposition~\ref{prop:t->whit};
we consider it as an $\Lx$-structure using the ring homomorphism $A \to k \inject K$.
Then the conclusion of Proposition~\ref{prop:t->whit} is that $(S_i)_i$ is a
$C^1$-Whitney stratification. To finish the proof of (1), we have to get rid of this
``$C^1$''. By taking $k = \bbC$, we obtain that each $\psi_i(\bbC)$ is an algebraic manifold;
since $\psi_i$ is quantifier free, it can be viewed as a variety which is smooth over $\bbC$;
in particular $\psi_i(\bbR)$ is a $C^\infty$-sub-manifold of $\bbR^n$.

It remains
to verify (2). We have to show that for each $d \le n$ and each $\nu \le \ell$, both
$S_d \cap X_\nu$ and $S_d \setminus X_\nu$ are open in $S_d$.
Since this is first order, we can instead prove the corresponding statement in $K^n$,
i.e., that $\nsS_d \cap \nsX_\nu$ and $\nsS_d \setminus \nsX_\nu$ are open in $\nsS_d$.

Let $x \in \nsS_d$ be given. We choose a ball $B \subseteq \nsS_{\ge d}$ containing $x$,
we choose an exhibition $\pi\colon B \to K^d$
of $V := \tsp_B((\nsS_i)_i)$, and we shrink $B$ such that each $\pi$-fiber intersects
$B \cap \nsS_d$ in a single point. Then $V$-translatability implies that
the set $B \cap \nsS_d$ is either disjoint from $\nsX_\nu$
or entirely contained in $\nsX_\nu$.
\end{proof}

\section{Sets up to isometry in \texorpdfstring{$\bbQ_p$}{\041\032p}}
\label{sect:QpB}

The main conjecture of \cite{i.QpB} essentially is a classification of
definable subsets of $\bbZ_p^n$ up to isometry. More precisely, it
classifies certain trees associated to definable sets, which are closely related
to isometry types.
The original motivation for the present article was to prove that conjecture
for $p$ sufficiently big. This is achieved with Theorem~\ref{thm:QpB}.
Let us recall the trees considered in \cite{i.QpB}.

\begin{defn}
For a set $X \subseteq \bbZ_p^n$, the \emph{tree} $\Tr(X)$ associated to $X$
is the partially ordered set of those balls $B \subseteq \bbZ_p^n$ which intersect $X$
non-trivially; the ordering is given by inclusion. (If $X \ne \emptyset$, then $\Tr(X)$
is indeed a rooted tree, with root $\bbZ_p^n$.)
\end{defn}

It is not difficult to check that for topologically closed sets $X$,
a tree encodes exactly the isometry type of $X$. In general, we
have the following (see \cite[Lemma~3.1]{i.QpB}).

\begin{lem}\label{lem:TrIso}
Let $X, X' \subseteq \bbZ_p^n$ be any sets and write $\bar X$, $\bar X'$ for their
topological closures. Then
there is a natural bijection between the set of isometries $\bar X \to \bar X'$
and the set of isomorphisms of partially ordered sets $\Tr(X) \to \Tr(X')$.
\end{lem}

If $X \subseteq \bbZ_p^n$ is a definable
set of dimension $d$, then according to \cite[Conjecture~1.1]{i.QpB}
$\Tr(X)$ should be a ``tree of level $d$''. We will
more or less recall this definition below, but instead of speaking about a
tree $T$ itself, we will formulate it in terms of sets $X$ with $\Tr(X) \cong T$.
More precisely, we will introduce the notion of a subset $X \subseteq \bbZ_p^n$ being
``of level $\le d$''; this will be slightly stronger than
$\Tr(X)$ being of level $d$. Our main result will then be that in sufficiently
big residue characteristic, every definable
set of dimension $\le d$ is also a set of level $\le d$.
(Note that instead of calling the corresponding trees ``of level $d$'', they should better
also have been called ``of level $\le d$''. This better terminology is used in \cite{i.treesICMS},
and we will also use it below.)

The differences between $\Tr(X)$ being of level $\le d$ and $X$ being of level $\le d$ are the following.
\begin{itemize}
\item Trees classify the topological closures of definable sets
up to isometry. The notion of a set of level $\le d$ captures the isometry type of the set itself.
\item Some definable sets have a more complicated 
isometry type in small residue characteristic. Since our present result
only speaks about sufficiently big residue characteristic, we omit these from the notion
of sets of level $\le d$.
\end{itemize}

We also take the opportunity to generalize the conjecture from \cite{i.QpB} as follows.
\begin{itemize}
\item In \cite{i.QpB}, only subsets of $\bbZ_p^n$ are considered. we allow subsets of $\bbQ_p^n$.
\item Instead of working only in $\bbQ_p$, we work in any henselian valued field whose
residue field is finite and has sufficiently big characteristic. In particular, this includes
finite extensions of $\bbQ_p$ and function fields $\bbF_{p^r}((t))$. Non-discrete value groups
are also allowed.
\end{itemize}

Here is the precise formulation of the main result of this section.

\begin{thm}\label{thm:QpB}
Suppose that $\phi(x, y)$ is an $\LHen$-formula (see Definition~\ref{defn:LHen}), where $x$ is a tuple of valued field variables
and $y$ is a tuple of arbitrary variables. Then there exists an $N \in \bbN$ with the following property.
If $K$ is a Henselian valued field whose residue
field is finite and has characteristic at least $N$ and if moreover $b$ is any tuple in $K$
of the same sort as $y$, then $X := \phi(K, b)$
is a set of level $\le \dim X$ in the sense of Definition~\ref{defn:lev}.
\end{thm}

For this to (almost) classify definable subsets of $\bbZ_p^n$ up to isometry, one also
needs a converse. Indeed, by \cite[Theorem~1.2]{i.QpB}, for any tree $T$ of level $\le d$, there exists a
definable set $X \subseteq \bbZ_p^n$ of dimension $\le d$ with $\Tr(X) \cong T$.

Now we introduce the notion of sets of level $\le d$.
The translation between this
and trees of level $\le d$ in the sense of \cite[Definition~4.1]{i.QpB} is pretty straight forward; a part
of this translation is written down in detail in the claim below Definition~4.3 of \cite{i.treesICMS}.
(The precise relation is: if $X$ is a subset of $\bbZ_p^n$ and it is of level $\le d$,
then $\Tr(X)$ is of level $\le d$.)

\begin{notn}
We write $\Loag = \{0, +, -, \le\}$ for the language of ordered abelian groups;
by ``$\Loagpar$-definable'', we mean definable in the language $\Loag$, where parameters are allowed.
\end{notn}

\begin{defn}\label{defn:lev}
Suppose that $K$ is a valued field with finite residue field.
A subset $X \subseteq K^n$ is a \emph{set of level $\le d$} if it can be obtained as follows.

Choose any $m \in \bbN$, any $s_1, \dots, s_m \in K^n$, and
set $S_0 := \{s_1, \dots, s_m\}$.
In the case $d = 0$, we (only) require $X \subseteq S_0$.

In the case $d \ge 1$, additionally choose, for each $\ell \le m$ and each $\lambda \in \Gamma$,
an enumeration $B_{\ell,1,\lambda},\dots, B_{\ell, |k|^n,\lambda}$ of the maximal strict subballs
of $\ball{s_\ell, \ge \lambda}$ (i.e., $\rado(B_{\ell,j,\lambda}) = \lambda$ and $\ball{s_\ell, \ge \lambda} = \bigdcup_j B_{\ell,j,\lambda}$).
Finally choose, for each $\ell \le m,j \le |k|^n,\lambda \in \Gamma$, a set $Y_{\ell,j,\lambda} \subseteq K^{n-1}$ of level $\le d-1$. We reqiure the following.
\begin{enumerate}
\item For each $\ell,j,\lambda$ as above, if $B_{\ell,j,\lambda} \cap S_0 = \emptyset$, then $X \cap B_{\ell,j,\lambda}$ is isometric to
$Y_{\ell,j,\lambda} \times \ball{0,>\lambda}$, where by $\ball{0,>\lambda}$ we mean a one-dimensional ball.
\item For each fixed $\ell,j$, the family $(Y_{\ell,j,\lambda})_{\lambda \in \Gamma}$ is of level $\le d-1$
uniformly in $\lambda$, in the sense described below.
\end{enumerate}
A family of sets $X_\kappa \subseteq K^n$ parametrized some $\kappa \in M \subseteq \Gamma^\nu$ is of level $\le d$
\emph{uniformly in $\kappa$} (where $M$ is $\Loagpar$-definable), if all of the above choices can be made for all $\kappa$
such that moreover the following holds. (We use the notation from above, without indices $\kappa$.)
\begin{enumerate}\setcounter{enumi}{2}
\item For each $m_0 \in \bbN$ and each $I_0 \subseteq \{1, \dots, m_0\}$, the
set $M'\subseteq M$ of those $\kappa$ such that $m = m_0$ and $X \cap S_0 = \{s_\ell \mid \ell \in I_0\}$ is $\Loagpar$-definable.
\item For each $M'$ as in (3) and each $\ell,\ell' \le m_0$, $\val(s_\ell - s_{\ell'})$ is an $\Loagpar$-definable function of $\kappa \in M'$.
\item For each $M'$ as in (3) and each $\ell \le m_0,j \le |k|^n$, the sets $Y_{\ell,j,\lambda}$ are of level $\le d-1$ uniformly in
$(\lambda,\kappa) \in \Gamma \times M'$ (and not just uniformly in $\lambda$).
\end{enumerate}
\end{defn}

Note that whether a set is of level $\le d$ only depends on the isometry type of the set.

\begin{proof}[Proof of Theorem~\ref{thm:QpB}]
Let $\Lx$ be a language consisting of $\LHen$, any set $C$ of constant symbols in any sorts,
an angular component map $\aci\colon K \to k$,
and Skolem functions inside $k$, and let $\Tx$ be the corresponding expansion of $\THen$.
Elimination of valued field quantifiers implies that every $\Lx$-formula $\psi(z)$,
where $z$ is a tuple of $\Gamma$-variables, is already equivalent (modulo $\Tx$) to an
$\Loag(C')$-formula for some suitable set of constants $C'$ (namely, $C' = \dcl_{\Lx}(\emptyset) \cap \Gamma$).

For $N \in \bbN$, let $\CN$ be the class of all 
henselian valued fields with finite residue field of characteristic $> N$,
considered as $\Lx$-structures.
(Note that by \cite[Corollary~1.6]{Pas.ac}, valued fields with finite residue field always admit an angular component map.)

In the following, we will work uniformly in all models $K \models \Tx$. In particular, unless specified otherwise,
by a ``definable set $X$ in $K$'', we mean an $\Lx$-formula $X$ and, abusing notation, we write $X$ instead of $X(K)$.
We will prove the following by induction.

\medskip

\textbf{Claim.} Suppose that we have $n,d,\nu \in \bbN$ and $S_{i,q}$, $X_q$, $Q$, $\chi$, $M$ such that
for every model $K \models \Tx$, the following holds:
\begin{itemize}
\item
$Q$ is a $\emptyset$-definable set and $(S_{i,q})_{q \in Q}$, $(X_q)_{q \in Q}$ are $\emptyset$-definable
families such that for every $q \in Q$, $\dim X_q \le d$ and
$(S_{i,q})_i$ is a t-stratification of a subball of $K^n$
reflecting $X_q$;
\item
$\chi\colon Q \surject M \subseteq \Gamma^\nu$ is a $\emptyset$-definable map such that
for every $q,q' \in Q$ with $\chi(q) = \chi(q')$, there exists a risometry $\alpha_{q,q'}\colon K^n \to K^n$ which
sends $((S_{i,q})_i, X_q)$ to $((S_{i,q'})_i, X_{q'})$ and which is definable (with parameters) separately for each $q, q' \in Q$ and $K \models \Tx$.
\end{itemize}
Then there exists an $N \in \bbN$ such that in every $K \in \CN$, we have the following.
For every (not necessarily definable) cross section $M \to Q, \kappa \mapsto q_{\kappa} \in \chi\1(\kappa)$,
the family $(X_{q_{\kappa}})_{\kappa \in M}$ is uniformly of level $\le d$.

\medskip

The claim implies the theorem using a singleton for $Q$. Indeed,
set $X := \phi(K,c)$, where $\phi(x,y)$ is the formula given in the theorem and $c$ is a tuple of
constants. Then, working in the language $\LHen \cup \{c\}$,
we can apply Corollary~\ref{cor:pos-char} to uniformly obtain t-stratifications $(S_i)_i$ reflecting $X$
in each $K \models \THen$. After that, we enlarge the language to $\Lx$ and apply the claim.

\medskip

\emph{Proof of the Claim.}
First note that we may suppose that $m := |S_{0,q}|$ is constant (for all models $K \models \Tx$ and all $q \in Q$).
Indeed, we can partition $Q$ according to the cardinality of $S_{0,q}$, and
the existence of the risometries $\alpha_{q,q'}$ implies that this partition
induces a partition of $M$. The partition of $M$ can be defined by $\Loag(C')$-formulas, so that it
also induces a finite, $\Loagpar$-definable partition of $M$ in any $K \in \CN$ for $N \gg 1$;
the notion of being uniformly of level $\le d$ is not affected by such a partition.

Next, we choose an enumeration $s_{1,q}, \dots, s_{m,q}$ of $S_{0,q}$
which is definable uniformly in $q$ (and uniformly in all models $K \models T$) and which satisfies
\[\tag{$\square$}
s_{\ell,q'} = \alpha_{q,q'}(s_{\ell,q})
\]
whenever $\chi(q) = \chi(q')$. To see that such an enumeration exists, first note that
since we have Skolem functions for finite subsets of $\Gamma$ and $k$, we also have
Skolem functions for finite subsets of $\RV^{(n)}$ (the angular component map yields a definable bijection $\RV^{(n)} \setminus \{0\} \to (k^n \setminus \{0\}) \times \Gamma$).
By Lemma~\ref{lem:fin&iso}, the map $\rho_q\colon x \mapsto \rv(x - S_{0,q})$ is injective on $S_{0,q}$ and
the Skolem functions can be used to enumerate the image $\im(\rho_q)$ in a way depending definably on $\code{\im(\rho_q)}$.
Then we automatically obtain ($\square$), since
we have $\rho_{q'} \circ \alpha = \rho_q$ for any risometry $\alpha\colon S_{0,q} \to S_{0,q'}$.

For $N \gg 1$, $(s_{\ell,q})_\ell$ is an enumeration of $S_{0,q}$ in any $K \in \CN$.
With this enumeration, Definitions~\ref{defn:lev}~(3) and (4) are satisfied. Indeed, ($\square$) implies that
the set $I_q := \{\ell \le m \mid s_{\ell,q} \in X_{q}\}$ only depends on $\chi(q)$, so fixing $I_q$ yields
a definable subset $M' \subseteq M$, which then yields the required $\Loagpar$-definability of $M'$ in any $K \in \CN$
for $N\gg0$ (independently of the cross section $\kappa \mapsto q_\kappa$).
Similarly, for each $\ell, \ell' \le m$, ($\square$) implies that $\val(s_{\ell,q} - s_{\ell',q})$
only depends on $\chi(q)$, which yields (4).

If $d = 0$, then $\dim X_q = 0$ implies $X_q \subseteq S_{0,q}$ for all $q$ and we are already done, so now assume $d \ge 1$.
Also fix $\ell \le m$ for the remainder of the proof.

For every $\lambda \in \Gamma$ and $q \in Q$, let $(B_{u,\lambda,q})_{u \in k^n}$ be
the family of maximal strict subballs of $\ball{s_{\ell,q}, \ge \lambda}$.
We assume that for $q, q' \in Q$ with $\chi(q) = \chi(q')$, we have $\alpha_{q,q'}(B_{u,\lambda,q}) = B_{u,\lambda,q'}$
and (using the map $\aci$) that $B_{u,\lambda,q}$ is definable uniformly in $u$, $\lambda$, and $q$.

Now also fix $u, \lambda, q$ for the moment, set $B := B_{u,\lambda,q}$, and suppose that $B \cap S_{0,q} = \emptyset$.
Choose a one-dimensional subspace $V \subseteq \tsp_B((S_{i,q})_i)$ and an exhibition $\pi\colon K^n \surject K$ of $V$.
Using that straighteners are definable uniformly for all $K \models \Tx$ (Corollary~\ref{cor:unif-straight}), we also
have a straightener in $K \in \CN$ for $N \gg 1$, which yields an isometry $B \cap X_q \to (F \cap X_q) \times \pi(B)$ for any fiber $F = \pi\1(x) \cap B, x \in \pi(B)$. Thus (1) holds and it remains to verify (5) with $Y_{\ell,j,\lambda} = F \cap X_q$ (which is a strengthening of (2)).

The exhibition $\pi$ can be chosen uniformly in $u, \lambda, q$, so we may as well assume that it
does not depend on $u, \lambda, q$ (only to simplify notation).
In each $K \models \Tx$, the set $F \cap X_q$ is reflected by the t-stratification $(F \cap S_{i+1,q})_{i \le n-1}$.
We set $\Lx' := \Lx \cup \{u\}$; our plan is to
apply the induction hypothesis in $\Lx'$ to $X'_{\lambda, q, x} := F \cap X_q$ and
$S'_{i, \lambda, q, x} := F \cap S_{i+1,q}$, considered as families parametrized by
$Q' := \{(\lambda, q, x) \mid \lambda \in \Gamma, q \in Q, x \in \pi(B_{u,\lambda,q})\}$. 
To this end, we first modify the definitions of $X'_{\lambda, q, x}$ and
$S'_{i, \lambda, q, x}$ in the cases where $B \cap S_{0,q} \ne \emptyset$
(to ensure that $(S'_{i, \lambda, q, x})_i$ always reflects $X'_{\lambda, q, x}$),
e.g.\ by setting $X'_{\lambda, q, x} := \emptyset$ and $S'_{n-1, \lambda, q, x} := F$ in these cases.
Next, we define
$\chi'\colon Q' \surject M' := \Gamma \times M, (\lambda, q, x) \mapsto (\lambda, \chi(q))$;
it remains to find suitable risometries
\[
\alpha'_{(\lambda,q,x),(\lambda,q',x')}: \pi\1(x) \cap B_{u,\lambda,q} \to \pi\1(x') \cap B_{u,\lambda,q'}
\]
for $\lambda \in \Gamma$, $q, q' \in Q$ with $\chi(q) = \chi(q')$,
$x \in \pi(B_{u,\lambda,q})$, and $x' \in \pi(B_{u,\lambda,q'})$. These are obtained
by applying Lemma~\ref{lem:tr-faser} (2) to $\alpha_{q,q'}(B_{u,\lambda,q}) = B_{u,\lambda,q'}$.

Now the conclusion of the induction
hypothesis yields exactly Definition~\ref{defn:lev}~(5). Indeed,
suppose that $K \in \CN$ is given (for $N \gg 1$) and that in $K$, we have a cross section $M \to Q, \kappa \mapsto q_{\kappa}$.
Fix $u \in k^n$; this turns $K$ into an $\Lx'$-structure. (Note that our $u$ here corresponds to the $j$ in Definition~\ref{defn:lev}.)
We choose a cross section $M' \to Q'$ of the form $(\lambda, \kappa) \mapsto (\lambda, q_{\kappa}, x_{\lambda,\kappa})$ with $x_{\lambda,\kappa} \in \pi(B_{u,\lambda,q_{\kappa}})$ arbitrary. By induction, $X'_{\lambda, q_\kappa, x_{\lambda,\kappa}}$
is of level $\le d-1$ uniformly in $\lambda$ and $\kappa$, which is what we had to prove.
\end{proof}

In \cite[Section~7]{i.QpB}, several strengthenings of the conjecture about the trees have been proposed;
to conclude this section, let me comment on these strengthenings.
\begin{itemize}
\item
The tree $\Tr(\Zp^n)$ can be considered as an imaginary sort; then, for any definable
$X \subseteq \Zp^n$, the tree $\Tr(X)$ is a definable subsets of $\Tr(\Zp^n)$.
Conjecture~7.1 of \cite{i.QpB} describes arbitrary definable subsets
$Y \subseteq \Tr(\Zp^n)$ instead of only those of the form $\Tr(X)$.
For big $p$, it should also be possible to prove that conjecture, using a t-stratification
reflecting the map
\[
\chi(x) := \code{\{\gamma \in \Gamma \mid \ball{x, \ge \gamma} \in Y \}}
.
\]
\item
In \cite[Section~7.2]{i.QpB}, a version of the conjecture has been proposed for
arbitrary Henselian valued fields of characteristic $(0,0)$ (without giving much details).
However, as noted there, the conjecture has far less meaning when the residue field
is infinite (since then, too many isometries exist), so instead of considering
pure abstract trees, one should consider trees with some additional residue field data.
Driving this idea further is what finally led to the definition of t-stratifications.

According to \cite{i.QpB}, these ``trees in characteristic $(0,0)$'' should imply
the conjecture in $\bbQ_p$ for big $p$ and should even yield some kind
of uniformly in $p$. Our proof of Theorem~\ref{thm:QpB} indeed yields uniformity
in the following sense. Given a formula $\phi(x,y)$ as in the theorem,
the $\Loagpar$-definable objects that we construct to prove that $\phi(K,b)$ is of level $\le d$
(for Definition~\ref{defn:lev}~(3), (4)) can be defined by $\Lx$-formulas
not depending on $K$ and $b$ (but taking $b$ as a parameter).
Moreover, since in this uniform
setting, the cardinality of the residue field grows, it becomes worthwile
to note that our proof moreover yields that in
Definition~\ref{defn:lev}~(5), $Y_{\ell,j,\lambda}$ is uniform also in $j$
and not just in $\lambda$ and $\kappa$.
\end{itemize}

\section{Open questions}
\label{sect:open}

There are several ways in which it might be possible to enhance the results of this article.

\subsection{A stronger notion of t-stratification}

Recall from the introduction that t-stratifications to not satisfy the straightforward translation of
Whitney's Condition~(a):
for two strata $S_d, S_j$ with $d < j$, $x \in S_d$, and $y \in S_j$,
we have that ``$T_yS_j$ is close to containing $T_xS_d$ when $y$ is close to $x$'', whereas Condition~(a)
requires that $T_yS_j$ converges to a space containing $T_xS_d$.
Also recall that in the $p$-adics, the existence of Whitney stratifications in this more classical sense has been proven in \cite{CCL.cones}. It seems plausible that there exists a common generalization of both kinds of stratifications (at least in equi-characteristic $0$ and if the value group is of rank one). Such a generalization might be defined as follows.

Let us define ``stronger
risometries'': maps $\phi$ such that $\val\big((\phi(x) - \phi(x')) - (x - x')\big) > \val(x - x') + \delta$
for some given $\delta \ge 0$. (For $\delta = 0$, this is just a usual risometry.) This
yields corresponding notions of translatability which we call ``$\delta$-strong translatability''.

Using this, a ``strong t-stratification'' should roughly require that
for any $\delta \ge 0$ and any ball $B$ ``sufficiently far away from $S_{\le d-1}$'',
we have $\delta$-strong $d$-translatability on $B$.
More precisely, it seems plausible that we can obtain $\delta$-strong $d$-translatability on any
ball $B$ which is contained
in a ball $B' \subseteq S_{\ge d}$ with $\rado B \ge \rado B' + \delta$.

Note that this indeed implies Condition~(a). For any $x \in S_d$ and any $\delta \ge 0$, there exists a ball $B$ around $x$ which is sufficiently far away from $S_{\le d-1}$ in the above sense, and $\delta$-strong translatability on $B$ then implies that for $y \in B \cap S_j$, $T_yS_j$ is $\delta$-close to a space containing  $T_xS_d$.

\subsection{Mixed characteristic}
\label{subsect:mixed}

It should be possible to prove the existence of a variant of t-stratifications
in mixed characteristic, but again, it is not entirely clear how this
variant has to be formulated. For a ball $B \subseteq S_{\ge d}$,
even $0$-strong (i.e., usual) $d$-translatability can only be expected on subballs $B'$ of $B$
with $\rado B' \ge \rado B + \delta$ for some fixed $\delta$ (depending only on the
t-stratification). This can be seen, for example, at the cusp curve in characteristic 2
(see \cite[Section~3.3]{i.QpB} or \cite[Section~5.4]{i.treesICMS} for a detailed computation).
In terms of the description of the trees of \cite{i.QpB}, this $\delta$ would
be exactly the maximal length of the finite trees appearing at the beginning
of side branches.

When the valuation of the residue characteristic $p$ is finite
(i.e., when there are only finitely elements of $\Gamma$ between $0$ and $\vali(p)$), then
in the previous paragraph, it should be also possible to require $\delta$ to be finite,
and the resulting notion of t-stratification might be the ``right one''.
However, if $\vali(p)$ is not finite, we are forced to allow finite multiples of $\vali(p)$ for $\delta$.
But then, I am afraid that then the notion of t-stratification
becomes too weak e.g.\ to imply Proposition~\ref{prop:sak}; in
particular, the induction in the proof of Theorem~\ref{thm:main} would fail.

\subsection{Getting classical Whitney stratifications more generally}

The fact that the existence of t-stratifications implies the existence of
Whitney stratifications should also work in languages other than the pure
(semi-)algebraic one. For this to work,
we need the existence of t-stratifications $(S_i)_i$ which are defined
without using the valuation. Probably Proposition~\ref{prop:unif-cl} can
be applied to prove such a result, but I did not check it. In the algebraic
language, we used this to deduce a posteriori that each $S_i$ is smooth.
This too, should work more generally, again with an argument
that manifolds in $\bbC^n$ which are $C^1$ in the sense of complex differentiation
are automatically smooth.

\subsection{Minimal t-stratifications}

It would be nice if,
for every definable set $X \subseteq K^n$, there would be a ``minimal'' t-stratification
$(S_i)_i$ reflecting $X$. ``Minimal'' could mean that for any other t-stratification
$(S'_i)_i$ reflecting $X$, we have $S_{\le i} \subseteq S'_{\le i}$ for all $i$.
Moreover (or alternatively), one might hope that for a minimal $(S_i)_i$,
a definable risometry $K^n \to K^n$ preserves $X$ if and only if it preserves $(S_i)_i$
(in general, there are less risometries preserving $(S_i)_i$).
In the case of Whitney stratifications of complex analytic spaces,
minimal stratifications in the first sense have indeed been constructed by Teissier; see \cite{Tei.singICM}.

There are (at least) two reasons for minimal t-stratifications not to exist, but for both
of them, all hope is not lost.
The first obstacle is the non-canonicity of t-stratifications explained in Example~\ref{ex:alm-sing}. This might be overcome
as follows. Instead of letting $S_{\le i}$ be a subset of $K^n$, we let
it be a subset of the set of subballs of $K^n$, where points are also considered as balls.
Then we require $d$-translatability on a ball $B \subseteq K^n$ iff no ball of $S_{\le d - 1}$ is (strictly?) contained in $B$. At least for Examples~\ref{ex:alm-sing} and \ref{ex:ball}, this seems to solve the problem.

A second problem is that one can construct a set $X$ such that whether $X$ is $d$-translatable
on some ball $B$ does not depend definably on $B$ (see Example~\ref{ex:tsp-ndef}). Since for t-stratifications, $d$-translatability
is always definable (Lemma~\ref{lem:trans-def}), any t-stratification reflecting $X$ will necessarily
have less risometries  than $X$ preserving it.
I do not think that it is possible to solve this problem in general, but it might
be possible to find a good condition on the residue field which avoids the problem.
A candidate which at least destroys Example~\ref{ex:tsp-ndef} is the following.
For any definable function $f\colon k^n \to k$, there exists a definable function
$\tilde{f}\colon \valring^n \to \valring$ such that $\res \circ \tilde{f} = f \circ \res$.


\begin{thebibliography}{10}

\bibitem{BCR.realGeom}
{\sc J.~Bochnak, M.~Coste, and M.-F. Roy}, {\em Real algebraic geometry},
  vol.~36 of Ergebnisse der Mathematik und ihrer Grenzgebiete (3) [Results in
  Mathematics and Related Areas (3)], Springer-Verlag, Berlin, 1998.
\newblock Translated from the 1987 French original, Revised by the authors.

\bibitem{CCL.cones}
{\sc R.~Cluckers, G.~Comte, and F.~Loeser}, {\em Local metric properties and
  regular stratifications of {$p$}-adic definable sets}, Comment. Math. Helv.,
  87 (2012), pp.~963--1009.

\bibitem{CL.analyt}
{\sc R.~Cluckers and L.~Lipshitz}, {\em {Fields with analytic structure.}}, J.
  Eur. Math. Soc. (JEMS), 13 (2011), pp.~1147--1223.

\bibitem{CL.bmin}
{\sc R.~Cluckers and F.~Loeser}, {\em b-minimality}, J. Math. Log., 7 (2007),
  pp.~195--227.

\bibitem{CL.mot}
\leavevmode\vrule height 2pt depth -1.6pt width 23pt, {\em Constructible
  motivic functions and motivic integration}, Invent. Math., 173 (2008),
  pp.~23--121.

\bibitem{Den.rat}
{\sc J.~Denef}, {\em The rationality of the {P}oincar\'e series associated to
  the {$p$}-adic points on a variety}, Invent. Math., 77 (1984), pp.~1--23.

\bibitem{iF.dim}
{\sc A.~Fornasiero and I.~Halupczok}, {\em Dimension in topological structures:
  topological closure and local property}, in Groups and Model Theory, vol.~576
  of Contemp. Math., Amer. Math. Soc., Providence, RI, 2012, pp.~89--94.

\bibitem{i.QpB}
{\sc I.~Halupczok}, {\em Trees of definable sets over the {$p$}-adics}, J.
  Reine Angew. Math., 642 (2010), pp.~157--196.

\bibitem{i.treesICMS}
\leavevmode\vrule height 2pt depth -1.6pt width 23pt, {\em Trees of definable
  sets in $\mathbb{Z}_p$}, in Proceedings of the conference ``{M}otivic
  {I}ntegration and its interaction with Model Theory and Non-{A}rchimedean
  Geometry'', R.~Cluckers, J.~Nicaise, and J.~Sebag, eds., Cambridge University
  Press, 2011, pp.~87--107.

\bibitem{Pas.ac}
{\sc J.~Pas}, {\em On the angular component map modulo {$P$}}, J. Symbolic
  Logic, 55 (1990), pp.~1125--1129.

\bibitem{Tei.singICM}
{\sc B.~Teissier}, {\em Sur la classification des singularit\'es des espaces
  analytiques complexes}, in Proceedings of the {I}nternational {C}ongress of
  {M}athematicians, {V}ol.\ 1, 2 ({W}arsaw, 1983), Warsaw, 1984, PWN,
  pp.~763--781.

\bibitem{Dri.dimDef}
{\sc L.~van~den Dries}, {\em Dimension of definable sets, algebraic boundedness
  and {H}enselian fields}, Ann. Pure Appl. Logic, 45 (1989), pp.~189--209.
\newblock Stability in model theory, II (Trento, 1987).

\bibitem{Whi.strat}
{\sc H.~Whitney}, {\em Tangents to an analytic variety}, Ann. of Math. (2), 81
  (1965), pp.~496--549.

\end{thebibliography}
\end{document}